 \documentclass[11pt]{article}   
\usepackage{amssymb,amscd,latexsym}   
\usepackage{amsmath}
\usepackage{amsthm}

\newtheorem{Theorem}{Theorem}
\newtheorem{Lemma}[Theorem]{Lemma}
\newtheorem{Proposition}[Theorem]{Proposition}
\newtheorem{Corollary}[Theorem]{Corollary}

\newtheorem{Question}[Theorem]{Question}

\theoremstyle{definition}
\newtheorem{Definition}[Theorem]{Definition}
\newtheorem{Remark}[Theorem]{Remark}
\newtheorem{Example}[Theorem]{Example}

\numberwithin{Theorem}{section} \numberwithin{equation}{section}

\def\demo{\noindent{\it Proof. }}
\def\QED{\hfill$\Box$}
\def\qed{\QED}

\newcommand{\rar}{\rightarrow}
\newcommand{\lar}{\longrightarrow}

\newcommand{\xx}{{\bf x}}
\newcommand{\yy}{{\bf y}}
\newcommand{\pp}{{\mathbb P}}

\newcommand{\Sing}{\operatorname{Sing}}
\def\spec#1{{\rm Spec}\, (#1)}

\renewcommand{\O}{{\mathcal O}}

\newcommand{\F}{{\mathcal F}}

\renewcommand{\L}{{\mathcal L}}
\newcommand{\Proj}{{\mathbb P}}

\newcommand{\MM}{{\mathcal M}}

\newcommand{\G}{{\mathbb G}}
\newcommand{\NN}{{\mathcal N}}

\newcommand{\C}{{\mathbb C}}
\newcommand{\p}{{\mathbb P}}
\newcommand{\A}{{\mathbb A}}
\renewcommand{\a}{{\`a}}

\newcommand{\red}{\operatorname{red}}

\newcommand{\rk}{\operatorname{rk}}

\newcommand{\map}{\dasharrow}

\def\surjects{\twoheadrightarrow}

\title{\LARGE\sc  Homaloidal hypersurfaces and hypersurfaces with vanishing
Hessian}

\author{
{\normalsize\sc Ciro Ciliberto\thanks{This author thanks CNPq for support and the
Departamento de Matem\'atica at Recife for hospitality during the
preparation of this paper.}}
\vspace{-0.75mm}\\
{\small Dipartimento di Matematica}\vspace{-1.4mm} \\
{\small Universit\'a di Roma ``Tor Vergata''}\vspace{-1.4mm}\\
{\small 00133 Roma, Italy}\vspace{-1.4mm} \\
{\small e-mail: {cilibert@mat.uniroma2.it}} \\
\and {\normalsize\sc Francesco Russo\thanks{Partially
 supported by CNPq, Brazil and by PRONEX--FAPERJ.}}
\vspace{-0.75mm}\\
{\small Dipartimento di Matematica e Informatica}\vspace{-1.4mm} \\
{\small Universit\'a degli Studi di Catania}\vspace{-1.4mm}\\
{\small 95125 Catania, Italy}\vspace{-1.4mm} \\
{\small e-mail: {frusso@dmi.unict.it}}\\
 \and {\normalsize\sc Aron Simis\thanks {Partially supported by CNPq, Brazil and by PRONEX--FAPERJ.}}
\vspace{-0.75mm}\\
{\small Departamento de Matem\'atica}\vspace{-1.4mm} \\
{\small Universidade Federal de Pernambuco}\vspace{-1.4mm}\\
{\small 50740-540 Recife, PE, Brazil}\vspace{-1.4mm} \\
{\small e-mail: {aron@dmat.ufpe.br}} }

\begin{document}

\date{}
\maketitle


\begin{abstract}

We {introduce} various families of
irreducible homaloidal hypersurfaces in projective space $\p^ r$,
for all $r\geq 3$. Some of these are families of homaloidal
hypersurfaces whose degrees are arbitrarily large as compared to the
dimension of the ambient projective space. The existence of such a
family solves a question that has naturally arisen from the
consideration of the classes of homaloidal hypersurfaces known so
far. The result relies on a fine analysis
of hypersurfaces {that are dual} to certain scroll
surfaces. We also introduce an infinite family of determinantal
homaloidal hypersurfaces based   on a
certain degeneration of a generic Hankel matrix. {The latter
family} fit non--classical versions of de
Jonqui\`eres transformations. As a natural counterpoint, we broaden
up aspects of the theory of Gordan--Noether hypersurfaces with
vanishing Hessian determinant, bringing
over some more precision into the present knowledge.

\end{abstract}

\section*{Introduction}\label {intro}

The study of Cremona  transformations of
$\p^ r$ is {a}  classical and fascinating
subject(s) in algebraic geometry. The Cremona group of $\p^ r$ is
well understood only for $r\leq 2$. By contrast, in dimension $r\geq
3$ it is even problematic to produce non--trivial examples of
birational transformations of $\p^ r$.
Therefore, {any relevant addition to the universe of these
transformations} is very welcome, especially if it bridges {up}
with other interesting concepts in the field.

In this perspective, {a good example}
 is {that of a}
\emph{homaloidal hypersurface}. This is a projective hypersurface
$X\subset \p^r$, not necessarily reduced or irreducible,
defined by a homogeneous {polynomial
$f=f(x_0,\dots,x_r)$ of degree ${d\geq 2}$ whose}
partial derivatives define a Cremona
transformation of $\p^r$. {Quite
generally}, the rational map
$\phi_f:\p^r\dasharrow \p^r$ defined by the partial derivatives of
$f$ is called the \emph {polar map}  of the hypersurface $X$, so
that, if  $X$ is reduced, the indeterminacy locus of $\phi_f$ is
precisely the singular locus of $X$. For instance,  a smooth quadric
is homaloidal, inasmuch as its polar map is the usual polarity,
which is a projective transformation. However, if $X$ is smooth of
degree $d\geq 3$, then $X$ is never homaloidal, since its polar map
has no indeterminacy locus and it is defined by forms of degree
$d-1>1$.  Indeed, a
relevant role in the understanding of homaloidal hypersurfaces is
played by the analysis of their singular locus.

On the other hand, an obvious necessary
condition {in order that ${X}$ be
homaloidal, is the non--vanishing of the Hessian determinant $h(f)$
of $f$.  Note that, if one measures the \emph{complexity} of a hypersurface
by the degree of its polar map, the hypersurfaces with vanishing Hessian
have to be considered as the \emph{simplest} ones, and the
homaloidal hypersurfaces are the
{simplest} among those for which the Hessian
is not identically zero. Thus, { a couple of}
natural questions arise: what can we say about hypersurfaces with
identically vanishing Hessian? What are the relations, if any,
between these and homaloidal hypersurfaces?

As {is generally} known,  both
problems-- the} classification of homaloidal {
hypersurfaces} and {of hypersurfaces with vanishing Hessian}
-- play a classical role
in the history of algebraic geometry, perhaps with homaloidal
running first, as subsumed into Cremona theory,
while vanishing  Hessian
winning in drama ever since Gordan and Noether (see
\cite{GN}) showed that Hesse (see \cite{Hesse1}, \cite{Hesse2}) had
previously misapprehended the question.

{Although fairly understood, even the
theory of plane Cremona transformations} is already quite
involved. The early results of Noether (see, e.g., \cite[Remark 2.3
and ff.]{cremona}), inspired on Cremona's original work, showed how
much more complicated is the theory in $\p^3$. However, it has
perhaps been common  thought that, {
notwithstanding}  the difficulties of the general Cremona
theory, homaloidal hypersurfaces {would
be} easier to understand and eventually be subject
to classification. For instance, the classification  of
reduced homaloidal curves in $\p^2$ by Dolgachev (see
\cite{Dolgachev}){--} which shows that there are only three
types  up to projective transformations, and, more
generally, the examples coming from the theory of pre--homogeneous
vector spaces (see \cite{EKP}), whose degree is bounded in terms of
the embedding dimension {--} could have generated the
expectation that the degree of an irreducible, or perhaps even only
reduced, homaloidal hypersurface in $\p^r$ is at most $r+1$. If this
were proved to be the case, one would perhaps be half--way from the
classification goal.

Alas, nature had the upper hand.  Indeed, one of the main objectives
of this paper is to show that, as a counterpart to the planar case,
in which a full classification is  {
fairly} easy to state, the situation is
much more complicated in higher dimension. In fact, one of our main
results here is to show the existence of families of irreducible
hypersurfaces in $\p^ r$, for $r\geq 3$, with arbitrarily large
degree with respect to $r$ (see {Section} \ref {infser}).
{We think that this uncovers} some complex
phenomenology which makes  the
classification of irreducible homaloidal hypersurfaces quite
intricate and therefore deserves a deeper
scrutiny, { beyond} our presently {inadequate}
understanding of the matter. A special role is of course
played by the complicated nature of the scheme structure of the base
locus of a homaloidal Cremona transformation. In particular, for a
homaloidal hypersurface $X$, this is due to the existence of
embedded components originating some infinitely near base points for
the linear system of polars of $X$, which are somehow unexpected
inasmuch as they are not singular points of $X$ or do not
 even belong to $X$ (see, e.g.,
\cite{Aluffi} and  {Section} \ref{sing}). {Incidentally,
this phenomenon is already present}
in one of the plane cases appearing in Dolgachev's
classification.

As for the second question envisaged here,
the problem {is after all} to find the
homogeneous polynomial solutions $f$ of the classical \emph
{Monge--Amp\`ere} differential equation $h(f)=0$. It is therefore
not surprising to see how far {an
outpost} this question has reached in subsequent geometric
developments and how strong a role it has played in various other
areas, such as differential geometry and approximation theory (see,
for example, \cite{Segrediff}, \cite{FW} and \cite{control}).

{In their celebrated work  \cite{GN},
Gordan and Noether constructed counterexamples}
to Hesse's original claim to the effect that $X$ has
vanishing Hessian if and only if it is a
cone. The examples have been later revisited and partly
extended by several authors (in chronological order,
\cite{hessianoSegre}, \cite{Fr1},
\cite{Fr2}, \cite{Permutti1}, \cite{Permutti2}, \cite{Permutti3},
\cite{Lossen}). {In spite of the difficulty of their original
paper, the examples themselves} are not {all that} difficult to understand
 and  {can actually} be
easily described in explicit algebraic terms (see also
\cite{Perazzo}).

A second goal of this paper is to give a modern overview of the
known methods to deal with the problem of vanishing Hessian and to
generalize results of Permutti and Perazzo { quoted}
 above. One of the challenges is to
determine the structure of the dual variety to Gordan--Noether or
Permutti hypersurfaces, for which we add a tiny contribution that
may help improving our understanding of these defective dual
varieties. { As it turns,}  there is a strong relationship between {the}
 families of homaloidal
hypersurfaces {described here} and some hypersurfaces with
vanishing Hessian. {We hope to pursue work along this line in the near future}.

We now describe the sections of the paper in somewhat more detail.

The first  section contains a {recap}
of known concepts and is primarily {meant as a
collection of} properties of scroll surfaces and their
dual varieties {that are either} spread out
or difficult to find in the current
literature. The main results are contained in a series of
propositions (see Proposition~\ref{non_developable}
through Proposition~\ref{sing_of_dual}). We also describe the
behavior of more general rational scroll surfaces containing a
so--called \emph{line directrix}, and their dual hypersurfaces (see
Propositions~\ref{directrix} and \ref{multi_scrolls}). This section
prepares the ground for the more thorough considerations of the
third  section, { for which the present
material}  is essential
 in the construction of {the
announced} examples.

The second section starts with an overview of the aforementioned
{polar map} $\phi_f$ associated to a non--zero homogeneous
 polynomial $f$. After a
brief introduction about the polars and the Hessian of $f$, we
switch to the problem of the vanishing Hessian. Just enough of the
Gordan--Noether construction is reviewed in order to state a
geometric description of its structure (see
Proposition~\ref{gn_hypersurface}), based on a notion of {\it
core\/} of such a hypersurface. We next discuss  the work of
Permutti extending the previous construction in a special situation,
and following the same ideas we also give some features of
Permutti's generalized hypersurface (see Proposition~\ref{Pcarat}).
We proceed to establishing both the structure of the dual variety to
a Permutti hypersurface and of its polar image (see
Propositions~\ref{Pdual} and \ref{Pzeta}). The section ends with a
generalization of a result of Perazzo (see Proposition~\ref{Dhyper})
establishing a bound for the dimension of the image of $\phi_f$ for
a so-called $H$-hypersurface $X\subset \p^r$ with equation $f=0$,
i.e. a reduced hypersurface which contains a subspace of dimension
$t$ such that the general subspace of dimension $t+1$ through it
cuts out on $X$ a cone with a vertex  of dimension at least $r-t-1$.
The dual hypersurfaces to scrolls with a line directrix are special
cases of $H$--hypersurfaces and come up in our examples.

In the third section we { introduce}
 families of irreducible
homaloidal hypersurfaces, including the case {in which they
have} arbitrarily large degree as compared to the ambient
dimension. As a preliminary, we state a general principle for {
a} Cremona transformation saying that such a map always {\it
contracts\/} its Jacobian, and ask whether, in the case of a polar
map $\phi_f$, contraction is also sufficient for birationality,
provided $f$, or the corresponding hypersurface $X$ with equation
$f=0$, is \emph{totally Hessian} in the sense that $h(f)=cf^{\frac
{(d-2)(r+1)}d}$ with $c\in k\setminus \{0\}$. Here a good deal of
examples of such forms arises from the theory of
\emph{pre-homogeneous vector spaces},  a notion introduced by Kimura
and Sato (see \cite{KS}, also \cite{ESB}, \cite{Dolgachev},
\cite{EKP} and Remark \ref {tothess}). In this setup $f$ is the
so-called \emph{relative invariant} of the pre--homogeneous space,
uniquely defined up to a non--zero factor from $\C$. If, moreover,
its Hessian is non--zero then it is in fact totally Hessian and $f$
is a homaloidal polynomial such that $\phi_f$ coincides with its
inverse up to a projective transformation (see \cite[Theorem 2.8]
{ESB}).

As {mentioned}, the
singularities of a hypersurface $X\subset \p^ r$ which is either
homaloidal or has vanishing Hessian are not arbitrary.  For example,
in the second case, if $r\geq 3$ then $X$ cannot have isolated
singularities. The same result regarding homaloidal hypersurfaces
is a conjecture of
Dimca--Papadima (see \cite{DP}). We give a slight evidence for this
conjecture in terms of a resolution of the indeterminacies of the
polar map of $X$ by successive blowups  {
along smooth centers, to wit}, if $X\subset \p^r$ is
homaloidal and its degree exceeds $r+1$ then, for some blowing--up
step, the multiplicity of the proper transform of the general first
polar of $X$ is at least the dimension of the center of the
blowup (see Proposition~\ref{codimen}). In other words, the polar
linear system of $X$ cannot be \emph {log--canonical} (see \cite[pg.
56] {KM}). This gives a measure of the complexity of the singular
locus of $X$. In particular it shows that a homaloidal hypersurface
in $\p^3$, of degree at least $5$, cannot have ordinary
singularities.

After these preliminaries, we produce, for every  $r\geq 3$, the
promised infinite series of irreducible homaloidal hypersurfaces in
$\p^ r$ of arbitrarily large degree $d\geq 2r-3$. They are the dual
hypersurfaces to certain scroll surfaces with a line directrix. It
is relevant to observe that the present examples are not related to
the ones based on pre--homogeneous vector spaces as
mentioned above. Also they show, perhaps
{against} the ongoing folklore, that there are
plenty of homaloidal polynomials around. They even seem to be in
majority as compared to polynomials with vanishing Hessian, though a
complete classification does not seem to be presently at hand.

The full results are a bit too technical to be narrated here - we
refer to the main theorem of the section Theorem \ref {serie}, in
which one shows that the dual hypersurfaces to certain rational
scroll surfaces $Y(r-2,d-r+2)\subset {\p^r}$ are homaloidal and have
degree $d\geq 2r-3$. These examples are obtained via a rather
intricate geometric construction linking in an unexpected way
hypersurfaces with vanishing Hessian and homaloidal hypersurfaces. A
central piece is Theorem~\ref{fibrespace}, whose proof is fairly
technical but keeps a strong geometric flavor.
 We then dwell quite a bit into the structure of these scroll surfaces,
looking at their construction from various different angles in order
to fully apprehend their properties. Finally, in Theorem
\ref{Dhyperi} we produce different infinite families of homaloidal
examples in $\p^ r$, $r\geq 4$. These, though still
related to some scroll surfaces, { do
not seem in general to relate to}  hypersurfaces with
vanishing Hessian, {which adds to}  the
feeling that {the} classification of
homaloidal hypersurfaces {has still a long way to go.}

In addition we give a refined analysis of the nature of the
singularities of the homaloidal examples in $\p^3$ along with
an insight into the degree of the inverse
map. {That is, here} we deal
with the scroll $Y(1, d-1)$ {which, for $d=3$}
turns out to be a particular case of a series of degenerate
determinantal Hankel hypersurfaces considered in the following and
last section \ref {algex}.

This  latter construction, which has a more algebraic flavor,  is
based on a certain specialization of the generic
\emph {Hankel matrix}. The interest of
these examples {lies in} that, besides being irreducible
and of degree $r$, they fit a recent construct generalizing the
classical de Jonqui\`eres transformations (see \cite{{Pan0}}) and
boil down in particular cases to projections of certain scroll
surfaces. The full development of the nature of these homaloidal
hypersurfaces relates to several typical concepts of commutative
algebra. It also relates to the method devised in \cite{cremona}.
These examples do not come (either) from the theory of
pre-homogeneous vector spaces { either} since,
 for example, they are not totally Hessian.
{A marked} feature of these homaloidal hypersurfaces
is that the corresponding degree is the dimension of the ambient
space, while in most examples coming from pre-homogeneous vector
spaces the degree of the invariant polynomial is small with respect
to the number of variables.

Though somewhat exceptional, all these examples share in common the
property of having large degree with respect to the number of
variables. Additional inquiry could be made as to whether there are
families of totally Hessian polynomials, not necessarily homaloidal,
of arbitrary   large degree for any $r\geq
3$. Or even { be} wondered if there exists a
characterization of all homaloidal polynomials whose Hessian is a
non-zero multiple of a linear form such as is the case for the
Hankel degeneration examples constructed in the last section.

\section{Dual varieties of scroll
surfaces}\label{scrolls}

In this section we recall, with no proofs, some general and perhaps
mostly well known facts about projective
duality and dual varieties of scroll surfaces.   Standing
reference for this part {are}  \cite{Kleiman}, \cite{Zak1}, \cite{Ru}.

\subsection{Generalities}

Throughout this paper $k$ denotes {an algebraically closed} field of characteristic zero --
though many contentions herein will hold more generally.

Let $\p^ r=\p(V)$ be a projective space over $k$, where $V$ is a
$k$-vector space of dimension $r+1$. The {\it dual\/} projective
space of $\p^r$ is ${\p^ r}^
*=\p(V^ *)$, where $V^*={\rm Hom}_k(V,k)$. If $\Pi=\p(W)\subseteq \p^ r$, with $W\subset V$
a vector subspace of dimension $m+1$, then  the {\it orthogonal\/}
projective subspace $\Pi^ \perp\subseteq {\p^ r}^ *$ to $\Pi$ is
defined to be $\p((V/W)^*)=\p({\rm Ann}(W))\subset \p(V^*)$, where
${\rm Ann}(W)=\{f\in V^*\,|\, f(w)=0\, , \,\forall w\in W\}$. Note
that, geometrically, if one identifies ${\p^ r}^*$  with the linear system of
all hyperplanes in $\p^ r$, then $\Pi^ \perp$ is identified with the linear system
of all hyperplanes in $\p^ r$ containing $\Pi$ and has dimension
$r-m-1$.

Let $X\subset \Proj^ r$ be an irreducible projective variety of dimension $n$. For a
smooth point $x\in X$,  $\,T_{X,x}$ will denote the {\it embedded tangent space\/}
to $X$ at $x$, a subspace of dimension $n$.

The {\it conormal variety\/} $N(X)$ of $X\subset \Proj^ r$ is the
incidence variety defined as the closure of the set of all pairs
$(x,\pi)\in \p^ r\times{ \p^r}^  *$, such that $x$ is a smooth point
of $X$ and $\pi\in T_{X,x}^ \perp$ -- each such a hyperplane $\pi$
is said to be {\it tangent\/} to $X$ at $x$. Since the fiber of the
first projection $N(X)\rar X$ over a smooth point $x\in X$ is the
projective subspace $T_{X,x}^\perp\simeq \p^{r-n-1}$ of
hyperplanes containing $T_{X,x}$, then {$N(X)$ is irreducible and} $\dim(N(X))=r-1$.

The image of the projection of $N(X)$ to the second factor is, by
definition, the {\it dual variety} $X^ *$ of $X$. Since $k$ has
characteristic zero, one has $N(X)=N(X^ *)$ via the natural
identification  $\p^ r = ({\p^r}^*)^ *$ - a property known as {\it
reflexivity} (see, e.g., \cite{Kleiman}). It
follows  that $(X^*)^*=X$.

The {\it dual defect\/} of $X\subset\p^r$ is the non-negative integer
$d(X):=r-1-\dim(X^ *)$ and $X\subset\p^r$ is said to be {\it (dual) defective} if $d(X)>0$, i.e. if
$X^*\subset\p^{r*}$ is not a hypersurface. Note that $d(X)$ is the dimension of
$(T_{X^*,\xi})^ \perp\subset \p^r$ for smooth $\xi\in X^*$; thus, if
$\xi$ corresponds to the general hyperplane $\pi$ tangent to a
point $x\in X$ then $\pi$ is tangent at all points of
$(T_{X^*,\xi})^ \perp\subset \p^r$.

 {Also recall that
$X\subset \p^ r$ is said to be {\it degenerate\/} if its linear span
$\Pi=<X>$ is a proper subspace of $\p^ r$, i.e., if its homogeneous
defining ideal contains some nonzero linear form.}

Let now $\Pi\subset \p^ r$ be a subspace of dimension $m$, and let
$$\sigma_{_\Pi}: \p^ r\map (\Pi^ \perp)^ *\simeq \p^ {r-m-1}$$
be the projection from $\Pi$, defined as
$\sigma_{_\Pi}(p)=(\ell_1(p):\ldots:\ell_{r-m}(p))$, where
$\ell_1,\ldots,\ell_{r-m}$ are linear forms cutting $\Pi$ as a
linear subspace of $\p^r$. If $X\subset \p^ r$ is not contained in
$\Pi$, the closure $X_\Pi$ of the image  of $X$ via $\sigma_{_\Pi}$
is called the {\it projection of $X$ from $\Pi$}.  { If $\Pi\cap X=\emptyset$, then
$\sigma_\Pi$, or  $X_\Pi$, is said to be an {\it external}
projection of $X$. If $\dim (X)<r-m-1$ then $X_\Pi$ is a proper
subvariety of $(\Pi^ \perp)^*\simeq \p^ {r-m-1}$ and one has the
following:

\begin{Proposition}\label{projection_compatibility}
 With the previous notation, suppose that $X\subset\p^r$ is
non-degen\-erate and that $\dim (X)<r-\dim(\Pi)-1$. Then:
\begin{itemize}
\item[{\rm (i)}] $(X_\Pi)^*\subseteq \Pi^\perp\cap X^*$;

\item[{\rm (ii)}] If $\Pi^\perp\cap X^*$ is irreducible and reduced and if $\dim(X_\Pi)^*=\dim(\Pi^\perp\cap X^*)$, then $(X_\Pi)^*=\Pi^\perp\cap X^*$ as a scheme.
\end{itemize}
\end{Proposition}

\begin{proof} A general tangent hyperplane to $X_\Pi$ pulls back, via $\sigma _\Pi$, to a
hyperplane containing $\Pi$ and tangent to $X$ at a general point, proving (i). Part (ii) follows from (i).
\end{proof}

\begin{Proposition}\label{cones}  Let $\Pi=\p(W)\subset \p^r=\p(V)$ stand for the linear
span of the variety $X\subset \p^r$ and let $\widetilde{X}$ denote
the variety $X$ as re-embedded into $\Pi$. Then  $X^*\subset
{\p^r}^*=\p(V^*)$ is the cone over $\widetilde{X}^*\subset \p(W^*)$
with vertex $\Pi^ \perp=\p((V/W)^*)$. Conversely the dual of a cone
is degenerate, lying {on}  the orthogonal of the vertex of
the cone.
\end{Proposition}

The proof follows immediately from the aforementioned
interpretation of $\Pi^{\perp}$ as the set of hyperplanes in $\p^r$
containing $\Pi$.

Therefore, a subvariety $X\subset \p^ r$ is a cone if and only if
its dual $X^*\subset {\p^ r}^*$ is degenerate. Thus, the study of
dual varieties may safely be restricted to non--degenerate
varieties.

Finally recall that the {\it Gauss map\/} of an embedding $X\subset
\p^r$ is the map
$$\gamma_{_X}: x\in X\setminus {\rm Sing}(X)\map T_{X,x}\in \G(n,r).$$
The {\it image of the Gauss map\/} is the closure of $\gamma_{_X}(
X\setminus {\rm Sing}(X))$;  $\gamma_{_X}$ is said to be {\it
degenerate\/} if the fiber of $\gamma_X$ over a general point of its
image has positive dimension, i.e., if the Gauss image of $X$ in
$\G(n,r)$ has dimension at most $n-1$.  If $X\subset\p^r$ is a
smooth variety, then $\gamma_X$ is well
known to be finite and birational {onto}  its image, see
\cite[Theorem I.2.3]{Zak1}. More generally, the closure of the
general fiber of the Gauss image is a projective subspace (see
\cite[2.10]{GH}  or \cite {Zak1}).

\subsection{Scrolls and their dual varieties}

{ As mentioned in the
Introduction, scrolls  will play a substantial role in the
construction of the homaloidal hypersurfaces. Thus, we next proceed
to define them.}
\medskip

\begin{Definition}\label{scrolldef} An
irreducible variety $X\subset\p^r$ of dimension $n$ is said to be a
{\em scroll} if it is swept out by an irreducible $1$-dimensional family
$\mathcal F(X)$ of linear subspaces of $\p^r$ of dimension $n-1$,
called {\em rulings}, in such a way that
through a general point of $X$ there passes a unique member of
$\mathcal F(X)$.

Equivalently, let
$C$ be  the normalization of the defining $1$-dimensional parameter
space $\mathcal F(X)\subset \mathbb G(1,n-1)$  and let
$\pi:Y\to C$ denote the pull-back of the universal family on
$\mathbb G(1,n-1)$ restricted to $\mathcal F(X)$.
Then  $\pi:Y\to C$ is a $\p^{n-1}$-bundle
over $C$ and there exists a birational morphism $\phi:Y\to
X\subset\p^r$, induced by the tautological morphism on $\mathbb
G(1,n-1)$, such that the fibers of $\pi$ are embedded as linear
subspaces of $\p^r$.

With this terminology the scroll $X\subset\p^r$ is said to be {\em
rational} if $C\simeq\p^1$  and  {\em elliptic} if
$C$  has genus one. More generally we can define the
{\em genus of $X$} to be the geometric genus of $C$.

A scroll $X\subset\p^r$ is said to be {\it a smooth scroll} if
$\phi:Y\to X$ is an isomorphism. As in the
classical literature, a (smooth) scroll $X\subset\p^r$ is said to be
{\it normal} if $X\subset\p^r$ is a linearly normal projective
variety, i.e. if $X\subset\p^r$ is not a isomorphic linear external
projection of a
variety $\widetilde X\subset\p^{r+1}$.
\end{Definition}
\medskip

{It is well known that $\pi:Y\to C$ can be naturally identified
with   $\pi:\p(\mathcal E)\to C$, where $\mathcal E$ is rank
$n$ locally free sheaf over $C$. Moreover,  up to twisting by the
pull back of a line bundle on $C$, we can assume that $\phi$ is
given by a base point free linear system contained  in $|\O_{\p(\mathcal
E)}(1)|$. This linear system is complete if and only
if $X\subset\p^r$ is a normal scroll. Thus we can also assume that
$\mathcal E$ is generated by global sections and, if $C\simeq \p^1$,
that $\mathcal E\simeq \oplus_{i=1}^n\O_{\p^1}(a_i)$ for suitable
integers $0\leq a_1\leq\ldots\leq a_n$. In this case, if
$d=a_1+\ldots+a_n$, then $S(a_1,\ldots, a_n)\subset\p^{d+n}$ will
denote the rational scroll obtained as the image of the birational
morphism $\phi:\p( \oplus_{i=1}^n\O_{\p^1}(a_i))\to \p^{d+n}$ given
by the complete linear system $|\O(1)|$.  In this situation,
$d$ is the degree of $S(a_1,\ldots, a_n)\subset\p^{d+n}$.

In the above setting,  a smooth non--normal scroll $X\subset\p^r$ is
an external projection  of a normal smooth scroll. From the point of
view of the theory of dual varieties these examples are particularly
interesting since every smooth scroll $X\subset\p^r$ has $d(X)=n-2$
(see, e.g., \cite{Kleiman}). The simplest of these examples is
perhaps the Segre embedding $X={\rm Seg}(1,n-1)=S(1,\ldots,
1)\subset \p^{2n-1}$ of $\p^1\times \p^{n-1}$ -- here $\dim (X^*)=n$
and $X^*\subset\p^{2n-1*}$ is projectively equivalent to the
original $X$, i.e.  these Segre varieties are \emph {self-dual}.

In dimension $2$ the picture turns out to be the following. Consider a
non-degenerate surface $X\subset \p^ r$, $r\geq 3$. Here
 $n=2$, and  $d(X)=1$ if and only if $X$ is {\it developable}.
This condition  is equivalent to $\gamma_X$
being degenerate which {in turn happens to be} the
case if and only if $X$ is either a cone with vertex a point $p\in
\p^r$ or the tangent developable to a curve $C$, i.e., its
\emph{tangential surface}
$$X=\overline {\bigcup_ {x\in C\setminus{\rm Sing}(C)}T_{C,x}}$$
(see \cite[3.19]{GH}).  By Proposition~\ref{cones}, the first alternative takes place if and
only if $X^ *$ is degenerate, contained in the hyperplane $p^
\perp\subset {\p^r}^*$.

We collect further remarks in the form of a proposition for ready
reference.

\begin{Proposition}\label{non_developable} Let $X\subset \p^ r$ be a
non-degenerate scroll surface, $r\geq 3$.
Let $d$ denote the degree of $X$, which we assume to be at least $3$.
\begin{enumerate}
\item[{\rm (i)}] If  $X$ is not developable, then $X^*$ is a hypersurface of degree
$d$ which is swept out by the $(r-2)$--dimensional subspaces $F^
\perp$, where $F$ varies in the algebraic family $\F(X)$ determined
by the rulings of $X${\rm ;}
\item[{\rm (ii)}] Conversely, if $\,Y\subset {\p^ r}^ *$ is a hypersurface which is
swept out by a one-dimensional family $\mathcal F(Y)$ of subspaces of
dimension $r-2$, then $Y^ *\subset \p^r$ is either a $2$-dimensional
scroll or else a curve. Moreover, $Y^ *$ is a curve if and only if
one of the following equivalent conditions holds{\rm $\,$:}
\begin{enumerate}
\item[{\rm (a)}] $Y$ is \emph{developable}, that is to say, the general fiber of
the Gauss map $\gamma_{_Y}$  coincides with the general element of
$\mathcal F(Y)${\rm ;}
\item[{\rm (b)}] $\mathcal F(Y)$ is the family of the $(r-2)$--dimensional
subspaces $(r-1)$--osculating a curve.
\end{enumerate}
\end{enumerate}
\end{Proposition}
\demo Part (i) follows from the fact that a hyperplane $\xi$ is
tangent to $X$ if and only if it contains a ruling so that a general pencil
of hyperplanes cuts $X^*$ exactly in $d$ points. As
for (ii), see \cite[{Section}
2]{GH}. \qed

\subsubsection{Smooth rational normal scroll surfaces}\label{RNS}

We now go deeper into the structure of  rational
scroll surfaces.

Let $X=S(a,b)\subset\p^{a+b+1}$, $0<a\leq b$, be a smooth rational
normal scroll surface of degree $d=a+b$, in its standard embedding.
Recall that $S(a,b)$ is swept out by all lines joining corresponding
points on rational normal curves of degree $a$ and $b$ spanning
$\p^{a+b+1}$. This makes sense even if $a=0$, in which case $S(0,b)$
is the cone over a rational normal curve of degree $b$.
The homogeneous defining ideal is generated by the $2$-minors
of the piecewise $2\times(a+b)$ catalecticant matrix
$$\left(\begin{array}{cccc|cccc}
x_0 & x_1 & \ldots & x_{a-1} & x_{a+1} & x_{a+2} & \ldots & x_{a+b}\\
x_1 & x_2 & \ldots & x_{a} & x_{a+2} & x_{a+3} & \ldots & x_{a+b+1}
\end{array}
\right) $$
(see \cite {EH}).

We collect the main features of these scrolls.

\begin{Proposition}\label{scroll_features} Let $S(a,b)\subset\p^{a+b+1}$, $0<a\leq b$
be as above.
\begin{enumerate}
\item[{\rm (i)}] $S(a,b)$ is a  linear section of  the Segre embedding
${\rm Seg}(1,a+b-1)$ of $\p^{1}\times\p^{a+b-1}$ into
$\p^{2(a+b)-1}$ by a subspace $\Pi$ of dimension $a+b+1${\rm ;}
\item[{\rm (ii)}] $S(a,b)^*$ is the projection to ${\p^{a+b+1}}^*$ of
$\,{\rm Seg}(1,a+b-1)^*$, from $\Pi^{\perp}$, where $\Pi$ is a
subspace as in {\rm (i)}. In particular, $S(a,b)^ *$ is a
hypersurface in ${\p^{a+b+1}}^*$ of degree
$$\deg (S(a,b)^*)=\deg ({\rm Seg}(1,a+b-1)^ *)=\deg (S(a,b))=a+b;$$
\item[{\rm (iii)}]
As an abstract surface $S(a,b)$ is isomorphic to the so-called
\emph{ Hirzebruch surface\/} $\mathbb{F}_{b-a}=\Proj(\O_{\Proj^
1}\oplus \O_{\Proj^ 1}(a-b))$, with $\pi:\mathbb{F}_{b-a}\to\p^1$
the structural morphism. In particular, $S(a,b)$ admits a section
$E$ of $\pi$ with $E^2=a-b\leq 0$, which is unique if $a<b$.
Moreover, if $H$ is a hyperplane
section class of $S(a,b)\subset\p^{a+b+1}$, then $H\equiv E+bF$,
where $F$ is the class of the fibers of $\pi${\rm ;}
\item[{\rm (iv)}] Let $\Pi$ be a hyperplane in $\p^r$
tangent to $S(a,b)$ at finitely many points $p_1,\ldots p_m\in
S(a,b)$, $m\geq 1$. Let $H=H_\Pi$ be the corresponding hyperplane
section divisor. Then $H=F_{p_1}+\ldots+ F_{p_m}+ C,$ where
$F_{p_i}$ is the  ruling of $S(a,b)$ through the point $p_i$ and
$C\equiv E+(b-m)F$ is the divisor of a curve in $S(a,b)$ of degree
$a+b-m$ passing through $p_1,\dots, p_m${\rm ;}
\item[{\rm (v)}] If $m\leq a$ and $3m\leq a+b+1$, then the general hyperplane section
of $S(a,b)$ tangent at $m$ general points has exactly $m$ ordinary
quadratic singularities there and it is smooth elsewhere{$\,$\rm ;}
\item[{\rm (vi)}] If either $m\leq a$ or $m=b$ then the general curve $C\in | E+(b-m)F|$
is smooth and irreducible and, together with $m$ distinct fibres
$F_1,\dots ,F_m$ of $\pi$, gives rise to a hyperplane section
tangent at the  intersection points $p_i$ of $F_i$ with $C$,
$i=1,\dots,m$, and nowhere else.
\end{enumerate}
\end{Proposition}
\demo (i) This is clear from the above algebraic description of
$S(a,b)$  and the corresponding
defining equations of the Segre embedding as given by the $2$-minors
of a generic $2\times (a+b)$ matrix over $k$.

(ii)  As we pointed out already,  these Segre varieties are
self--dual, i.e. ${\rm Seg}(1,a+b-1)^ *$ is projectively equivalent
to ${\rm Seg}(1,a+b-1)$. Since $\Pi\cap {\rm Seg}(1,a+b-1)=S(a,b)$
is reduced and irreducible, then, by Proposition ~\ref{projection_compatibility}},
$S(a,b)^ *$ coincides with the projection of ${\rm Seg}(1,a+b-1)^ *$  from $\pi^ \perp$.
Note that $\Pi^\perp\cap
{\rm Seg}(1,a+b-1)^*=\emptyset$ since $\Pi\cap {\rm Seg}(1,a+b-1)$
is smooth.  Then the degree of $S(a,b)^ *$ is
the same as the degree of ${\rm Seg}(1,a+b-1)^ *$, which is the same
as the degree of $S(a,b)$, namely $a+b$.

(iii) The first part is well-known (see, e.g., \cite{EH})  and the
rest follows from this.

(iv) Describing  $S(a,b)^*\subset\p^{a+b+1}$ is the same as
describing the singular hyperplane sections of $S(a,b)$, i.e. those
given by hyperplanes $\Pi$ containing tangent planes of $S(a,b)$. If
$\Pi\supseteq T_{S(a,b),x}$, then $\Pi\supseteq F_x$, the line of
the ruling through $x$. Thus, if  $\Pi$ is a hyperplane tangent to
$S(a,b)$ at finitely many points $p_1,\ldots p_m$, $m\geq 1$ and
$H=H_\Pi$ denotes the corresponding hyperplane section divisor, it
is clear that $H=F_{p_1}+\ldots+ F_{p_m}+ C$, where $C\equiv
E+(b-m)F$ is the divisor of a curve of degree $a+b-m$ in $S(a,b)$.
Moreover $H_\Pi$ has to be singular at $p_1,\dots,p_m$, hence $C$
contains $p_1,\dots,p_m$.

(v) If $a\geq m$, then $(E+(b-m)F)^ 2=b+a-2m\geq 0$. Let $C$ be the
general curve in $| E+(b-m)F| $. The exact sequence

$$0\to \O_{S(a,b)}\to \O_{S(a,b)}(C)\to \O_C(C)\to 0$$
implies that the linear system $|E+(b-m)F| $ is base point free of
dimension $b+a-2m+1$ and its general curve $C$ is smooth and
rational. If $b+a\geq 3m-1$, the general curve in $|E+(b-m)F| $
contains $m$ general points of $S(a,b)$. This proves the assertion.

(vi) If either  $m\leq a$, or if $m=b$, the general such curve $C$
is smooth and irreducible (see the above argument). The assertion
follows. \qed

\bigskip

Next we highlight the nature of the singularities of the dual
$S(a,b)^*$. Let $E\subset S(a,b)$ be as in
Proposition~\ref{scroll_features}, (iii).


\begin{Proposition}\label{sing_of_dual} Let $S(a,b)\subset\p^{a+b+1}$, $0<a\leq b$
be as above.
\begin{enumerate}
\item[{\rm (i)}] The points of $S(a,b)^*$ corresponding to
hyperplanes tangent to $S(a,b)$ at $m$ distinct points are points of
multiplicity at least  $m$ of $S(a,b)^*${\rm ;}

\item[{\rm (ii)}] The singularities $\,\Sing(S(a,b)^ *)$ have a natural
stratification into locally closed sets $S^*_{\alpha}(a,b)$, with $2\leq \alpha\leq
a$ and $\alpha=b$, consisting of points of multiplicity at most
$\alpha$; as for $\alpha=b$, one has $S^*_{b}(a,b)=<E>^\perp\subset
S(a,b)^ *$, a linear space of dimension $b$ contained in
$\Sing(S(a,b)^ *)${\rm ;}
\item[{\rm (iii)}] {\rm (}$a=1${\rm )} The stratum
$\Sing(S(1,d-1)^*)=S^ *_{d-1}(1,d-1)$ is the subspace $<E>^ \perp$
of dimension $d-1$, whose general points are points of multiplicity
$d-1$ of the hypersurface $S(1,d-1)^*\subset {\p^{d+1}}^*$.
\end{enumerate}
\end{Proposition}
\demo  (i) Quite generally, the points of $S(a,b)^ *$ corresponding
to hyperplanes tangent to $S(a,b)$ at $m$ distinct points, with
$m\leq a$ or $m=b$, are points of multiplicity at least $m$. One
sees that $S^*_{b}(a,b)=<E>^\perp\subset S(a,b)^ *$ is a linear space of
dimension $b$ contained in $\Sing(S(a,b)^ *)$. We now suppose that $2\leq
m\leq a$.

We know that $S(a,b)^*$ is a hypersurface. If a point of $S(a,b)^*$
corresponds to a hyperplane $H=F_{p_1}+\ldots+F_{p_m}+C$ tangent to
$S(a,b)$ at the $m$ points $p_1,\ldots, p_m$, one sees that there
are at least $m$ distinct branches of $S(a,b)^*$ passing
through $H$, namely the ones corresponding to hyperplane sections of
the form $F_{p_i}+C_i$, $C_i$ irreducible and smooth, proving the
assertion.

Assertions (ii) and (iii) follow from (i).
 \qed

\medskip

Notice that the scheme structure on $\Sing(S(a,b)^*)$ defined by the
partial derivatives of the defining equation of $S(a,b)^*$ has
embedded points (see \cite{Aluffi} for some interesting
considerations on this scheme structure on $\Sing(S(a,b)^*)$).

It is { classically known}
 that a non--developable scroll surface in $\p^3$  is self--dual. We prove this result anew in the case
where the scroll is rational, which is our  main focus. The proof
contains elements for later use.

\begin{Proposition}\label{autodual} Let $X\subset \p^ 3$ be a rational scroll which
is not developable. Then
$X$ is self--dual, i.e. there is a projective transformation sending $X$ to $X^ *$.
\end{Proposition}

\begin{proof}
 {By
definition} $X\subset\p^3$ is the
birational projection to $\p^3$ of a smooth rational normal scroll
surface $S(a,b)\subset  \p^ {a+b+1}$, with $0<a\leq b$, from a
subspace $\Psi$ of dimension $a+b-3$ such that $S(a,b)\cap
\Psi=\emptyset$.

By part (i) of Proposition \ref {projection_compatibility}, we have

\begin{equation}\label {eq:inclus}
X^ *\subseteq \Psi^ \perp\cap S(a,b)^ *.
\end{equation}

The right hand side is a hypersurface of degree $a+b$ in $\p^ 3$.
Moreover $X^ *$ is also a hypersurface, since $X$ is not
developable, and its degree is $a+b$ (see part (i) of Proposition
\ref{non_developable}). Then equality holds in \eqref {eq:inclus},
i.e.

\begin{equation}\label {eq:inclus2}
X^ *= \Psi^ \perp\cap S(a,b)^ *.
\end{equation}

By (i) of Proposition \ref
{scroll_features}, $S(a,b)=\Pi\cap {\rm Seg}(1,a+b-1)$ with $\Pi$ a
subspace of dimension $a+b+1$ of $\p^{2(a+b)-1}$. Thus

$$
\begin{aligned}
&X^ *=\Psi^ \perp\cap ({\rm Seg}(1,a+b-1)\cap\Pi)^ *=&\cr
&=\Psi^ \perp\cap \sigma_{\Pi^ \perp}({\rm Seg}(1,a+b-1)^ *).&\cr
\end{aligned}
$$

Therefore, up to a projective transformation

$$
\begin{aligned}
&X^ *=\Psi^ \perp\cap \sigma_{\Pi^ \perp}(({\rm Seg}(1,a+b-1))=&\cr
&=\sigma_{\Pi^ \perp}( \langle \Pi^ \perp, \Psi^ \perp\rangle\cap({\rm Seg}(1,a+b-1))=X.&\cr
\end{aligned}
$$
\end{proof}

\subsubsection{Multiple line directrix on scrolls}\label{multidir}

We now consider  {another interesting class of} scroll surfaces. Any
non--developable rational scroll surface $X\subset \p^ r$ of degree $d$
is a {birational external} projection
of a scroll $S(a,b)$ with $d=a+b$ and $X^
*$ is a section of $S(a,b)^ *$ by
Proposition \ref {autodual}.
If there is a line $L\subset X$
such that $X$ is smooth at the general point of $L$  and $L$ meets
the general ruling of $X$ at one single point, then $X$ is the
projection of $S(1,d-1)$. In such a case $X^ *$ is a hypersurface
 of degree $d$ for which the $(r-2)$--dimensional subspace $L^
\perp$ has multiplicity $d-1$.

Such a line $L$ is called a {\it simple line directrix}. More
generally, a line $L\subset X$ is a {\it line directrix of
multiplicity\/} $e:=e(X)$ if the general point $x\in L$ has
multiplicity $e$ for $X$ and there is some line in $\F(X)$,
different from $L$, passing through $x$. Note that $L$ may, or may
not, belong to $\F(X)$. It is clear that a scroll with a line
directrix is not developable, unless it is a plane. Therefore in
what follows we will implicitly assume that a scroll with a line
directrix is not developable.

\medskip

\begin{Proposition}\label{directrix}
Let $L$ be a line directrix of multiplicity $e$ on a rational scroll
surface $X\subset \p^ r$ of degree $d$. Let $\mu:=\mu(X)$ denote
the number of rulings in $\F(X)$ {not
coinciding with} $L$ and passing through a
general point ${ x}$ {of} ${L}$ and let
$F_{x,i}, i=1,...,\mu$, be these rulings. Let $\nu:=\nu(X)$ be the
dimension of the span  $<L, F_{x,1},...,F_{x,\mu}>$.

One has:
\begin{enumerate}
\item[{\rm (i)}] $\mu\leq e${\rm ;}
\item[{\rm (ii)}]
$\mu<e$ if and only if $L$ is a ruling in $\F(X)${\rm ;}
\item[{\rm (iii)}] The dual $X^ *\subset {\p^r} ^*$ is a hypersurface of degree $d$
and contains the $(r-2)$-dimensional subspace $\Pi=L^ \perp$.
Moreover $X^*$ has multiplicity $d-\mu$ at the general point of
$\Pi$ and the general hyperplane through $\Pi$ cuts out on $X^ *$
off $\Pi$, the union of $\mu$ codimension 2 subspaces whose intersection
with $\Pi$ is a subspace of dimension $r-\nu-1$.
\end{enumerate}
\end{Proposition}

\demo (i) Let $f: \bar X\to X$ be the normalization morphism. The
surface $\bar X$ is ruled and its rulings are mapped to the lines in
$\F(X)$.  Let $p_1,\dots,p_h$ be the points on
$\bar X$ mapping to $x$.  It is clear that $e\geq h$. Note that
$\bar X$, which is normal and therefore smooth in codimension one,
is smooth at $p_1,\dots, p_h$. Hence there is a unique ruling of
$\bar X$ through each of the points $p_1,\dots, p_h$. Moreover the
$\mu$ rulings in $\F(X)$, different from $L$, and passing through
$x$ are the image, via $f$, of rulings on $\bar X$ passing through
one of the points $p_1,\dots, p_h$. Thus $\mu\leq h\leq e$.

(ii)  Let us prove that, if $\mu<e$, then $L$ is a ruling of $\F(X)$.
The converse is similar and can be left to the reader.

Suppose first that $\mu\leq h<e$. This is equivalent to say that $f$ is
ramified at some of the points $p_1,\dots, p_h$, which we denote
by $y$. Let $F$ be the ruling of $\bar X$ passing through $y$.
Since $F$ maps to a line via $f$, the only possibility is that $F$
maps to $L$, hence $L$ is a ruling of $\F(X)$ in this case.

Suppose that $h=e$, i.e. $f$ is unramified at a general point $x\in L$. One
has therefore $e$ distinct points $p_1,..,p_e$ on $\bar X$ mapping
to $x$. Let $F_i$ be the ruling through $p_i, i=1,...,e$. If $L$ is
not a ruling in $\F(X)$, then the images on $X$ of $F_1,...,F_e$
are all distinct from $L$. Moreover they are also $e$ distinct lines,
since $f$ is a finite birational map. Hence $\mu=e$, proving the assertion.

(iii) The general hyperplane $\Xi=x^ \perp$ containing $\Pi=L^ \perp$
corresponding to the general point $x\in L$ cuts out $X^*$, off
$\Pi$, along the union of the $\mu$ codimension $2$ subspaces
$F_{x,i}^ \perp, i=1,...,\mu$. The intersection
$$\Pi^  \perp\cap F_{x,1}^ \perp\cap...\cap
F_{x,\mu}^ \perp=<L, F_{x,1},...,F_{x,\mu}>^ \perp,$$
has dimension $r-\nu-1$. \qed

\medskip

We will see  later how to construct scrolls with $\mu<e$ (see
Lemma~\ref{ratcurve} and ff.). As for the case $\mu=e$, the
following construct works: consider $S(a,b)\subset \p^ {d+1}, d=a+b,
a\geq 2$, and project it down to $\p^ {b+2}$ from a general linear
space of dimension $a-2$ which sits in $<E>$. In this way the image
$X(a,b)$ of the projection has still degree $d$ and the image of $E$
is the line $\Lambda$ to which $<E>$ maps. Notice that by projecting
$X(a,b)$ from $\Lambda$ to $\p^ b$ one gets a rational
normal curve $C$ of degree $b$. Thus
$X(a,b)$ sits {on} the $3$--dimensional cone of degree $b$
projecting $C$ from $\Lambda$.

Since $E$ has degree $a$, one has that $\Lambda$ is a line directrix
of multiplicity $a$ and clearly $\mu=a$. In this case $\nu=a+1$
(see the argument in the proof of parts (i) and (ii) of  Proposition \ref{directrix} above).
Notice that $X(a,b)$ is contained in a cone $Z(a,b)$ of dimension $a+2$
which is swept out by the subspace $<\Lambda, F_{x,1},...,F_{x,a}>$
of dimension $a$  as $x$ varies on $\Lambda$. The cone $Z(a,b)$ is a
rational normal scroll of degree $b-a+1$ (see \cite {EH}).

One can also obtain the previous example in terms of the dual
variety of certain projections of more general
 scrolls, {as follows.}

Let $1\leq a=a_1\leq a_2\leq a_3\leq\ldots\leq a_{r-1}$ be integers
and set $d=\sum_{i=1}^{r-1}a_i$. Consider the rational normal scroll
$X_1=S(a,a_2,\ldots,a_{r-1})\subset\p^{d+r-2}$ of degree $d$ and
dimension $r-1$, with
$$S(a,a_2,\ldots,a_{r-1})\simeq \Proj(\O_{\Proj^ 1}(a)\oplus\O_{\Proj^
1}(a_2)\oplus... \O_{\Proj^ 1}(a_{r-1}))$$ embedded via the $\O(1)$
bundle. Algebraically, the homogeneous defining ideal of this
embedding is generated by the $2$-minors of a multi-piecewise
catalecticant matrix as in (\ref{RNS}). Consider also the rational
normal scroll $X_2=S(a_2,\ldots,a_{r-1})\subset\p^{d+r-a-3}$ of
degree $d-a$ and dimension $r-2$. By a suitable identification, one
can consider $X_2$ as a subvariety of $X_1$. Let $\Omega$ be a
sufficiently general linear space of dimension $d-3$ which cuts the
linear space $<X_2>$ along a subspace of dimension $d-a-2$, and set
$Y=\sigma_{_{\Omega}}(X_1)\subset\p^r$, where $\sigma_{_{\Omega}}$
as before stands for the projection from $\Omega$.

\begin{Proposition}\label{multi_scrolls} Let the notation be as above,
with $Y=\sigma_{_{\Omega}}(X_1)\subset\p^r$. The dual $Y^*\subset
{\p^r}^*$ is a scroll surface of degree $d$, with $e(Y^
*)=\mu(Y^*)=a$ and with  line directrix $\sigma_{_{\Omega}}(<X_2>)$ of
multiplicity $a$.
\end{Proposition}
\demo Clearly $\deg(Y)=\deg(X_1)=d$ and $\Sing(Y)$ contains the
linear space $\Pi=\sigma_{_{\Omega}}(<X_2>)$ of dimension $r-2$. The
general point of $\Pi$ has multiplicity $d-a$ for $V$ because
$\deg(X_2)=d-a$.

Since $Y$ is swept out by a 1--dimensional family of projective
spaces of dimension $r-2$, then $X=Y^ *$ is a scroll surface of
degree $d$ with line directrix $L=\Pi^ \perp$.  Since $Y$  has
multiplicity $d-a$ along $\Pi$, we see $\mu(X)=a$. Actually the
multiplicity of the line directrix $L$ on $Y^*$ is also $a$ because
clearly $L$ is not a ruling in $\F(X)$. \qed

\section{The polar map of a projective hypersurface} \label {polarmap}

Let $f=f(\mathbf{x})=f(x_0,\ldots,x_r)\in k[x_0,\ldots,x_r]$ be a non-zero
homogeneous polynomial of degree $d$ in the $r+1$ variables $x_0,\ldots,x_r$ over an
algebraically closed field $k$ of characteristic zero.

Then $V(f)\subset\p^r$ will denote the hypersurface  scheme theoretically defined by
the equation $f(x_0,\ldots,x_r)=0$, so $V(f)$ might not be reduced. Its {\it
support\/} ${\rm Supp}(V(f))$ is the set of points of $\Proj^ r$ where $f$ vanishes.

We will often denote by $f_i$ the partial derivative $\frac {\partial f} {\partial
x_i}$, $i=0,...,r$.

Let $\mathbf {p}=(p_0,\ldots,p_r)\in k^ {r+1}\setminus \{0\}$, and let
$p= (p_0,\ldots,p_r)$ denote the corresponding point in $\Proj^ r$. For
every positive integer $s< d$ consider the polynomial

$$\Delta^ s_\mathbf {p} f(\mathbf{x})=(p_0\frac {\partial}
{\partial x_0}+...+ p_r\frac {\partial} {\partial x_r})^ {(s)}f(\mathbf{x})$$ where
the exponent $s$ in brackets means, as usual, a symbolic power involving products
and derivatives. The polynomial $\Delta^ s_\mathbf {p}f$ has degree $d-s$ and, for
any $t\in k^ *$, one has:

$$\Delta_{t\mathbf{p}}^ sf(\mathbf{x})=t^s \Delta_\mathbf {p}^
sf(\mathbf{x}).$$
If $\Delta^ s_\mathbf {p} f$ is not identically
zero, then it makes sense to consider  the hypersurface $V(\Delta^
s_\mathbf {p} f)$ which depends only on $p$ and on $V(f)$ and is
called the $s$--{\it th polar} of $V(f)$ with respect to $p$. We
will denote it by $V_p^ s(f)$. If $\Delta^ s_\mathbf {p} f$ is
identically zero, one says that the $s$--{\it th polar} $V_p^ s(f)$
of $V(f)$ with respect to $p$ {\it vanishes identically}. In this
case we consider $V_p^ s(f)$ to be the whole $\Proj^ r$.

For general properties of polarity, which we will freely use later
on, we refer to \cite{Segrepr}. Among these we mention here the so
called {\it reciprocity theorem\/}:

\begin{Proposition}\label{recip} Given the hypersurface $V(f)$ in $\p^ r$
and two points $p= (p_0,\ldots,p_r)$, $q= (q_0,\ldots,q_r)$, one has:
$$\frac 1{s!} \Delta_\mathbf {p}^ s f(\mathbf{q})=\frac 1{(d-s)!}\Delta_\mathbf {q}^ {d-s}f(\mathbf{p}).$$
Thus $q\in V_p^ s(f)$ if and only if $p\in V_q^ {d-s}(f)$.
\end{Proposition}

As $p$ varies in $\Proj^ r$, the polars $V_p^ s(f)$ do not vary in a
linear system, unless $s=1$. The base locus scheme of the linear
system $\mathcal {P}(f)$ of the first polars of $V(f)$ is the {\it
singular locus\/} ${\rm Sing}(V(f))$ of $V(f)$, defined by the {\it
Jacobian\/} (or {\it gradient\/}) ideal generated by the partial
derivatives $f_0(\mathbf{x}),\cdots,f_r(\mathbf{x})$.

A consequence of the reciprocity theorem is that the polar
hyperplane $\pi_p(f):=V_p^ {d-1}(f)$ has equation:
$$f_0(\mathbf{p})x_0+\ldots+f_r(\mathbf{p})x_r=0$$
which vanishes identically if and only if $p\in {\rm Sing}(V(f))$.
If $p\in V(f)$ and it is not singular, then $\pi_p(f)$ is the {\it
tangent hyperplane} $T_{V(f),p}$ to $V(f)$ at $p$.

The ({\it first\/}) {\it polar map\/} of $f$ or of $V(f)$ is the
rational map
$$\phi_f: x=(\mathbf{x})\in \p^r\map
(f_0(\mathbf{x}),\ldots,f_r(\mathbf{x}))\in \p^{r}.$$ It can be
interpreted as mapping the point $p$  to its polar hyperplane
$\pi_p(f)$ and, as such, its target is ${\p^ r}^ *$.

In terms of linear systems $\phi_f$ is the map defined by the system
$\mathcal {P}(f)$ of the first polars. Thus, if $V(f)$ is reduced,
as we will now assume,  the indeterminacy locus
of $\phi_f$ is ${\rm Sing}(V(f))$. The restriction  of $\phi_f$ to
$V(f)\setminus\Sing(V(f))$ is the {Gauss map} of $V(f)$, hence the
corresponding image is the {dual variety} $V(f)^
*$ of $V(f)$. We will set $v(f)=\dim (V(f)^ *)$.

Denote by $Z(f)$ the closure of the image of $\Proj^ r$ via $\phi_f$
-- called the {\it polar image\/} of $f$ -- and set
$z(f)=\dim(Z(f))$. Clearly $v(f)\leq z(f)$, but we shall see in a moment
that strict inequality holds (see Remark \ref {hessim})

We denote by $\delta(f)$ the degree of the map $\phi_f$, which is
meant to be $0$ if and only if $z(f)<r$, otherwise it is a positive
integer. We will call $\delta(f)$ the {\it polar degree\/} of
$V(f)$. Let $f_{\red}$ be the {\it radical of $f$}, i.e.  $(f_{\red})=\sqrt{(f)}$.

We record the following result from \cite[Corollary 2]{DP} which
proves a conjecture stated in \cite{Dolgachev}:

\begin{Theorem} \label{nored} Let notation be as above.
Then $\delta(f)=\delta(f_{\red})$,  i.e. the polar degree of $V(f)$
depends only on ${\rm Supp}(V(f))$.
\end{Theorem}

This result enables us to restrict our attention to reduced
hypersurfaces if we are interested in studying the polar degree. The
argument in \cite{DP} depends on topological considerations.
For different proof, see \cite {FP} , whereas an
algebraic  proof of the case where the irreducible factors of $f$
are of degree one has been established in \cite{Bruno}.

\subsection{The Hessian of a projective hypersurface}\label{hess_of_hyp}

Consider now the $(r+1)\times(r+1)$ {\it Hessian matrix\/} of
$f(\mathbf{x})$

$$h(f)(\mathbf{x}):=\det\left(\frac{\partial^2 f(\mathbf{x})}{\partial x_i\partial
x_j}\right)_{i,j=0,...,r}.$$ Its determinant $h(f)\in
k[x_0,\ldots,x_r]$ is the {\it Hessian polynomial\/} of
$f(\mathbf{x})$.

Sometimes we will abuse notation and denote by $h(f)$ also the
Hessian matrix rather than its determinant, hoping no ambiguity will
 be caused.

We note that the Hessian is covariant by a linear change of
variables. If $h(f)$ is a non-zero polynomial,  the {\it Hessian\/}
of the hypersurface $V(f)\subset \p^r$ is the hypersurface
$H(f):=V(h(f))$. Otherwise we say that $V(f)$ has {\it vanishing\/},
or {\it indeterminate\/} Hessian, in which case we consider $H(f)$
to be the whole of $\Proj^ r$.

A couple of basic remarks is in order.

\begin{Remark}\label {hesspol} A point $p\in \p^r$ belongs to $H(f)$ if and only if
the polar quadric $Q_p(f):=V_p^ {d-2}(f)$ is either singular or
vanishes identically.

Thus, in particular $p\in V(f)\cap H(f)$ if and only if either $p\in
{\rm Sing}(V(f))$ or $p$ is a {\it parabolic point\/} of $V(f)$ in
the sense that the tangent cone $A_p(f)$ at $p$ of the intersection
of $V(f)$ with the tangent hyperplane $\pi_p(f)$ (necessarily
singular at $p$) has a vertex of positive dimension (see \cite[p.
71]{Segrepr}). This cone is called the {\it asymptotic cone\/} of
$V(f)$ at $p$. More precisely, a point $p\in V(f)$ is said to be
$h$--{\it parabolic}, $h\geq 0$, if the vertex of the asymptotic
cone $A_p(f)$ has dimension $h$. In that case $p$ is a point of
multiplicity $h$ for $H(f)$ (see \cite{Segre}). Note that
$0$--parabolic means \emph {non--parabolic}. If $f$ is irreducible
and the general point of $V(f)$ is $h$--parabolic, then $f^ {h}$
divides $h(f)$; in particular, if $h>0$, then $V(f)$ is contained in
$H(f)$.

Conversely, if $f$ is irreducible and $V(f)$ is contained in $H(f)$
then the Gauss map of $V(f)$ is degenerate, i.e. $v(f)<r-1$ and the
general point $p\in V(f)$ is $h$--parabolic with $h=r-v(f)-1>0$ (see
\cite[4--5]{Seg},  \cite{Segre}, \cite{Ciliberto}). In this case,
since the general fibre of the Gauss map is a linear space, then
$V(f)$ is described by an $(r-h-1)$--dimensional family of
$h$--dimensional linear subspaces
 of $\Proj^ r$, parameterized by $V(f)^*$.
 Moreover $H(f)$  contains  $V(f)$ with multiplicity at least $h=r-v(f)-1$.

The question as to when  $H(f)$  contains  $V(f)$ with higher
multiplicity than the expected value $r-v(f)-1$ has been considered
in \cite {Segre}, \cite{Segre2}, \cite{Fr1}, \cite{Ciliberto}.
\end{Remark}

\begin{Remark}\label{hessim} A point $p\in \p^r$ belongs to $H(f)$ if
and only if the rank of the map $\phi_f$ at $p$ is not maximal, i.e.
if and only if ${\rm rk}_p(\phi_f)<r$. Hence $z(f)<r$ if and only if
$V(f)$ has vanishing Hessian. Set $\rho(f):={\rm rk}(h(f))$,  where
the rank of $h(f)$ is computed as a matrix over the field
$k(x_0,...,x_r)$, or, what is the same, at a general point of
$\Proj^ r$. Then one has:
$$z(f)=\rho(f)-1$$
Indeed, if $p=(p_0,\ldots,p_r)$ is a point in $\Proj^r$ not on ${\rm
Sing} (V(f))$, and if $\xi=\phi_f(p)$, then $T_{Z(f),\xi}$ is
spanned by $\xi$ and by the points $(f_{i0}(\mathbf{p}),\ldots,
f_{ir}(\mathbf{p})),$ $ i=0,...,r$.
Notice that $T_{Z(f),\xi}^ \perp$ is the vertex of $Q_p(f)$.
A vastly more general principle
holds in this connection (see \cite[Proposition 1.1]{Simis} for a
detailed argument).

Notice that, if $V(f)$ is irreducible and its general point $p$ is $h$--parabolic,
then $v(f)+2=r-h+1=\rk(Q_p(f))\leq \rho(f)=z(f)+1$, i.e. $v(f)<z(f)$,
namely the dual $V(f)^ *$  of $V(f)$ is properly contained in the
polar image $Z(f)$.

Note that, by Theorem~\ref{nored}, the property of having vanishing
Hessian only depends on the support of a hypersurface. Thus, if one
is interested in hypersurfaces with vanishing Hessian, one can
restrict the attention to the reduced ones.
\end{Remark}

\subsection{Hypersurfaces with vanishing Hessian}\label{hess}

The hypersurface $V(f)$ has vanishing Hessian if and only the
derivatives $f_0,\ldots, f_r$ are {\it algebraically dependent},
i.e. if and only if there is some non--zero polynomial
$g(x_0,\ldots,x_r)\in k[x_0,\ldots,x_r]$ such that
$g(f_0,\ldots,f_r)= 0$.

Note that $V(f)$ is smooth if and only if $f_0,\ldots, f_r$ form a
{\it regular sequence}; in particular, if $V(f)$ is smooth then
$h(f)\neq 0$. Thus, having vanishing Hessian implies at least that
${\rm Sing}(V(f))\neq \emptyset$ and one then asks how big is this
locus.

The following result due to Zak (see \cite[Proposition 4.9]{ZakHesse})
partially answers this question. Part of it can be traced back to
Gordan--Noether (see \cite {GN}).

\begin{Proposition}\label{hessian} Let $X=V(f)\subset\p^r$ be a reduced
hypersurface with vanishing Hessian and let $Z(f)\subset {\p^r}^*$
denote the polar image of $f$. Then
\begin{itemize}
\item[{\rm (i)}] The closure of the fiber of the map $\phi_f$ over a general point
$\xi\in Z(f)$ is the union of finitely many linear subspaces of
dimension $r-z(f)=r-\rho(f)+1$, passing through the subspace
$(T_{Z(f),\xi})^\perp${\rm ;}
\item[{\rm (ii)}] $Z(f)^ *$ is contained in ${\rm Sing}(V(f))$.
\end{itemize}
\end{Proposition}

The careful reader will notice that the argument in
\cite[Proposition 4.9]{ZakHesse}
actually proves the above statement
(i) rather than the corresponding part (ii)
of the statement there.

A clear-cut case of vanishing Hessian is when, $f_0,\dots,f_r$
are linearly dependent, i.e. up to a linear change
of variables, $f$ does not depend on all the variables, i.e., when
$V(f)$ is a cone (Proposition~\ref{cones}). One could naively ask
for the converse:

\begin{Question}\label{Hessequestion}{\rm ({\bf Hesse problem})}
\rm Does $h(f)= 0$ imply that the derivatives $f_0,\ldots,f_r$ are
linearly dependent?
\end{Question}

Hesse claimed this twice (see \cite{Hesse1}, \cite{Hesse2}), however the
proofs had a gap. The question was taken up by Gordan and Noether
in \cite {GN}, who showed that the question has an affirmative
answer for $r\leq 3$, but is false in general for $r\geq 4$. Their
methods have been revisited in more recent times by Permutti in
\cite {Permutti1}, \cite {Permutti3} and \cite{Lossen}.

Using Proposition~\ref{hessian} we can give an easy proof of this
fact for $r\leq 2$. The case $r=3$ is slightly more complicated and
will not be dealt with it here - we refer to \cite {GN},  \cite
{Fr1} or \cite{Lossen}. A simple proof is also contained in
\cite {GR}.

\begin{Proposition}\label{lowerN} Let $V(f)\subset\p^r$, $1\leq r\leq 2$,
be a  reduced hypersurface of degree $d$. Then $V(f)$ has vanishing
Hessian if and only if $V(f)$ is a cone. More precisely, $V(f)$ has
vanishing Hessian if and only if either $r=1$ and $d=1$, or else
$r=2$ and $V(f)$ consists of $d$ distinct lines through a point.
\end{Proposition}
\begin{proof}
If $r=1$, then $Z(f)\subset\p^1$ must be a point, so the partial
derivatives of $f$ are constant and $d=1$.

Suppose $r=2$. Then $z(f)\leq 1$. As above, $Z(f)$ is a point if and
only if $d=1$. Let $z(f)=1$. From  part (ii) of
Proposition~\ref{hessian}, we have that $Z(f)^*\subset\Sing(V(f))$.
Since we are assuming $V(f)$ to be reduced,  we have that $Z(f)^*$
is a point, so that $Z(f)$ is a line, hence degenerate. This is
equivalent to saying that $V(f)$ is a cone.
\end{proof}

\begin {Remark} \label{arrangm} It is interesting to note that the only
hyperplane arrangements with vanishing Hessian are cones (see \cite[Cor. 2 and Cor. 4]{DP}).
\end {Remark}

\subsection{Gordan--Noether counterexamples to Hesse's problem}\label {GNP}

 We will now briefly
recall the results of Gordan--Noether and Permutti in connection
with the Hesse problem, which showed that Hesse's argument was
faulty for dimension $r=4$ and higher.

Thus, assume that $r\geq 4$ and fix integers $t\geq m+1$ such that
$2\leq t\leq r-2$ and $1\leq m\leq r-t-1$. Consider forms
$h_i(y_0,...,y_m)\in k[y_0,...,y_m]$, $i=0,...,t$, of the same
degree, and also forms $\psi_j(x_{t+1},...,x_r)\in
k[x_{t+1},...,x_r]$, $j=0,...,m$, of the same degree. Introduce the following
homogeneous polynomials all of the same degree:

$$
Q_\ell(x_0,...,x_r)=\det \left (
\begin{matrix}
x_0&...&x_t \\
\frac {\partial h_0}{\partial \psi_0}&...&\frac {\partial h_t}{\partial \psi_0}\\
... & ...& ... \\
\frac {\partial h_0}{\partial \psi_m}&...&\frac {\partial h_t}{\partial
\psi_m}\\[5pt]
a^ {(\ell)}_{1,0}&...&a^ {(\ell)}_{1,t} \\
... & ...& ... \\
a^ {(\ell)}_{t-m-1,0}&...&a^ {(\ell)}_{t-m-1,t} \\
\end{matrix}
\right)
$$
where ${\ell}=1,...,t-m$. Here $a^ {(\ell)}_{u,v}$ ($u=1,...,t-m-1,
v=0,...,t$) are elements of the base field $k$, while ${\partial
h_i}/{\partial \psi_j}$ stands for the derivative ${\partial
h_i}/{\partial y_j}$ computed at $y_j=\psi_j(x_{t+1},...,x_r)$, for
$i=0,...,t$ and $j=0,...,m$. Let $n$ denote the common degree of the
polynomials $Q_\ell$. Taking Laplace expansion along the first row,
one has an expression of the form:
$$Q_\ell=M_{\ell,0} x_0+...+M_{\ell,t}x_t$$
where $M_{\ell,i}$, ${\ell}=1,...,t-m, i=0,...,t$, are homogeneous
polynomials of degree $n-1$ in $x_{t+1},...,x_r$.

Fix an integer $d>n$ and set $\mu=[d/n]$. Fix biforms
$P_k(z_1,\ldots,z_{t-m};x_{t+1},\ldots,x_r)$ of bidegree $(k,d-kn)$,
$k=0,\ldots,\mu$. Finally, set

\begin{equation}\label{GNexpression}
f(x_0,...,x_r):=\sum_{k=0}^ \mu P_k(Q_1,\ldots Q_{t-m},x_{t+1},\ldots,x_r),
\end{equation}
a form of degree $d$ in $x_0,...,x_r$. It will be called a {\it
Gordan--Noether polynomial\/} (or a {\it {\rm GN}--polynomial\/}) of
type $(r,t,m,n)$,  and so will also any polynomial which can be
obtained from it by a projective change of coordinates. Accordingly,
a {\it Gordan--Noether hypersurface\/} (or {\it {\rm
GN}--hypersurface\/}) of type $(r,t,m,n)$ is the hypersurface $V(f)$
where $f$ is a non--zero GN--polynomial of type $(r,t,m,n)$.

\medskip

The main point of the Gordan--Noether construction is the following
result:

\begin{Proposition}\label {GNvanhess} Every {\rm GN}--polynomial has vanishing
Hessian.  \end{Proposition}

\begin{proof} Let $f(\mathbf x)$ be a GN-polynomial. Its first $t+1$ partial derivatives
$$f_i=\sum_{\ell=1}^ {t-m}
\frac {\partial f}{\partial Q_\ell}M_{\ell,i}, \quad i=0,...,t,$$
can be expressed in the form of a column vector

\begin{equation}\label{GNN}
\left ( \begin{matrix}
f_0\\
\ldots\\
f_t \\
\end{matrix}\right)
=\left ( \begin{matrix}
M_{1,0}&\ldots&M_{t-m,0} \\
\ldots & \ldots& \ldots \\
M_{1,t}&\ldots&M_{t-m,t} \\
\end{matrix} \right)
\cdot
 \left (\begin {matrix}
 \frac {\partial f}{\partial Q_1}\\
\ldots\\
\frac {\partial f}{\partial Q_{t-m}} \\
\end{matrix}\right).
\end{equation}

Consider the rational map
$$\phi_\ell: (x_{t+1}:\ldots :x_r)\in \p^ {r-t-1}\map
(M_{\ell,0}(x_{t+1},\ldots , x_r),\ldots,M_{\ell,t}(x_{t+1},\ldots,x_r))\in \p^
t$$ with $\ell=1,...,t-m$. Its image has dimension at most $m$ since
the polynomials $h_i, i=0,\ldots,t$, appearing in the determinants
which define the polynomials $M_{\ell,i}$, $\ell=1,...,t-m,
i=0,...,t$, depend on $m+1$ variables.

The rational map
$$\phi: (x_0:\ldots :x_r)\in \p^ r\map
(f_0(x_0,\ldots,x_r):\ldots :f_t(x_0,\ldots,x_r))\in \p^ t,$$ is the composite
of the polar map $\phi_f$ with the projection $(x_0:\ldots :x_r)\map
(x_0:\ldots :x_t)$. Therefore, if we let
$$\sigma: (x_0:\ldots :x_r)\in \p^ r\map (x_{t+1}:\ldots :x_r)\in \p^ {r-t-1}$$
denote the complementary  projection then, for a general point $p\in
\p^ r$,  equation \eqref{GNN} shows that $\phi(p)$ sits in the span
of $\phi_1(\sigma(p)),...,\phi_{t-m}(\sigma(p))$. Thus we see that
the image of $\phi$ has dimension at most $m+t-m-1=t-1$. This proves
that $f_0,...,f_t$ are algebraically dependent, hence so are
$f_0,...,f_r$.
\end{proof}

For a proof of the previous proposition which is closer to
Gordan--Noether's original approach, see \cite {Lossen}.

Following \cite {Permutti3} we give a geometric description of a
GN--hypersurface of type $(r,t,m,n)$, as follows. For this we
introduce the following notion.

\begin{Definition}\label{core} \rm Let $f$ be {\rm GN}--hypersurface of type $(r,t,m,n)$.
The {\it core\/} of $V(f)$ is the $t$-dimensional subspace
$\Pi\subset V(f)$ defined by the equations $x_{t+1}=...=x_r=0$.
\end{Definition}

We agree to call a GN--hypersurface of type $(r,t,m,n)$ {\it
general\/} if the defining data have been chosen generically, namely,
the polynomials $h_i(y_0,...,y_m)$, $i=0,...,t$, the polynomials
$\psi_j(x_{t+1},...,x_r)$, $j=0,...,m$, the constants $a^
{(\ell)}_{u,v}, {\ell}=1,...,t-m, u=1,...,t-m-1, v=0,...,t$, and
the biforms $P_k$, $k=0,\dots,\mu$,
are sufficiently general.

\begin{Proposition}\label{gn_hypersurface} Let  $V(f)\subset\p^r$ be a
{\rm GN}--hypersurface of type $(r,t,m,n)$ and degree $d$. Set
$\mu=[\frac dn]$. Then
\begin{enumerate}
\item[{\rm (i)}] $V(f)$ has multiplicity at least $d-\mu$ at the general point
of its core $\Pi\,${\rm ;}
\item[{\rm (ii)}] The general $(t+1)$-dimensional subspace $\Pi'$ through $\Pi$ cuts out on
$V(f)$, off $\Pi$, a cone of degree at most  $\mu$ whose vertex is an
$m$-dimensional subspace $\Gamma subset \Pi\,${\rm ;}
\item[{\rm (iii)}]  If $V(f)$ is general, then it has multiplicity exactly $d-\mu$ at the general point
of $\Pi$,  the general $(t+1)$-dimensional subspace $\Pi'$ through $\Pi$  cuts out on
$V(f)$, off $\Pi$, a cone of degree exactly  $\mu$, and,
as $\Pi'$ varies the corresponding subspace
$\Gamma$ describes the family of tangent spaces to an
$m$--dimensional unirational subvariety $S(f)$ of $\Pi\,${\rm ;}
\item [{\rm (iv)}] If $V(f)$ is general and $\mu>r-t-2$ then $V(f)$ is not a cone{\rm $\,$;}
\item [{\rm (v)}] The general {\rm GN}--hypersurface is
irreducible.
\end{enumerate}
\end{Proposition}
\begin{proof}
Let $\bar{\Pi}\subset \p^r$ denote the subspace defined by the
equations $x_0=...=x_t=0,$ the coordinate complementary subspace to
$\Pi$. For any non--zero $\xi=(0:\ldots: 0,\xi_{t+1}:\ldots:
\xi_r)\in \bar{\Pi}$, set $\Pi_{\xi}=<\Pi,\xi>\,\subset \p^r$. Then
as $\mathbf{\xi}$ varies, $\Pi_\xi$ describes the set of all
$(t+1)$--dimensional subspaces containing $\Pi$. For a fixed such
$\xi$ the points of $\Pi_\xi$ are parameterizable as $(x_0:\ldots
:x_t : z\xi_{t+1}:\ldots : z\xi_r)$, where $z$ is a parameter. Hence
we can take $(x_0:\ldots :x_t : z)$ as homogeneous coordinates in
$\Pi_\xi$ and $(x_0:\ldots :x_t)$ as coordinates in $\Pi$.

Fix such a $\xi$. The intersection $V(f)\cap \Pi_\xi$ is a
hypersurface $V_\xi$ of $\Pi_\xi$ with equation:

$$\sum_{k=0}^ \mu z^ {d-k}P_k\left(
M_{1,0}(\mathbf{\xi})\,x_0+...+ M_{1,t}(\mathbf{\xi})\,x_t,...,
M_{t-m,0}(\mathbf{\xi})\,x_0+...+ M_{t-m,t}(\mathbf{\xi})\,x_t,
\,\mathbf{\xi}\right)=0.$$

The presence of the factor $z^ {d-\mu}$ shows that the general point
of $\Pi$ has multiplicity at least $d-\mu$ for $V(f)$. This proves (i). The residual
hypersurface $W_\xi$ contains the subspace $\Gamma_\xi$ of $\Pi$
with equations:
\begin{equation}\label{vertice}
 M_{1,0}(\mathbf{\xi})\,x_0+...+
M_{1,t}(\mathbf{\xi})\,x_t=0,\,\ldots\,,
M_{t-m,0}(\mathbf{\xi})\,x_0+...+ M_{t-m,t}(\mathbf{\xi})\,x_t=0,
z=0.
\end{equation}
Furthermore $W_\xi$ is a cone with vertex $\Gamma_\xi$. Indeed, if
$p=(p_0:\ldots : p_t : p)\in W_\xi$ and $q=(q_0:\ldots :q_t:0)\in
\Gamma_\xi$, the line joining $p$ and $q$ is parameterizable by
$x_i=\lambda p_i+\nu q_i, \quad z= \lambda p,\quad  i=0,...,t$,
where $(\lambda :\nu)\in \p^ 1$ is a parameter. By restricting the
equation of $W_\xi$ to this line, we find that the resulting
equation is identically verified in $\lambda$ and $\nu$, because
$(q_0:\ldots :q_t)$ is a solution of the system \eqref {vertice}.
This proves (ii).

Assume now $V(f)$ is general.
We note that, by setting $x_0=1, x_2=\dots=x_t=0$, the coefficient
of $z^ {d-k}$, $k=0,\dots,\mu$, in the resulting polynomial can be seen
as a general polynomial of degree $d-k$ in the variables
$\xi_{t+1},...,\xi_r$.  By taking into account the proofs of parts (i) and (ii),
the first part of (iii) immediately follows.

Next note that $\dim(\Gamma_\xi)\geq m$ and the equality holds if
$V(f)$ is a general GN--hypersurface. Now $\Gamma_\xi$ contains the
$m+1$ points $p_j(\xi)=(p_{j,0}(\xi):\ldots :p_{j,t}(\xi))$, where
$$p_{j,i}(\xi)= \frac {\partial h_i}{\partial y_j}(\psi_0(\mathbf{\xi}),...,
\psi_m(\mathbf{\xi})), \quad i=0,...,t,  \quad j=0,...,m.$$
 Since $V(f)$ is general, the points $p_j(\xi), j=0,...,m$,
 are linearly independent. Hence $\Gamma_\xi=<p_0(\xi), ...,
p_m(\xi)>$. Consider the unirational subvariety $S(f)$ of $\Pi$
which is the image of the map $h:\bar \Pi\map \Pi$ sending the
general point $\xi\in \bar\Pi$ to the point $(\eta_0:\ldots :\eta_t)$
where
$$\eta_i=h_i(\psi_0(\mathbf{\xi}),...,
\psi_m(\mathbf{\xi})), \quad i=0,...,t.$$

It is clear now that $S(f)$ has dimension $m$ and that $\Gamma_\xi$ is the
tangent space to $S(f)$ at $h(\xi)$. This concludes the proof of (iii).

As for (iv), we notice that,  if $\mu+1>r-t-1$, then for no $\mathbf
\xi$ does the hypersurface $V_{\mathbf \xi}$ vanish identically.
Thus, if $V(f)$ is a cone, the vertex of $V(f)$ should lie on $\Pi$.
In this case all the tangent spaces to $S(f)$ should contain the
vertex of the cone, hence $S(f)$ itself ought to be a cone (cf. e.g.
\cite [Proposition 1.2.6]{Ru}). This is clearly not  the case for a
general GN--hypersurface.

To prove (v) let $f$ be a general GN--polynomial of type
$(r,t,m,n)$ and degree $d$ as in (\ref{GNexpression}) and let
$V(f)\subset\p^r$ be the corresponding hypersurface.
Cutting $V(f)$ with $\bar \Pi$ gives the hypersurface
with equation $P_0(x_{t+1},...,x_r)=0$, which does not involve the
variables $x_1,\ldots,x_{t-m}$. Indeed, for every $k=1,...,\mu$, the
homogeneous polynomial $P_k$ involves these variables and, being
general, must vanish for $x_0=...=x_t=0$ as do the $Q_\ell$'s. Now,
since $P_0(x_{t+1},...,x_r)$ is also general, the original
hypersurface is reduced. In addition, if $t<r-2$, the hypersurface
$V(P_0)$ is also irreducible, implying the irreducibility of the
original hypersurface. For $t=r-2$, the zero set of the polynomial
$P_0(x_{r-1},x_r)$ is a finite set of points, which are the
intersection points of the line $\bar \Pi$ with $V(f)$. However in this
case we can appeal to the fact that the polynomial
$P_0(x_{r-1},x_r)=0$ is a general equation of degree $d$ and
therefore its Galois group is the full symmetric group. Thus we see
that the intersection of $V(f)$ with a general line consists of $d$
distinct points, which are exchanged by monodromy when the line
moves, proving the irreducibility also in this case.
\end{proof}

The proposition admits a converse statement to the effect that if
$V(f)\subset\p^r$ is a hypersurface of degree $d$ satisfying a suitable
reformulation of the above enumerated properties, then it is a
GN--hypersurface of type $(r,t,m)$ (see \cite[pp. 104--105]{Permutti3}).

\subsection{Permutti's generalization of Gordan--Noether machine}

Permutti (see \cite {Permutti3}) has extended Gordan--Noether construction
in the case $t=m+1$. Let us briefly recall this too.

Fix integers $r,t$ such that $r\geq 2$, $1\leq t\leq r-2$. Fix $t+1$
homogeneous polynomials $M_0(x_{t+1},...,x_r),...,
M_t(x_{t+1},...,x_r)$ of the same degree $n-1$ in the variables
$x_{t+1},...,x_r$ and assume that they are algebraically dependent
over $k$ -- which will be automatic if $r\leq 2t$ because then the
number $r-t$ of variables is smaller than the number $t+1$ of
polynomials.

Set $Q=M_0x_0+...+ M_tx_t$, a form of degree $n$. Fix an integer
$d>n$ and set $\mu=[\frac dn]$. Further fix forms
$P_k(x_{t+1},...,x_r)$ of degree $d-kn$ in $x_{t+1},...,x_r$ ,
$k=0,...,\mu$. The form of degree $d$

$$f(x_0,...,x_r)=\sum_{k=0}^ \mu Q^ kP_k(x_{t+1},...,x_r),$$
or any form obtained thereof by a linear change of variables, will
be called a {\it Permutti polynomial}, or a
 {\it {\rm P}--polynomial\/} of type $(r,t,n)$. Accordingly, the corresponding
hypersurface $V(f)\subset\p^r$ will be called a {\it Permutti
hypersurface\/} or {\it {\rm P}--hypersurface\/} of type $(r,t,n)$,
with \emph{core}  the $t$--dimensional subspace $\Pi$ with equations
$x_{t+1}=\dots=x_r=0$. It is immediate to see that a GN--polynomial
of type $(r,t,t-1,n)$ is a P--polynomial of type $(r,t,n)$.

\begin{Proposition} \label{vanhessperm} Every {\rm P}--hypersurface has vanishing
Hessian.  \end{Proposition}

\begin{proof} Let $f$ be a P--polynomial. Then it is immediate to
see that

$$\frac {\partial f}{\partial x_i}=\frac {\partial f}{\partial Q}M_i, \quad
i=0,...,t,$$
 where $\partial f/\partial Q$ denotes the formal derivative of $f$ with respect to $Q$.
 Since by assumption $M_0,...,M_t$ are
algebraically dependent, it is clear that $\partial f/\partial x_0
,...,\partial f/\partial x_r$ are algebraically dependent too.
\end{proof}

One can easily prove the analogue of
Proposition~\ref{gn_hypersurface} in Permutti's setup. We use the
same terminology and notation employed in the previous section.

\begin{Proposition}\label{Pcarat} Let $V(f)\subset\p^r$ be a general {\rm P}--hypersurface
of type $(r,t,n)$ and degree $d$. Set $\mu= [\frac dn]$. Then
\begin{enumerate}
\item[{\rm (i)}] $V(f)$ has multiplicity  $d-\mu$ at the general point
of its core $\Pi \,$ {\rm ;}
\item[{\rm (ii)}]  The general $(t+1)$-dimensional subspace $\Pi'$ through
$\Pi$ cuts out on $V(f)$, off $\Pi$, a cone of degree at most $\mu$,
consisting of  $\mu$ subspaces of dimension $t$ which all pass
through a subspace $\Gamma$ of $\Pi'$ of dimension
$t-1\,${\rm ;}
\item[{\rm (iii)}]   As $\Pi'$ varies the corresponding $\Gamma$ describes a
unirational family of dimension $\chi\leq \min\{t-1,r-t-1\}\,${\rm
;}
\item [{\rm (iv)}]  If $\mu>r-t-2$, then $V(f)$ is a cone
if and only if the forms $M_0,...,M_t$ are linearly dependent over
$k$. This in turn happens as soon as either $t=1$, or $n=1,2${\rm ;}
\item [{\rm (v)}] $V(f)$ is
irreducible.
\end{enumerate}
\end{Proposition}

\begin{proof}
One verifies that the general
subspace $\Pi_\xi$ cuts out on $V(f)$ a hypersurface $V_\xi$ which
contains $\Pi$ with multiplicity $d-\mu$. The residual hypersurface
$W_\xi$ is the union of $\mu$ subspaces of dimension $t$ which all
pass through the subspace $\Gamma_\xi$ of $\Pi_\xi$ of dimension
$t-1$ with equation:

$$M_{0}(\mathbf{\xi})x_0+...+
M_{t}(\mathbf{\xi})x_t=0.$$

Note that, since $M_0,...,M_t$ are algebraically dependent,
then $\chi\leq t-1$. The inequality $\chi\leq r-t-1$ is
obvious. Parts (i)--(iii) follow by these considerations.

As for part (iv), like in the proof of Proposition~\ref{gn_hypersurface}, we
see that the hypersurface $V(f)$ is a cone if and only if, as $\xi$
varies, the subspace $\Gamma_\xi$ contains a fixed point. This
happens if and only if the polynomials $M_0,...,M_t$ are linearly
dependent. The rest of the assertion is trivial.

The proof of (v) is completely analogous to the proof of the
corresponding statement in Proposition \ref{gn_hypersurface} and
shall be omitted.
\end{proof}

It has been proved in \cite[pp. 100--101] {Permutti3} a converse to
the effect that if $V(f)\subset\p^r$ is a hypersurface of degree $d$ enjoying
the above properties -- with the core replaced by a subspace with
the same property -- then it is a P--hypersurface of type $(r,t,n)$.

For P--hypersurfaces $V(f)\subset\p^r$ one can describe the dual variety
$V(f)^*\subset\p^{r*}$. Note that, as $\xi$ varies in the subspace $\bar \Pi$
with equations $x_0=\dots=x_t=0$, then the subspace
$\Gamma_\xi^ \perp$ of dimension $r-t$ varies
describing a cone $W(f)\subset\p^{r*}$ , of
dimension $r-t-1$ with vertex $\Pi^ \perp$ which contains
the subspace $\Pi_\xi^ \perp$ of dimension $r-t-2$.
More precisely, we have the:

\begin{Proposition}\label{Pdual} Let $V(f)\subset\p^r$ be a general {\rm P}--hypersurface
of type $(r,t,n)$ and degree $d$. Let $\mu=[\frac dn]$. Then:
\begin{enumerate}
\item[{\rm (i)}]  $V(f)^*\subset W(f)$, where $W(f)\subset\p^{r*}$ is a cone over a
unirational variety of dimension $\chi\leq \min\{t-1,r-t-1\}$ whose
vertex is the orthogonal of the core $\Pi$ of $V(f)${\rm ;}
\item[{\rm (ii)}]  The general ruling
of the cone $W(f)\subset\p^{r*}$ is an $(r-t)$-dimensional subspace through $\Pi^
\perp$ which cuts $V(f)^ *$, off $\Pi^ \perp$, in $\mu$ subspaces of
dimension $r-t-1$  all passing through the same subspace of $\Pi^
\perp$ of dimension $r-t-2$. Hence $v(f)=\min\{r-2,2(r-t-1)\}$.
\end{enumerate}
Conversely, if $V(f)\subset\p^r$ is the dual of such a variety, then
$V(f)\subset\p^r$ is a {\rm P}--hypersurface.
\end{Proposition}
\begin{proof} It follows by dualizing the contents of Proposition \ref{Pcarat}.
\end{proof}

From this we also see that a general {\rm P}-hypersurface is not a
cone. In addition, one has:

\begin{Proposition}\label{Pzeta} Let $V(f)\subset\p^r$ be a general {\rm P}--hypersurface
of type $(r,t,n)$. Then  $Z(f)=W(f)\subset\p^{r*}$, and therefore
$z(f)=\min\{r-1,2(r-t)-1\}$.
\end{Proposition}

\begin{proof} For $\xi\in \bar\Pi$ general, $\Pi_\xi$ cuts out on $V(f)$ a
hypersurface $V_\xi$ which is a
union of hyperplanes of $\Pi_\xi$ and is a
cone with vertex $\Gamma_\xi$. If $p\in \Pi_\xi$ is a general point, then
the polar hyperplane $\pi_{\xi,p}$ of $p$ with respect to $V_\xi$ contains $\Gamma_\xi$.
By Remark \ref {arrangm}, when $p$ varies in $\Pi_\xi$, then $\pi_{\xi,p}$ varies
describing an open dense subset of the set of
all hyperplanes of $\Pi_\xi$ containing
$\Gamma_\xi$. If $\pi_p(f)$ is the polar hyperplane of $p$ with respect to $V(f)$,
then $\pi_p$ cuts out $\pi_{\xi,p}$ on $\Pi_\xi$. Hence the subspace
$<\phi_f(p),\Pi_\xi^ \perp>$ sits in the ruling $\Gamma_\xi^ \perp$ of $W(f)$ and,
as $p$ varies, it describes a dense open subset of
$\Gamma_\xi^ \perp$. This proves that $W(f)=Z(f)$.
\end{proof}

\begin{Remark} {\rm The case $t=r-2$ is particularly interesting. Then
 $V(f)^ *$ is a scroll surface with a line
directrix $L=\Pi^ \perp$ of multiplicity $e\geq \mu$, where $\mu$ is
the invariant introduced in {Section} \ref{multidir}. It is a
subvariety of the $3$--dimensional rational cone $W(f)$  over a
curve with vertex $L$, and the general plane ruling of the cone cuts
$V(f)^*$ along $\mu$ lines of $V(f)^ *$, all passing through the
same point of $L$.  In particular, for $\mu=1$, the dual $V(f)^ *$
is a rational scroll  (see {Sections}
\ref{RNS} and \ref{multidir}). According to Proposition \ref
{Pzeta}, we have $Z(f)=W(f)$, hence $z(f)=3$.}\end{Remark}

If $t=2$ the two constructs of GN--hypersurfaces and
P--hypersurfaces coincide. For $r=4$ this is the only value of $t$
which leads to hypersurfaces which are not cones. The case $r=4$ is
well understood due to a result of Franchetta (see \cite {Fr2}; see also
Proposition~\ref{Pdual}; according to \cite {Lossen}  this result
is contained in \cite {GN}; for another proof see \cite {GR}):

\begin{Theorem}\label {frank} Let $V(f)\subset \p^4$ be a reduced hypersurfaces
of degree $d$. The following conditions are equivalent:
\begin{enumerate}
\item[{\rm (i)}]  $V(f)$ has vanishing Hessian {\rm ;}
\item[{\rm (ii)}]  $V(f)$ is a {\rm GN}--hypersurface of
type $(4,2,1,n)$, with $\mu=[\frac dn]$, which has a plane of
multiplicity $d-\mu\,${\rm ;}
\item[{\rm (iii)}]  $V(f)^* $ is a
scroll surface of degree $d$, having a line directrix $L$ of
multiplicity $e$, sitting in a rational cone $W(f)$ of dimension $3$
with vertex $L$, and the general plane ruling of the cone cuts
$V(f)^ *$ off $L$ along $\mu\leq e$ lines of the scroll, all passing
through the same point of $L$.
\end{enumerate}
In particular, $V(f)^*$ is smooth if and only if $d=3$,
$V(f)^*$ is a rational normal scroll and $V(f)$ contains a plane, the
orthogonal to the line directrix of $V(f)^ *$, with multiplicity $2$.
\end{Theorem}

\subsection{Variations on some results of Perazzo}\label {Peraz}

Let $V(f)\subset\p^r$ be a hypersurface of degree $d$ with $r\geq
4$. If $d=2$, it is clear that $V(f)$ has vanishing Hessian if and
only if it is a cone. So the first meaningful case is the one $d=3$,
in which, as we saw, there are examples which are not cones (see
Theorem~\ref{frank}). The case of cubic hypersurfaces has been
studied in some detail by U. Perazzo (see \cite{Perazzo}). We will
partly generalize Perazzo's results. Inspired by the construction of
P--hypersurfaces and by Perazzo's results, we will give new examples
of hypersurfaces with vanishing Hessian, which are extensions of
some P--hypersurfaces.

Consider a hypersurface $V(f)\subset\Proj^ r$ which contains a
subspace $\Pi$ of dimension $t$ such that the general subspace
$\Pi_\xi$ of dimension $t+1$ through $\Pi$ cuts out on $V(f)$ a cone
with a vertex $\Gamma_\xi$ of dimension $s$. Assume that $s\geq r-t-1$.
By extended analogy, we will call $\Pi$ the {\it core\/} of $V(f)$
and call $V(f)$ an H--{\it hypersurface\/} of type $(r,t,s)$. Notice
that a P--hypersurface of type $(r,t,n)$ with $r\leq 2t$ is also an
H--hypersurface of type $(r,t,t-1)$.

As for P--hypersurfaces, we can introduce the cone $W(f)\subset\p^{r*}$ with vertex
$\Pi^ \perp$, which is swept out by the $(r-s-1)$--dimensional subspaces $\Gamma_\xi^ \perp$
as $\Pi_\xi$  varies among all subspaces of dimension $t+1$ containing
$\Pi$.

A special case of an H--hypersurface is that of a hypersurface
$V(f)\subset\p^r$ of degree $d$ containing a subspace  $\Pi$ of dimension $t$
whose general point has multiplicity $d-\mu>0$ for $V(f)$, such that
the general subspace $\Pi_\xi$ of dimension $t+1$ through $\Pi$ cuts
out on $V(f)$, off $\Pi$, a union of $\mu$ subspaces of dimension
$t$, with $\mu\leq 2t-r+1$.  In this situation, we will call $V(f)\subset\p^r$
an R--{\it hypersurface\/} of type $(r,t,\mu)$.

\begin{Proposition}\label{Dhyper} An {\rm H}--hypersurface $V(f)\subset\p^r$ of type
$(r,t,s)$ has vanishing Hessian. Moreover $Z(f)=W(f)\subset\p^{r*}$.
\end{Proposition}

\begin{proof} Let $p$ be a general point in $\Proj^ r$ and let $\Pi'$ be
the span of $\Pi$ and $p$. Since the intersection of $V(f)$ with $\Pi'$
is a cone with vertex a subspace $\Gamma$ of dimension $s$, the polar
quadric $Q_p(f)$ cuts out on $\Pi'$ a quadric singular along $\Gamma$.
If $Q_p(f)$ is smooth we have $s=\dim(\Gamma)\leq r-t-2$, a
contradiction. This proves that $Q_p(f)$ is singular hence $V(f)$ has
vanishing Hessian.

The argument for the second assertion is similar to the one in the
proof of Proposition \ref {Pzeta} and therefore can be omitted. \end{proof}

\begin{Remark} \label {dualrig} It is interesting to look at duals of
R--hypersurfaces of degree $d$ and type $(r,r-2,\mu)$. If
$V(f)\subset\p^r$ is such a hypersurface, its dual $V(f)^
*\subset\p^{r*}$ is a scroll surface with a line directrix $L$ of
multiplicity $e\geq \mu$,  where $\mu\leq
r-3$ is as in { Section} \ref {multidir}. We assume
$V(f)\subset\p^r$ not a cone and therefore $V(f)^ *\subset\p^{r*}$
is non--degenerate.

In this case the invariant $s$ is related to the number $\nu$
introduced in  {Section}
\ref{multidir}: $\nu=r-s-1$ {and, moreover, one has
$Z=W(f)\subset\p^{r*}$ where $Z$ is the cone in the same section,
and} $\dim W(f)=r-s$.

By  Proposition~\ref {Dhyper}, one has  $Z(f)=W(f)=Z$. This means
that $\rho(f)=r-s+1$,  hence the vertex of the general polar quadric
has dimension $s-1$.

Let $p\in \p^ r$ be a general point. The quadric $Q_p(f)$ cuts the
hyperplane $\Pi'=<\Pi,p>$ in a quadric singular along the subspace
$\Gamma$ of dimension $s$.  Set $\xi=\phi_f(p)$. The vertex
of $Q_p(f)$, which coincides with $T_{Z(f),\xi}^ \perp$
(see Remark \ref {hessim}), has dimension $s-1$,
hence it is  contained in $\Gamma$.  \end{Remark}

An R--hypersurface with $\mu=1$ is a hypersurface of degree $d$ with
a core  $\Pi$ of dimension $t$ whose general point has multiplicity
$d-1$ for the hypersurface, and moreover $2t\geq r$. This is the
case considered by Perazzo in \cite[p. 343]{Perazzo}, where he
proves that these hypersurfaces have vanishing Hessian.

\section{Homaloidal polynomials}\label {homal}

A hypersurface $V(f)\subset \p^ r$, or the form $f$, of degree $d$
is said to be {\it homaloidal} if $\delta(f)=1$, i.e. if the polar
map $\phi_f$ is birational. According to Theorem \ref {nored}, this
property depends only on ${\rm Supp}(V(f))$, therefore we will
mainly refer to the case $V(f)$ reduced.

The simplest example is when $V(f)$ is a smooth quadric: in this
case the polar map $\phi_f$ is the usual polarity, which is an
invertible linear map. This is also the only case of a reduced
homaloidal polynomial if $r=1$.

Reduced homaloidal curves in $\p^ 2$ have been classified by
Dolgachev in \cite {Dolgachev}:

\begin{Theorem}\label {dolgclass}  A reduced plane curve $V(f)\subset \p^ 2$
of degree $d$ is homaloidal if and only if either
\begin {itemize}
\item[{\rm (i)}]  $V(f)$ is a smooth conic, or
\item[{\rm (ii)}] $d=3$ and $V(f)$ consists of three non concurrent lines, or
\item[{\rm (iii)}] $d=3$ and $V(f)$ consists of the union of a smooth conic
with one of its tangent lines.
\end {itemize}
\end{Theorem}

Note that in case (ii) the polar map $\phi_f$ is a {\it standard
quadratic transformation\/} based at three distinct points, whereas
in case (iii) the map $\phi_f$ is a {\it special quadratic
transformation\/} based at a curvilinear scheme of length three
supported at one single point. More algebraically, in cases (ii) and
(iii) the base locus ideal of $\phi_f$ is a codimension $2$ perfect
ideal (Hilbert--Burch) - see \ref{hankelmatrices} and also
\cite{cremona}, \cite{Simis2} for the ubiquitous role of
Hilbert--Burch ideals in the theory of Cremona transformations.

\begin{Remark}\label{ext} We note that the three cases in Theorem~\ref{dolgclass}
can be naturally extended to any dimension $r\geq 2$, thus yielding
an infinite series of homaloidal hypersurfaces in $\p^ r$, with
$r\geq 2$ (see  \cite {Dolgachev}). Namely, the following reduced hypersurfaces
$V(f)\subset\p^r$ of degree $d$ are homaloidal in $\p^ r$ for any
$r\geq 2$:
\begin {itemize}
\item[{\rm (i)}]  A smooth quadric{\rm ;}
\item[{\rm (ii)}] The union of $r+1$ independent hyperplanes{\rm ;}
\item[{\rm (iii)}] The union of a smooth quadric with
one of its tangent hyperplanes.
\end {itemize}
Note that (ii) gives the only example of arrangements of hyperplanes
which are homaloidal (see \cite {DP} , \cite {Dolgachev}).
\end{Remark}
\medskip

There is a general principle for rational maps
$\phi:\p^r\dasharrow\p^r$. In what follows, we adopt the terminology
{\it the image of\/}  $\phi$ to mean the closure in the target of
the image of the points of the source $\p^r$ at which $\phi$ is
well-defined. Accordingly, we use the notation $\phi(\p^r)$. This
convention sticks to subvarieties as well.

\begin{Proposition}\label{jacobian}
Let $\phi=(F_0:\cdots :F_r):\p^r\dasharrow \p^r$ denote a rational
map, where $F_i\in k[x_0,\ldots,x_r]$ are forms of the same degree
without proper common factor. Let $J\in k[x_0,\ldots,x_r]$ denote
the Jacobian determinant of these forms. Consider the following
conditions:
\begin {enumerate}
\item[{\rm (i)}]  $J\neq 0${\rm ;}
\item[{\rm (ii)}] $\dim (\phi(\p^r))=r$.
\end {enumerate}
Then {\rm (i)} $\Leftrightarrow$ {\rm (ii)}.

 If $\phi$ is birational, then $\dim \phi(V(J))\leq r-2$.
\end{Proposition}
\begin{proof}
  First note that $\phi$ is well-defined at a general point of
$V(J)$, otherwise $F_0,\ldots ,F_r$ would be multiples of a single
form which contradicts the assumption on these forms.

(i) $\Leftrightarrow$ (ii) Note that, up to a degree
renormalization, the homogeneous coordinate ring of
$\phi(\p^r)\subset \p^r$ is $k[F_0,\ldots,F_r]$. One then draws upon
the known fact saying that, in characteristic zero, the dimension of
$k[F_0,\ldots,F_r]$ is the rank of the Jacobian matrix of
$F_0,\ldots ,F_r$ (see, e.g., \cite{Simis}).

If $\phi$ is birational, it is dominant so that $J\neq 0$. Moreover,
it has to contract the hypersurface $V(J)$
since this is the locus where $\phi$ drops rank.
\end{proof}

\medskip

If $\dim(V(J))\leq r-2$, we shall say that $V(J)$ is {\it contracted
by\/} $\phi$.
\medskip

\begin{Corollary}\label{shom} If a hypersurface $V(f)\subset \p^ r$
is homaloidal,  then $h(f)$ does not vanish identically and $H(f)$ is contracted by the polar map
$\phi_f$.
\end{Corollary}
\medskip

\begin{Remark}\label{tothess} An interesting case of Corollary~\ref{shom}
is when $f\in k[x_0,...,x_r]$ is a non--zero reduced, homaloidal  form of degree $d$ such
that $h(f)=cf^{\frac {(d-2)(r+1)}d}$ with $c\in k^*$.
In this case we will say that such an $f$ is {\it totally Hessian\/}
and use the same terminology for the corresponding hypersurface.
Note that it entails the equality ${\rm Supp}(V(f))={\rm
Supp}(H(f))$ -- hence $V(f)$ is also contracted by $\phi_f$ -- and
any smooth point of $V(f)$ is parabolic (see Remark~\ref{hesspol}).
 It would be interesting to find whether a totally Hessian form is
 homaloidal.

 A good deal of examples of totally Hessian forms
arises from the theory of pre-homogeneous vector spaces (see
\cite{Segre}, \cite{Segre2}, also \cite{Mukai}), a notion introduced
by Kimura and Sato (see \cite{KS}, see also \cite {ESB}, \cite
{Dolgachev}, \cite {EKP}), which we now briefly recall for the
reader's convenience.

A {\it pre--homogeneous vector space\/} is a triple $(V,G,\chi)$
where $V$ is a complex vector space of finite dimension, $G$ is a
complex algebraic group, $V$ is an algebraic linear representation
of $G$, $\chi: G\to \C^*$ is a non--trivial character, and there is
a non--zero homogeneous polynomial $f: V\to \C$, with no multiple
factors, such that $f(g\cdot v)=\chi(g)f(v)$ for all $g\in G$ and $v\in
V$, and such that the complement of the hypersurface $\{f=0\}$ is a
$G$--orbit.

The polynomial $f$, called the {\it relative invariant\/} of the
pre--homogeneous space, is unique up to a non--zero factor from $\C$.
The pre--homogeneous vector space  $(V,G,\chi)$ is said to be {\it
regular\/} if $h(f)\neq 0$. In this case the relative invariant $f$
is totally Hessian (see \cite{KS} and \cite{ESB}) and $f$ is a
homaloidal polynomial such that $\phi_f$ coincides with its inverse,
modulo a projective transformation (see \cite[Theorem 2.8]{ESB}). In
\cite{EKP} (see also in \cite{Mukai}, \cite[Ch. III]{Zak1} and
\cite{ESB}), there is a description of several regular homogeneous
vector spaces related to smooth projective varieties with extremal
geometric properties (Severi and Scorza varieties, some varieties
with one apparent double point, varieties whose dual is small, see
{\it loc. cit.\/}). The first instances among these examples were
described in the classic literature (see \cite{Cayley},
\cite{Segre}, \cite{Fr1}).
\end{Remark}
\medskip

Being homaloidal or having vanishing Hessian implies strong
constraints on the singularities of the hypersurface $V(f)$.  Thus,
if $\dim V(f)\geq 2$ and if $V(f)\subset\p^r$ has vanishing Hessian, then $V(f)$ cannot have
isolated singularities. Also there is a conjecture in \cite{DP} to
the effect that a hypersurface of dimension at least $2$ with isolated
singularities cannot be homaloidal. We now prove a result which
points somewhat in this direction.

First we need to introduce some notation. Suppose $V(f)\subset\p^r$ is a
reduced hypersurface  of degree $d$. Let us resolve the
indeterminacies of the polar map $\phi_f$ by iteratively blowing up
$\p^ r$
$$X:=X_n\to X_{n-1}\to... \to X_1\to X_0=\p^ r,$$
thus getting $p: X\to \p^ r$ so that $\phi_f\circ p: X\to {\p^ r}^
*$ is a morphism. Here the map $X_i\to X_{i-1}, i=1,...,n,$ is a blowup with
center a smooth variety of codimension $a_i+1$, with $1\leq a_i\leq
r-1$. We denote by $E_i$ the total transform on $X$ of the
exceptional divisors of the blowup $X_i\to X_{i-1}, i=1,...,n$.
Further let $H$ stand for the proper transform on $X$ of a general
hyperplane of $\p^r$ and $\Phi$ for the proper transform on $X$ of
the first polar hypersurface of $V(f)$ with respect to a general
point of $\p^r$. Then
$$\Phi\equiv (d-1)H-\sum_{i=1}^ n \mu_iE_i,$$
where the $\mu_i$ 's are the {\it multiplicities\/} of $\Phi$ along
the various centers of the iterated blowups. By an obvious
minimality assumption, we may assume $\mu_i>0, i=1,...,n$.

\medskip

The following result can be seen as a
consequence of the so--called \emph{Noether--Fano
inequality} for Mori fibre spaces (see \cite {Corti}). We give here
a short direct proof. Let us recall that $\delta(f)=\deg(\phi_f)$
with the usual convention that $\deg(\phi_f)=0$ if and only if
$\phi_f$ is not dominant.

\begin{Proposition}\label {codimen} In the above setting, if
$\delta(f)\leq 1$ then either $d\leq r+1$ or $\mu_i>a_i$ for some
$i=1,...,n$, i.e. either $d\leq r+1$ or  the singularities of the
general first polar of $V(f)$ are not \emph {log--canonical}
{\rm (see \cite {KM}, p. 56)}.
\end{Proposition}

\begin{proof} As above, let $\Phi$ denote the proper transform
on $X$ of the general first polar hypersurface of $V(f)$. Note that
$\Phi$ is smooth, because the linear system $|\Phi|$ is base point
free. If $\delta(f)\leq 1$  then $\Phi$ is either rational or ruled
(see Proposition~\ref{hessian}). Since
$$K_X\equiv -(r+1)H+\sum_{i=1}^ n a_iE_i,$$
one has
$$K_\Phi\equiv (d-r-2)H_{\big| \Phi}+\sum_{i=1}^ n (a_i-\mu_i)E_{i\big| \Phi}.$$

If $d>r+1$ and $a_i\geq \mu_i$ for every $i=1,...,n$, this divisor
is effective, which would contradict the ruledness of $\Phi$.
\end{proof}
\medskip

\begin{Remark}\label{eeff} Although the proper transform $\tilde{V}$ of $V(f)$ on
$X$ admits a similar expression
$$\tilde{V}\equiv dH-\sum_{i=1}^ n m_iE_i,$$
here, in spite of the previous minimality assumption, some of the
$m_i$'s may vanish (see Section~\ref{sing} below).
\end{Remark}

 The problem of understanding the relationship between the
$m_i$'s and the $\mu_i$'s is longstanding, dating back to M.
Noether, and is far from being solved in general. For further
contributions in the plane case see \cite{segp} and \cite{ves} (see
also Remark~\ref{effbeh} below). Roughly speaking, one would expect
$\mu_i=m_i-1$ but this is not always the case.

\begin{Corollary} In the above setting suppose that
 $\mu_i=m_i-1$ for all $i=1,...,n$. If $d\geq r+2$ and $\delta(f)\leq 1$, then $m_i>a_i+1$ for some $i=1,\ldots, n$.
In particular, a surface $V(f)\subset \p^ 3$ with $d\geq 5$ and
$\delta(f)\leq 1$ cannot have {\it ordinary\/} singularities.
\end{Corollary}

\subsection{Irreducible homaloidal polynomials of arbitrarily large degrees}\label
{infser}

In this section we produce, for every  $r\geq 3$, an infinite series
of irreducible homaloidal hypersurfaces in $\p^ r$ of arbitrarily
large degree, thus settling a question that has been going around
for some time. These polynomials are the dual hypersurfaces to
certain scroll surfaces. It is relevant to observe, as we indicate
below, that these examples are not related to the ones based on
pre--homogeneous vector spaces as in \cite{EKP} and in \cite{ESB}.

The examples show that, perhaps opposite to the ongoing folklore,
there are plenty of homaloidal polynomials around. They even seem to
be in majority as compared to polynomials with vanishing Hessian,
though a complete classification does not seem to be presently at
hand.

In this respect Dolgachev's classification Theorem~\ref{dolgclass}
might be considered in counterpoint to Hesse's result to the effect
that the only hypersurfaces with vanishing Hessian in $\p^ r$,
$r\leq 3$,  are cones (see { Section}
\ref {hess}).

We wonder whether a counterpart of Franchetta's Theorem~\ref{frank}
could be a result to the effect that in $\p^ 3$ there are only
finitely many projectively distinct types of (irreducible)
homaloidal polynomials, apart from the ones constructed in this
section.

\medskip

We start with lemmas of general content.

\begin{Lemma}\label{lem:ess} Let $V(f)$ be a hypersurface in $\p^ r$.
Suppose there is  a point  $p\in V(f)$ and $s$ linearly independent hyperplanes
$H_i$, $i=1,\dots,s$, passing through $p$ and each cutting $V(f)$ in a hypersurface
having a point of multiplicity at least $s$ in $p$. Then $V(f)$ has multiplicity at
least $s$ in $p$.
\end{Lemma}

\begin{proof} Assume $p$ is the origin in affine coordinates and that $H_i$ has
equation $x_i=0$, $i=1,\dots,s$. Write $f=f_0+f_1+\dots+f_d$, where $f_j$
is the homogeneous component of degree $j$, with $j=0,\dots,d$.
By the assumption $f_0,\dots,f_{s-1}$ have to
be divisible by $x_i$, $i=1,\dots,s$. Hence $f_0,\dots,f_{s-1}$ are identically zero,
proving the assertion.\end{proof}

Recall now that the polar map
of a form $f\in k[x_0,\ldots, x_r]$ is denoted $\phi_f$ or
$\phi_{V(f)}$ to stress the corresponding hypersurface $V(f)\subset
\p^r$.

\begin{Lemma}\label {projpoint} Let $V(f)\subset \p^ r$ be a
hypersurface. Let $H\subset \p^ r$ be a hyperplane  not contained in
$V(f)$, let $\xi=H^ \perp$ be the corresponding point in ${\p^ r}^*$
and let $\sigma_{\xi}$ denote the projection from $\xi$. Then
$$\phi_{V(f)\cap H}= \sigma_\xi\circ ({\phi}_{V(f)})_{|H}.$$
\end{Lemma}

\begin{proof} The proof is
straightforward by assuming, as one can, that $H$ is a coordinate
hyperplane.\end{proof}

\begin{Lemma}\label{ratcurve} Let $C$ be  a rational normal curve in
$\p^ n$. Let $\L$ be a $g^ 1_m$, with $m<n$ and consider the
rational normal scroll $W(\L)=\cup_{D\in \L}<D>$, of dimension $m$
and degree $n-m+1$ {\rm (}see, e.g., {\rm \cite{EH})}. Let
$p_1,...,p_{n-m}$ be points of $C$. Then:
\begin{itemize}
\item[{\rm (i)}] $W(\L)$ intersects the
$(n-m-1)$--dimensional subspace $\Pi=<p_1,...,p_{n-m}>$
transversally only at $p_1,...,p_{n-m}${\rm ;}
\item[{\rm (ii)}] the general tangent space to $W(\L)$ does not intersect the
$(n-m-1)$--dimensional subspace $\Pi=<p_1,...,p_{n-m}>$.
\end{itemize}
\end{Lemma}

\begin{proof}  Project from $\Pi$ to $\p^ {m}$.
The image of $C$ is a rational normal curve $C'$ and the image of
$W(\L)$ is the analogous scroll $W(\L')=\cup_{D\in \L'}<D>$, where
$\L'$ is the  $g^ 1_m$ on $C'$ which is the image of $\L$. Since
this fills up $\p^ {m}$, both assertions follow.
\end{proof}

We introduce now the promised examples. Recall
 from { Section} \ref{multidir} that we have
rational scrolls $X:=X(a,b)$ of degree $d=a+b$ with a multiple line
directrix $\Lambda$ of multiplicity $a$ in $\p^ {b+2}$, for $1\leq
a\leq b$. The dual hypersurface $X^
*=X(a,b)^ *$ has vanishing Hessian as soon as $1\leq a<b$ and the
image of the corresponding polar map is the cone $Z:=Z(a,b)$
containing $X$,   introduced in {
Section}  \ref{multidir} (see Remark \ref{dualrig}). This is a
rational normal cone of degree $b-a+1$ and dimension $a+2$ with
vertex $\Lambda$. More specifically, let $C$ be the rational normal
curve in $\p^ b$ which is the projection of $X$ from $\Lambda$. One
has the general linear series $\L=g_a^ 1$ on $C$ whose general
divisor is the projection on $C$ of the $a$ lines of $X$ passing
through the general point of $\Lambda$. The scroll $Z$ is the cone
with vertex $\Lambda$ over $W(\L)$, which, by  the generality
assumption about $X(a,b)$, is a general rational normal scroll of
degree
 $b-a+1$ and dimension $a$ in $\p^ b$.
 \medskip

An essential piece of information for the construction of our examples
is the following:

\begin{Theorem}\label{fibrespace} If $1\leq a<b$, the closure
of the general fibre of the polar map $\phi:=\phi_{X(a,b)^ *}$ is a
projective subspace of dimension $b-a$ of $\,\p^ {b+2}$.
\end {Theorem}

\begin{proof} Let $p$ be a general point of $\p^ {b+2}$. Then $\xi=\phi(p)$
is a general point of $Z$. Recall that the closure $F_p$ of the fibre of $\phi$
over $\xi$ is the union of finitely many $(b-a)$--dimensional subspaces
containing $T_{Z,\xi}^ \perp$, which in turn is the $(b-a-1)$--dimensional
vertex $V_p$ of the polar quadric
$Q_p$ of $p$ with respect to $X^*=X(a,b)^ *$ (see Proposition
\ref {hessian}  and Remark \ref {hessim}).  What we have to prove is that
$F_p$ consists of the single subspace $\langle p,V_p\rangle$.

Recall Proposition \ref {directrix} and
Remark \ref  {dualrig} and keep the notation introduced
therein. In particular $\Pi=\Lambda^ \perp$ is a subspace of dimension $b$
in ${\p^ {b+2}}^ *$, which has multiplicity $b$ for $X^ *$.
The hyperplane $\Pi':=\Pi'_p=<\Pi,p>$ is dual to the general point $x\in \Lambda$.
Let $F_{x,1},\dots, F_{x,a}$ be the rulings of $X$ passing through $x$, hence
$\Pi'$ cuts $X^ *$ along $\Pi$, with
multiplicity $b$, and along the $a$ subspaces $\Sigma_i:=F_{x,i}^ \perp$,
$i=1,\dots,a$, of dimension $b$.
The intersection

$$\Gamma_p=\Pi \cap  \Sigma_1\cap...\cap
\Sigma_a=<\Lambda, F_{x,1},...,F_{x,a}>^ \perp$$
has dimension $b-a$. Note that $\Gamma_p=W_\xi^ \perp$, where
$W_\xi=<\Lambda, F_{x,1},...,F_{x,a}>$ is the ruling of $Z$ containing $\xi$.

Let $p'$ be another point in $\p^ {b+2}$ where $\phi$ is defined, and
set $\xi'=\phi(p')$. The above description implies that $\Gamma_p=\Gamma_{p'}$
if and only if $W_\xi=W_{\xi'}$. By recalling the structure of the
scroll $Z$ we see that this happens if and only if $\Pi_p'=\Pi'_{p'}$.

As we saw in Remark \ref  {dualrig},
the vertex $V_p$ of the quadric $Q_p$ is contained in $\Gamma_p$ because $W_\xi\subseteq T_{Z,\xi}$.
We claim now that there is no point $p'$ such that $\Gamma_{p'}\neq
\Gamma_p$ and $\Gamma_p\cap\Gamma_{p'}=V_p$. In fact
if this happens, then $T_{Z,\xi}=V_p^ \perp$ contains $W_{\xi'}=\Gamma_{p'}^ \perp$.
This means that, if $W$ is a general ruling of $Z$, then
 the tangent space to $Z$ at the general point of $W$
contains some other ruling $W'$ of $Z$. By projecting $Z$ from $\Lambda$ onto
the $a$--dimensional rational normal scroll $W(\L)\subset\p^b$,
we would have that, for a general point $q\in W(\L)$,
the tangent space $T_{W(\L),q}$  would contain
some ruling of $W(\L)$ different from the one of $q$.
This is impossible. Indeed, by cutting with $a-1$ general
hyperplanes, we would have
the general curve section $C$ of $W(\L)$,
a rational normal curve, with
the property that its general tangent line $T_{C,q}$ intersects $C$ at a point
$q'\neq q$, which is clearly not the case.

Let now $p,p'$ be points such that $\phi(p)=\phi(p')$. Then $V_p=V_{p'}$,
therefore $\Gamma_p=\Gamma_{p'}$
and $\Pi_p'=\Pi'_{p'}$.  Thus $F_p=F_{p'}$
is contained in $\Pi_p'$. To simplify notation,
we set $V=V_p$, $\Gamma=\Gamma_p$,
$\Pi'=\Pi'_p$, $F=F_p$.

We claim next that if  $\phi(p)=\phi(p')$,
then $\langle p,\Gamma\rangle=\langle p',\Gamma\rangle$
and this $(b-a+1)$--dimensional
subspace $\Gamma'$ contains $F$.
In fact, $F$ is the closure of the intersection, off the singular
locus of $X^ *$, of all
first polars of $X^ *$ containing $p$. In particular $F$ is contained
in the intersection of $\Pi'$ with of all first polars  of points of $\Pi'$
containing $p$ (or $p'$).
Remarks \ref {arrangm} and \ref {ext}, imply that this intersection
is exactly $\Gamma'$. Our claim thus follows.

Note now that the linear system cut out on $\Gamma'$ by of all first polars of the
points of $\Pi'$ is $0$--dimensional, consisting of $\Gamma$, counted
with multiplicity $a+b-1$. Thus
the linear system $\NN$ of hypersurfaces of degree $a+b-1$
cut out on $\Gamma'$ by all first polars
of $X^ *$ is a pencil, i.e.  $\dim(\NN)=1$.  Note that the fixed locus of
$\NN$ certainly contains $\Gamma$ with multiplicity $b-1$,
since the general first polar contains $\Pi$ with this multiplicity.
To finish our proof, we have to show that the movable part
of $\NN$, whose degree is bounded by $a$, is actually a pencil of hyperplanes.

To see this,  look at the linear system $\MM$ cut out by all first polars
on $\Pi'$ off $\Pi$, which, as we said, appears with multiplicity $b-1$ in the base locus.
The general member $M$ of $\MM$ is a hypersurface of degree $a$. Let us
consider its intersection with the hyperplanes $\Sigma_i$, $i=1,\dots,a$. Note that the
intersections $\Sigma_i\cap\Sigma_j$, $1\leq i<j\leq a$, all of dimension $b-1$ and
containing $\Gamma$, sit in the singular locus of $X^ *$, since they are intersection
of rulings of the scroll $X^ *$. Hence the
intersection of $M$ with $\Sigma_i$ has multiplicity $a-1$ along $\Gamma$
for all $i=1,\dots,a$. By Lemma \ref {lem:ess}, $M$ has multiplicity $a-1$ along
$\Gamma$. This implies  that the movable part
of $\NN$ has degree one, thus ending the proof of the theorem.
\end{proof}

Let now $F_1,...,F_{b-a}$ be general
rulings of $X(a,b)$. Together with $\Lambda$ they span a projective
space $\Phi$ of dimension $b-a+1$. Choose a general subspace $\Psi$
of dimension $b-a-1$ in $\Phi$ and project down $X(a,b)$ from $\Psi$
to $\p^ {a+2}$. The projection is a scroll surface $Y(a,b)\subset
\p^ {a+2}$ of degree $d=a+b\geq 2a+1$ which has a directrix $L$, the
image of $\Lambda$, of multiplicity $e=b$. However we have here
$\mu=a$ because, if $x\in L$ is the general point, only $a$ among
the $b$ lines of the ruling through $x$ vary, the other $b-a$ stay
fixed and coincide with $L$.

\begin{Theorem}\label{serie} For every $r\geq 3$ and for every $d\geq
2r-3$ the hypersurface $Y(r-2,d-r+2)^ *\subset {\p^r}^*$ of degree
$d$ is homaloidal.
\end{Theorem}

\begin{proof} We keep the above notation.
A repeated use of  \eqref  {eq:inclus2} gives
 $$Y(r-2,d-r+2)^ *= X(r-2,d-r+2)^*\cap \Psi^ \perp.$$ To
simplify the notation, set $X=X(r-2,d-r+2)$ and $Y=Y(r-2,d-r+2)$.

According to Remark \ref{dualrig}, $X^ *$ has vanishing Hessian and
the image of its polar map is a rational normal scroll
$Z=Z(r-2,d-r+2)$ of dimension  $r$ and degree $d-2r+5$. By Theorem
\ref {fibrespace} the general fibre of the polar map is a projective
subspace of dimension    $d-2r+4$.

Let us
repeat all pertinent dimensions translating from above $a,b$ to
present $d,r$:
$$\dim (X)=2,\quad \dim (X^*)=d-r+3\geq r \quad \mbox{\rm (from the
assumed inequality)} $$
$$\dim (Y)=2, \quad \dim (Y^*)=r-1, \quad \dim(\Phi)=d-2r+5,  \quad \dim (\Psi)=d-2r+3$$
\begin{eqnarray*}\dim (Z)&=&\dim (\Psi^\perp)=d-r+4-\dim (\Psi) -1\\
&=&d-r+4-(d-r+2-(r-2)-1)-1=r.
\end{eqnarray*}

By a repeated use of Lemma~\ref{projpoint} in a dual form, one has:
\begin{equation}\label{bello}\phi_{_{Y^*}}=\sigma_\Psi\circ
(\phi_{_{X^*}})_{|\Psi^\perp}.\end{equation}

We claim that the map $(\phi_{_{X^*}})_{|\Psi^\perp}:\Psi^ \perp\map
Z$ is birational. By part (ii) of Lemma~\ref{ratcurve} and duality,
if $z\in Z$ is a general point, then $\Phi^ \perp\cap T_{Z,z}^
\perp=\emptyset$. Let $\xi\in {\p^r}^*$ be an inverse image of $z$
by $\phi_{_{X^*}}$. Then
$<\xi,\Phi^ \perp>\cap<\xi,T_{Z,z}^\perp>=\{\xi\}$. Assuming, as we may,
that $\Psi^ \perp$ is a general subspace of dimension $r$ through
$<\xi,\Phi>$, then $\Psi^ \perp\cap <\xi,T_{Z,z}^\perp>=\{\xi\}$ and
moreover $\Psi^\perp$ intersects the fiber of $\phi_{X^*}$ over $z$
only at $\xi$ (see Propositions \ref{hessian} and Theorem \ref {fibrespace}).

By \eqref{bello} and Theorem \ref {fibrespace},
 the degree of the polar map $\phi_{_{Y^*}}$ is the
same as the degree of the restriction of the projection
$\sigma_\Psi$ to $Z$. To compute this latter degree, note that
$\Psi$ intersects $Z$ exactly in $d-2r+4$ distinct points, namely
the intersections of $\Psi$ with each one of the $d-2r+4$ planes
spanned by $\Lambda$ and by one of the $d-2r+4$ rulings spanning
$\Psi$ together with $\Lambda$. We claim that the intersection of
$\Psi$ with $Z$ at these points is transversal. Indeed, by
projecting from $\Lambda$ to $\p^ {d-r+2}$, we see that $X$ maps to
a rational normal curve $C$, the lines $F_1,...,F_{d-2r+4}$ map to
points $p_1,...,p_{d-2r+4}$ on $C$ and $\Psi$ maps to
$\Pi=<p_1,...,p_{d-2r+4}>$. By part (i) of Lemma \ref{ratcurve},
$\Pi$ intersects the projection of $Z$ transversally at
$p_1,...,,p_{d-2r+4}$. The claim follows.

Thus the restriction of the projection $\sigma_\Psi$ to $Z$
coincides with the projection of $Z$ from $d-2r+4$ independent
points on it. Since, as seen, $\deg(Z)=d-2r+5$, the restriction of
the projection $\sigma_\Psi$ to   $Z$ is a birational map of $Z$ to
$\p^ r$, thus completing proof.
\end{proof}

\begin{Remark}\label{duale}\rm As in
{ Section} \ref{multidir} and in the description
before Proposition \ref{multi_scrolls}, one can take the dual
viewpoint to describe the homaloidal hypersurfaces we constructed
above.

More precisely, to obtain the dual of $Y(a,b)$, one can proceeds as
follows. Consider the scroll $X_1=S(1^ a,b)\subset\p^{b+2a}$ of
degree $d=a+b$ and dimension $a+1$, with $S(1^ a,b)\simeq
\Proj(\O_{\Proj^ 1}(1)^ {\oplus a}\oplus  \O_{\Proj^ 1}(b))$
embedded via the $\O(1)$ bundle. Consider also the rational normal
scroll $X_2=S(1^ a)\subset\p^{2a-1}$ of degree $a$ and dimension
$a$. Clearly $X_2\subset X_1$. Take $b-a$ general rulings
$F_1,...,F_{b-a}$ of $X_1$. The span $\Sigma=<X_2,F_1,...,F_{b-a}>$
has dimension $2a-1+b-a=a+b-1$. Take a sufficiently general subspace
$\Sigma'$ of dimension $b+a-3$ intersecting $\Sigma$ in a general
subspace of dimension $b-2$, and project form $\Sigma'$ down to $\p^
{a+2}$. The image of $X_1$ is a hypersurface $V$ with a subspace
$\Pi$ of dimension $a$, the image of $\Sigma$, of multiplicity $b$,
since it is the image of $X_2$ and of $F_1,...,F_{b-a}$. The
hypersurface $V$ is the dual of $Y(a,b)$.

In this way we see that $Y(a,b)$ is a section
of $S(1^ a,b)^ *$ made with a suitable linear space of dimension $a+2$.
Since $S(1^ a,b)$ is a suitable linear section of  $S(1^ {a+b})={\rm
Seg}(1,a+b-1)$, which is self--dual, we see that $S(1^ a,b)^ *$ is a
suitable projection in $\p^{b+2a}$ of ${\rm Seg}(1,a+b-1)$.
\end{Remark}

\begin{Remark}\label{general} \rm In the construction of
$Y(a,b)$, it is not necessary that the rulings $F_1,...,F_{b-a}$ be
distinct. Indeed, one can consider an effective divisor
$D=m_1F_1+...+m_hF_h$ of degree $b-a$ formed by lines of the ruling
of $X(a,b)$. Then in the above construction one replaces the
subspace $\Phi$ with the span of $\Lambda$ and of the osculating
spaces of order $m_i$ at the points $p_i\in C$ projections of the lines
$F_i$, $i=1,...,h$. We will denote the resulting surface by
$Y(a,b;m_1,...,m_h)$.

The corresponding hypersurfaces in Theorem~\ref{serie} are still
homaloidal, since the proof of the theorem works even in this
special situation: indeed the intersection of $\Psi$ with
$Z(r-2,d-r+2)$ is no longer formed by $d-2r+4$ distinct points, but
by a $0$--dimensional scheme of length $d-2r+4$, formed by $h$
points $x_1,...,x_h$ , with length $m_1,...,m_h$ respectively, hence
the projection of $Z(r-2,d-r+2)$  to $\p^ r$ from $\Psi$ is still
birational.

As we will see in the next section however, this specialization influences
the degree of the inverse of the resulting polar map.\end{Remark}

\begin{Remark}\label{moduli}\rm The scroll $S(a,b)$, with $0<a<b$, has a
group of dimension $b-a+5$ of projective transformations which fixes
it and all scrolls $S(a,b)$ are projectively equivalent, i.e.
$S(a,b)$ has no projective moduli.

The scroll $X(a,b)$ has a group of dimension
$\max\{0,b-3a+7\}$ of projective transformations which fixes it,
and the scrolls $X(a,b)$ have no projective moduli if and only
if $b-3a+7\geq 0$.

Assume $b-a\leq 3$. Then
the subgroup fixing $F_1,...,F_{b-a}$ has
dimension $\max\{0,7-2a\}$. In conclusion $Y(a,b)$ has a group of
dimension $\max\{0,7-2b\}$ of projective transformations which
fixes it, and there are no projective moduli if and only if
$2b\leq 7$, i.e. $b\leq 3$.

If $b-a\geq 4$ the subgroup fixing $F_1,...,F_{b-a}$ has
dimension $\max\{0,b-3a+4\}$, and $Y(1,b)$ has no continuous
group of projective transformations because
$b-3a+4-2(b-a)=4-a-b<0$. In this case $Y(a,b)$
has projective moduli.

This implies that, except for $Y(1,2)^*$, the
homaloidal hypersurfaces we constructed here cannot be related to
pre--homogeneous vector spaces. The same holds for $Y(1,2)^*$,
as we will see later (see Theorem \ref {subhankel}, (iii) and (iv), and Example \ref {subhankel34}).
\end{Remark}

We finish this section by producing families of homaloidal
hypersurfaces in $\p^ r$, which are different from the above ones
as soon as $r\geq 4$.  They do not seem to be related in general
to hypersurfaces with vanishing Hessian. For $r=3$ instead one
essentially recovers the above examples.

Let $X\subset \p^ r$ be a non--degenerate scroll surface of degree
$d$ with a line directrix $L$ of multiplicity $e=r-2$, with $\mu=e=r-2$,
i.e. such that there are $r-2$ variable rulings $F_{x,1},\dots, F_{x,r-2}$ of
$X$ passing through the general point $x\in L$. According to Proposition
\ref {Dhyper}, if $\nu\leq r-2$, i.e. if $F_{x,1},\dots, F_{x,r-2}$ and $L$ do not span a
hyperplane, then $X^ *$ has vanishing Hessian. We will assume instead
that $\nu=r-1$ and that the hyperplane $<L,F_{x,1},\dots, F_{x,r-2}>$
homographically varies in a pencil when $x$ moves on $L$.

\begin{Example}\label{hypscroll}\rm ({\it Scrolls in $\p^r$ with
 line directrix, having $e=\mu=r-2$ and $\nu=r-1$}).
 Take a curve $C$ of degree $n\geq 2r-5$ in a $(r-2)$--dimensional subspace
 $\Pi$ of $\p^ r$, having a  $(n-r+2)$--secant $(r-4)$--dimensional subspace $\Pi'$.
 Assume also that the general hyperplane in $\Pi$ through $\Pi'$
 cuts $C$, off $\Pi$, in $r-2$ independent points.
 Curves of this sort are not difficult to construct. The first instance, is for
 $r=4$ in which case $C$ is a plane curve of degree $n\geq 3$ with a
point $O$ of multiplicity $n-2$. Note that, for $r\geq 4$ these curves
need not to be rational.

 Take a line $L$ in $\p^ r$ skew with $\Pi$  and
 set up an isomorphism between
$L$ and the pencil of hyperplanes through $\Pi'$ in $\Pi$. Fix a general point
$x\in L$, let $\Pi_x$ be the corresponding hyperplane in $\Pi$ through $\Pi'$
and let $x_1,\dots, x_{r-2}$ be the intersection points of $\Pi_x$ with $C$ off $\Pi$.
Then let $F_{x,i}$ be the line joining $x$ with $x_i$, $i=1,\dots,r-2$.
As $x$ varies on $L$, the lines $F_{x,1},\dots, F_{x,r-2}$ describe a scroll $X$
of the aforementioned type: the hyperplane
$<L,F_{x,1},\dots, F_{x,r-2}>=<L,\Pi_x>$ varies in the pencil of hyperplanes
through $<L,\Pi'>$.

The degree of such a scroll is $d=n+r-2$, as one sees by cutting it with a
general hyperplane through $L$.
\end{Example}
\medskip

\begin{Theorem}\label{Dhyperi}
Let $X\subset \p^ r$, be a non--degenerate scroll surface of degree
$d$ with a line directrix $L$ of multiplicity $e=r-2$ and with $\mu=e=r-2$.
Let $F_{x,1},\dots, F_{x,r-2}$ be the variables rulings of
$X$ passing through the general point $x\in L$. Suppose that  $<L,F_{x,1},\dots,F_{x,r-2}>$
is a hyperplane in $\p^r$ varying
homographically  in a pencil when $x$ moves on $L$.
Then $X^ *\subset {\p^r}^*$ is a
homaloidal hypersurface.
\end{Theorem}

\begin{proof} The space $\Sigma=L^ \perp$ has multiplicity $d-r+2$ for
$V=X^ *$ and the general hyperplane $\xi=x^ \perp, x\in L$, through
$\Sigma$, cuts out on $V$ a hypersurface $V_\xi$ formed by
$\Sigma$ with multiplicity $d-r+2$ and $r-2$ more $(r-2)$--dimensional subspaces
$\Sigma_i:=F_{x,i}^ \perp$, $i=1,\dots,r-2$, such that
$\Sigma\cap\Sigma_1\cap\dots\cap \Sigma_{r-2}=\{p\}$, where
$p=<L,F_{x,1},\dots,F_{x,r-2}>^ \perp$. Hence, as $\xi$ varies, $p$
homographically describes a line $\Lambda$ in $\Sigma$.

Consider the subspaces $T_{i,j}=\Sigma_i\cap
\Sigma_j, 1\leq i<j\leq r-2$, which have multiplicity $2$
for $V_\xi$, whereas $p$ has multiplicity $d$ for $V_\xi$. Note that
all $T_{i,j}$, with $1\leq i<j\leq r-2$, belong to the singular locus of $V$ since they are
intersections of two rulings of the scroll $V$. Furthermore, they all contain $p$.

Let now $z$ be a general point in $\xi$, hence a general point in ${\p^ r}^ *$.
The polar hyperplane $\pi_z$ of $z$ with respect to $X^ *$ contains $p$.
However it cannot contain the line $\Lambda$, otherwise all the polar
hyperplanes would contain $\Lambda$ and, by the reciprocity theorem, the
first polars of $X^ *$ with respect to the points of $\Lambda$ would vanish
identically, i.e. the points of $\Lambda$ would all have multiplicity $d$ for
$X^ *$, which would be a cone, a contradiction, because $X$ is
non--degenerate.

This proves that if $\pi_z=\pi_{z'}$ then $z'$ lies in $\xi=<z,\Sigma>$.
To finish our proof, we have to prove that all first polars through $z$
intersect $\xi$ only at $z$, off the singular locus of $V$.
To see this, first consider the polars with respect to points $y\in
\xi$, and containing $z$. By Remarks \ref {arrangm} and \ref {ext}, the closure
of their  intersection off the singular locus of $V$ is the line $\ell=\langle z,p\rangle$.
There is finally one more independent polar  through $z$ which we have to
take into account. It passes however with multiplicity $d-r+1$ through
$\Sigma$, hence it cuts $\xi$ along $\Sigma$ counted with
multiplicity $d-r+1$, plus another hypersurface $V'$ of degree $r-2$, which,
as we saw,
contains all the subspaces $T_{i,j}$, $1\leq i<j\leq r-2$.
By applying Lemma \ref  {lem:ess}, we see that $V'$ has multiplicity $r-3$ at
$p$. Hence $V'$ intersects $\ell$ only in $z$
and $p$, thus proving that $z$ is the only point having the polar
hyperplane $\pi_z$, i.e. the assertion.
\end{proof}

\begin{Remark}\label{reflection} Note that the homaloidal hypersurfaces
in Proposition \ref{Dhyperi} are reminiscent, in its structure, to
the homaloidal hypersurface $F_4=V(f^{(4)})$ in
Theorem~\ref{subhankel} below.
\end{Remark}

\subsection{Examples in $\p^ 3$ revisited}\label{sing}

In this section we want to revisit the examples of homaloidal
surfaces $Y(1,d-1)^ *$ in $\p^ 3$ of degree $d\geq 3$ constructed in
Theorem \ref {serie}. We want to  analyze the singularities of these
surfaces and understand the degree of the inverse of the polar map.

First of all,  the scrolls $Y(1,d-1)$ are self--dual, i.e.
$Y(1,d-1)^*$ is projectively equivalent to $Y(1,d-1)$ (see
Proposition \ref {autodual}). The surface $Y(1,d-1)$ has a line $L$
of multiplicity $d-1$ and no other singularity. One obtains the
desingularization $S(1,d-1)$ of $Y(1,d-1)$ by simply blowing--up
$L$. The pull--back of $L$ on $S(1,d-1)$ consists of the line
directrix $E$ plus $F_1,...,F_{d-2}$ rulings. This means that $L$ is
the intersection of $d-1$ distinct, generically smooth, branches,
$X, X_1,...,X_{d-2}$ respectively corresponding to $E,
F_1,...,F_{d-2}$. The branches $F_1,...,F_{d-2}$ intersect
transversally at a general point of $L$, whereas the branch $X$
glues with the branch $X_i$ at the point $O_i$, which is the image
of the intersection point $O\,'_i$ of $E$ with $F_i$, $i=1,...,d-2$.

We want to resolve the singularities of
the polar maps. We will see that, in order to do so, it is not sufficient
to blow--up $L$, but one has to perform further blow--ups.

In order to illustrate this, we analyze in detail the case $d=3$.
The other cases can be treated similarly, and we will briefly
discuss them later.

Consider the surface $F$ with equation:

$$f=x_2^ 3-2x_1x_2x_3+x_0x_3^ 2=0.$$
It will be shown later that $F=Y(1,2)$ (see Example \ref {subhankel34}). The partial derivatives of
$f$ are:

$$f_0=x_3^ 2, f_1=x_2x_3, f_2=3x_2^ 2-2x_1x_3, f_3=-2x_1x_2+2x_0x_3.$$

The double line $L$ has equation $x_2=x_3=0$. Now we pass to affine
coordinates $x=\frac {x_1}{x_0}, y=\frac {x_2}{x_0}, z=\frac
{x_3}{x_0}$, so that the equation of the surface becomes:
$$z^ 2+y^ 3-2xyz=0,$$
whereas the first polars $V(f_i)$, $i=0,\dots,3$, become:
\begin{equation}\label{firstpo} z^ 2=0,\;\;yz=0,\;\; 3y^ 2-2xz=0,\;\;z-xy=0
\end{equation}
and $L$ is the line $y=z=0$. Blow--up this line. To do this,
introduce coordinates $(x,y,\xi)$, the blow--up map being:
$$(x,y,\xi)\to (x,y,y\xi)$$
The exceptional divisor $M$ of the blow--up has equation $y=0$. The
strict transform $F\,'$ of the surface $F=Y(1,2)$ has equation:
$$\xi^ 2+y-2x\xi=0$$
which is smooth. A similar analysis at the infinity, shows that the
singularities of  $Y(1,2)$ can be resolved with one single blow--up
along the double line $L$. The proper transform of $L$ has now
equation:
$$y=0,\;\; \xi(\xi-2x)=0,$$
which is the union of two smooth rational curves on $M$, meeting at
the point $O\,'=(0,0,0)$ which maps to the origin $O$ in $\A^ 3$,
which is where the two branches of $Y(1,2)$ through $L$ glue.
Consider now the proper transform of the first polars
\begin{equation}\label{firstblow}
y\xi^ 2, y\xi, 3y-2x\xi, \xi-x.
\end{equation}
 We see that all these pass through  $O\,'$.
Thus, in order to resolve the singularities of the polar map, one
still has to blow--up $O\,'$ -- though, we emphasize, this is no
longer a singular point of $F'$.

This tells us that the scheme $S={\rm Sing}(Y(1,2))$ is not reduced:
it consists of the line $L$ with an embedded point at $O$. There is
no need to blow up in order to understand the structure of this
embedded point -- it suffices to analyze the affine equations
\eqref{firstpo} of the first polars. The scheme in question is a
subscheme of the surface of equation $z=xy$, whose coordinate ring
is $k[x,y,z]/(z-xy)\simeq k[x,y]$. Hence we interpret the scheme $S$
as the subscheme of $\A^ 2$ defined by the equations:
$$x^ 2y^ 2=0,\;\;xy^ 2=0,\;\;3y^ 2-2x^ 2y=0.$$
The line $L$, which has now equation $y=0$ splits off, leaving a
zero--dimensional scheme $S\,'$ supported at the origin $O$, which
is responsible for the embedded point of $S$. The equations  of
$S\,'$ are:
$$x^ 2y=0,\;\;xy=0,\;\;3y-2x^ 2=0.$$
This is now a subscheme of the smooth curve $C$ of equation $3y-2x^
2=0$, which is simply tangent to $L$ at $O$. The coordinate ring of
$C$ is $k[x,y]/(3y-2x^ 2)\simeq k[x]$ and the scheme $S'$ has now
the equations $x^ 3=0$. Summing up, the embedded point at the
origin is due to the fact that all the polars have multiplicity of
intersection $4$ with the curve $C$ at $O$. We thus see that we will
have to blow--up along $L$ and then three more times at subsequent
infinitely near points to resolve the singularities of the polar
map.

\begin{Remark}\label{effbeh} Note that, after the first blow--up, the
polar system is given by the system  \eqref{firstblow}. The base
point scheme is now zero dimensional supported at $O\,'$. Indeed it
is a curvilinear scheme $T$ of length $3$ along the proper transform
$C\,'$ of the curve $C$, defined by the equations $x=\xi, 3y=2x^ 2$.
Note however that $F'$ has only intersection multiplicity $2$ with
$C\,'$ at $O\,'$. This means that $F'$ does not contain $T$. In
other words, the fourth (and last) point, infinitely near to $O$, to
be blown--up in order to resolve the singularities of the polar map,
does not even belong to the original surface  $F$.
\end{Remark}

The above analysis gives another reason why $Y(1,2)$ is a homaloidal
surface. Indeed the polar system is formed by quadrics through $L$.
The general such quadric is smooth, as we see from the equations
of the polars or from \eqref {firstpo}. The residual
intersections of two general polars off $L$ are rational normal
cubics, i.e.,  curves of type $(1,2)$ on the general such quadric.
The self--intersection of these curves is therefore $4$. However,
the curves in question have to contain the $0$--dimensional scheme
of length $3$ supported at $O$, which is responsible for the
embedded point of ${\rm Sing}(Y(1,2))$ on $L$. This drops the
self--intersection of the system of cubic curves to $1$ and explains
why the polar map is birational.

We emphasize that the polar map is a quadratic transformation of
$\p^ 3$ which is a degenerate case of the well known quadratic
transformation defined by all quadrics passing through a given line
$L$ and three distinct general points $p_1,p_2,p_3$ (see
\cite{Conf}, \cite{Pan}).
\medskip

\begin{Remark} \label{genscroll} It is worth comparing the behavior of
the polar map of $Y(1,2)$ with the one of the general projection $Y$
of $S(1,2)$ to $\p^ 3$. We may think of $Y$ as the surface defined
by the equation $x_1(x_2^ 2+x_3^ 2)-2x_0x_2x_3=0$, whose double line
$L$ has the equations $x_2=x_3=0$. The resolution of the
singularities of $Y$ is obtained by blowing up along $L$. In this
way one recovers $S(1,2)$, and the proper transform of $L$ is a
conic $C$, which projects $2:1$ to $L$, with two branch points,
located at the points $O_1, O_2$ with affine coordinates $(0,0,0)$
and $(1,0,0)$. The scheme ${\rm Sing}(Y)$ consists of $L$ with two
embedded points of length $2$ at $O_1$ and $O_2$. This yields degree
$2$ for the polar map. The surface $Y(1,2)$ can be thought of as
obtained from $Y$ when $O_1$ and $O_2$ collapse together. Indeed the
conic $C$ then splits as the union of the line directrix $E$ of
$S(1,2)$ and a ruling. This also clarifies why $Y(1,2)$ coincides
with the surface $F$, which is a member of a series of homaloidal
hypersurfaces under the general name of sub--Hankel hypersurfaces,
to be dealt with in the next section -- in the notation of that
subsection, one has $F=V(f^{(3)})$ (see Example~\ref{subhankel34}).
\end{Remark}

The analysis of the general case $Y(1,d-1)$ is similar. The general
point $p$ of $L$ has multiplicity $d-1$ and it is the intersection
of $d-1$ smooth branches of $Y(1,d-1)$ containing $L$ and pairwise
intersecting transversally along $L$ around $p$. There are however
$d-2$ points $O_1,\ldots, O_{d-2}$ on $L$ around which $Y(1,d-1)$
looks like the union of $d-3$ branches which intersect transversally
along $L$ around $O_i$, plus another branch which is analytically
equivalent to $Y(1,2)$ at $O$ and which is generically located with
respect to the previous $d-3$ branches. The singularities of
$Y(1,d-1)$ can be resolved by blowing-up along $L$: in this way one
obtains $S(1,d-1)$ and the blowing-up map is nothing but the
projection $S(1,d-1)\to Y(1,d-1)$.

The general polar has a point of multiplicity $d-2$ at a general
point of $L$. It is again resolved when we blow--up $L$. However,
for the same reason as in the case of $Y(1,2)$, after blowing up,
there are $d-2$ curvilinear schemes of length $3$ supported at each
of the points $O'_1,...,O'_{d-2}$, which belong to the base locus of
the proper transform of the polar system. Another way of saying this
is that there are $d-2$ embedded points $O_1,...,O_{d-2}$ in the
scheme structure of ${\rm Sing}(Y(1,d-1))$ along $L$ supported at
$O_1,...,O_{d-2}$.

Again this explains the reason why the polar map is birational. Let
$\Phi$ be the proper transform of the general first polar after
having blown-up $L$ . This is a rational scroll. Let us denote by
$R$ the general ruling and by $D$ the proper transform of $L$, which
is a section. If $H$ is the pull--back of a general plane section,
we have $H\equiv D+R$. Since $H^ 2=d-1$, we find $D^ 2=d-3$. If
$\Gamma$ is the trace on $\Phi$ of the proper transform $\L$ of the
polar system, we have $\Gamma\equiv (d-1)H-(d-2)D\equiv (d-1)R+D$,
thus $\Gamma^ 2=3d-5$. Notice however that the trace of $\L$ on
$\Phi$ has $d-2$ base point schemes each of length $3$. After having
further blown-up these base point schemes, this reduces the
self--intersection of $\L$ to $1$.

The analysis is more complicated for the surfaces
$Y(a,b\,;\,m_1,...,m_h)$, $m_1+...+m_h=d-2$,  described in
Remark~\ref{general}. We will merely outline the results that can be
checked by a careful treatment. The singularity  can still
be resolved with a simple blow--up along $L$ thus getting
$S(1,d-1)$. The proper transform of $L$ is now
$E+m_1F_1+...+m_hF_h$. This means that $L$ is the intersection of
$h+1$ branches, a smooth one $X$, corresponding to the line
directrix $E$ of $S(1,d-1)$, the other branches $X_1,...,X_h$ are
instead cuspidal of orders $m_1,...,m_h$ corresponding to the
rulings $F_1,...,F_h$ respectively.

As for the degree of the inverse map, one has the following. First,
the degree of the inverse map of the polar map $\phi$ of $Y(1,d-1)$
coincides with the degree of the image of a general plane $\pi$ via
$\phi$. The linear system cut out on $\pi$ by the system of the
first polars, is a $3$--dimensional linear system of curves of
degree $d-1$ with only one ordinary base point $x$ of multiplicity
$d-2$, i.e. $x$ is the intersection of $L$ with $\pi$. Thus the
image of $\pi$ has degree $(d-1)^ 2-(d-2)^ 2=2d-3$. Note that, for
$d=3$, one retrieves the expected degree of the inverse to the polar
map of the specialized Hankel determinant (see
Remark~\ref{jonquieres}, (c))

Consider now the surface $Y(a,b\,;\,m_1,...,m_h)$ with
$m_1+...+m_h=d-2$. The general first polar has again a point of
multiplicity $d-2$ at a general point $x\in L$. Moreover, a local
computation shows that it has tangency of order $m_i-1$ along the
plane $\pi_i$ tangent to the branch $X_i$, $i=1,...,h$. Arguing as
above, we see that this decreases the degree of the inverse of
$\phi$ by $\sum_{i=1}^ h (m_i-1)=d-2-h$. In particular, if $h=1,
m_1=d-2$, then we have the maximal drop of the degree of the
inverse, namely $d-3$, i.e. the degree of the inverse of $\phi$ is
$d$.

It would be interesting to have a similar analysis in $\p^ r$, for
$r>3$.

\section{Some determinantal homaloidal polynomials}\label{algex}

In this section we bring up a series of
examples of homaloidal polynomials which can be treated in an
algebraic fashion. Some of the proofs, though elementary in spirit,
are nevertheless quite involved.

\subsection{Degenerations of Hankel matrices}\label{hankelmatrices}

 First we need a few algebraic concepts (see \cite{cremona}
for more contextual details).

\begin{Definition}\rm  Let $R$ be a Noetherian ring and  let $I\subset R$ be an
ideal.

(1) Let ${\cal S}_R(I)\surjects {\cal R}_R(I)$ denote the structural
graded $R$-algebra homomorphism from the symmetric algebra of $I$ to
its \emph {Rees algebra}, i.e. the $R$-algebra that defines the
blowup along the subscheme corresponding to the
 ideal $I$ (see \cite[{Section} 5.2]
{Eisenbook}). We say that $I$ is of {\it linear type\/} if this map
is injective{\rm ;}

(2) If $R$ is a Noetherian local ring (or a standard graded ring
over a field) the ideal $I$ is said to be {\it perfect\/} if it has
finite homological (i.e., projective) dimension over $R$ and this
attains its minimal possible value, namely, the codimension of $I$
(see \cite[p. 485]{Eisenbook}). It is known that if $R$ is moreover
a Cohen--Macaulay ring (e.g., regular) then an ideal $I$ is perfect
if and only if $R/I$ is Cohen--Macaulay{\rm ;}

(3) An ideal $I\subset R$ of linear type satisfies  the {\it
Artin--Nagata condition $G_{\infty}$\/}  (see \cite{AN}) which states
that the minimal number of generators of $I$ locally at any prime
$p\in \spec R$ is at most the codimension of $p$. This condition is
equivalent to a condition in terms of a free presentation
$$R^m\stackrel{\varphi}{\lar} R^n \lar I \lar 0 $$
of $I$, namely:
\begin{equation}\label{f1}
{\rm cod} (I_t(\varphi))\geq {\rm rank}(\varphi)-t+2,\quad \mbox{\rm
for}\quad 1\leq t\leq {\rm rank}(\varphi),
\end{equation}
where $I_t(\varphi)$ denotes the determinantal ideal of the $t\times
t$ minors of a representative matrix of $\varphi$
 (see, e.g., \cite[{
Section} 1.3]{wolmer}){\rm ;}

(4) Suppose that $R$ is standard graded over a field $k$ and $I$ is
 generated by forms of a given degree $s$. In this
case, $I$ is more precisely given by means of a free graded
presentation
$$R(-(s+1))^{\ell}\oplus\sum_{j\geq 2} R(-(s+j)) \stackrel{\varphi}{\lar} R(-s)^n\rar
I\rar 0$$
 for suitable $\ell$.
We call the image of $R(-(s+1))^{\ell}$ by $\varphi$ the {\it linear
part\/} of $\varphi$ and say that the corresponding submatrix
$\varphi_1$  has {\it maximal rank\/} if its rank is $n-1$ ($={\rm
rank}(\varphi)$). Clearly, the latter condition is trivially
satisfied if $\varphi_1=\varphi$, in which case $I$ is said to have
{\it linear presentation\/} (or is {\it linearly presented\/}).

We remark that such an ideal, if it is of linear type, then it is
generated by algebraically independent elements over $k$. In
particular, if $R=k[\xx]=k[x_0,\ldots,x_r]$ and $I$ happens to be of
linear type and generated by $r+1$ forms of the same degree then
these forms define a dominant rational map $\p^r\dasharrow \p^r$.
\end{Definition}
\subsubsection{Arithmetic of sub-Hankel matrices}

So much for generalities. We now introduce the main object of this
part, which is a degeneration of a generic Hankel matrix over a
polynomial ring by specializing convenient entries to zero (see
\cite{cremona2} for further classes of specializations of square
generic matrices whose determinants are often homaloidal, treated
within the general framework of the theory of ideals).

Let $y_1,\ldots,y_{r+1}$ be variables over a field $k$ and
set

$$
M^{(r)}=M^{(r)}(y_1,\ldots,y_{r+1})= \left(
\begin{matrix}
y_1&y_2&y_3&...&y_{r-1}&y_{r}\\
y_2&y_3&y_4&...&y_{r}&y_{r+1}\\
y_3&y_4&y_5&...&y_{r+1}&0\\
.&.&.&...&.&.\\
.&.&.&...&.&.\\
.&.&.&...&.&.\\
y_{r-1}&y_{r}&y_{r+1}&...&0&0\\
y_{r}&y_{r+1}&0&...&0&0\\
\end{matrix}
\right )
$$
Note that the matrix has two tags: the upper index $(r)$ indicates
the size of the matrix, while the variables enclosed in parentheses
are the total set of variables used in the matrix. We call attention
to the notation as several of these matrices will be considered with
variable tags throughout, though we will often omit the list of
variables if they are sufficiently clear from the context.

This matrix will be called a {\it generic sub--Hankel matrix\/};
more precisely, $M^{(r)}$ is the {\em generic sub--Hankel matrix of
order $r$ on the variables} $y_1,\ldots,y_{r+1}$. Its determinant, a
form of degree $r$, will be the central object of this section.
Throughout we fix a polynomial ring $k[x_0,\ldots,x_r]$ which will
be the source of all lists of variables appearing in the various
such matrices considered heretofore. We will denote by
$f^{(r)}(x_0,\ldots,x_r)$ the determinant of
$M^{(r)}(x_0,\ldots,x_r)$ for any $r\geq 1$, and we set $f^{(0)}=1$.
We also set
$\phi^{(j)}=\phi^{(j)}(x_{r-j},\ldots,x_r):=f^{(j)}(x_{r-j},\ldots,x_r)$.

We now head on to the main result concerning generic sub-Hankel
matrices. First we need the following algebraic structural lemmas
about the partial derivatives of $f^{(r)}$.

\begin{Lemma}\label{basic_lemma} Let $r\geq 2$. Then:

\begin{itemize}
\item[{\rm (i)}] One has

\begin{equation}\label{phased_out_r}
\frac{\partial f^{(r)}}{\partial x_i}=(-1)^r\, x_r\, \frac{\partial
\phi^{(r-1)}}{\partial x_{i+1}}, \quad 0\leq i\leq r-2{\rm ;}
\end{equation}

\item[{\rm (ii)}] For $0\leq i\leq r-1$, one has

\begin{equation}\label{g.c.d.s}
\frac {\partial f^{(r)} } {\partial x_0}, \ldots, \frac{\partial
f^{(r)} } {\partial x_i}  \in k\left[x_{r-i},\ldots,x_{r}\right]
\end{equation}
and the g.c.d. of these partial derivatives is $x_r^{r-i-1}${\rm ;}
\item[{\rm (iii)}] For any $i\,$ in the range $1\leq i\leq r-1$, the following holds:
\begin{equation}\label{basic_linear_relation}
x_r\,\frac{\partial f^{(r)} }{\partial
x_i}=-\sum_{k=0}^{i-1}\frac{2i-k}{i}\,\, x_{r-i+k}\, \frac{\partial
f^{(r)} }{\partial x_k}.
\end{equation}
Moreover,
\begin{equation}\label{perfect_linear_relation}
x_r\, \frac{\partial f^{(r)} }{\partial x_r}=(r-1)x_0\,\frac{\partial
f^{(r)} }{\partial x_0}+ (r-2)x_1\, \frac{\partial f^{(r)}}{\partial
x_1}+\cdots + x_{r-2}\,\frac{\partial f^{(r)}} {\partial x_{r-2}}
\end{equation}
\end{itemize}
\end{Lemma}

\demo (i) We induct on $r$. For $r=2$, the relation is readily seen
to hold. To proceed, introduce the following sign function on
integers: $\xi(r)=1$ if $r\equiv 1,2 \pmod{4}$ and $\xi(r)=-1$ if
$r\equiv 0,3 \pmod{4}$. The following identity is easily
established:
\begin{equation}\label{identity1}
\xi(j)\,\xi(j-1)=(-1)^j.
\end{equation}
Equivalently one has
\begin{equation}\label{identity2}
(-1)^{j+1}\xi(j)=-\xi(j-1).
\end{equation}

Assume that $r\geq 3$. Expanding $f^{(r)}$ by Laplace along the
first row one finds
\begin{equation}\label{Laplace_first_row}
f^{(r)}=-\xi(r)\,\sum_{j=0}^{r-1}\xi({j}) x_j
x_r^{r-j-1}\phi^{(j)}.
\end{equation}
By the same token, expanding $\phi^{(r-1)}$ by Laplace along the first
row one finds
\begin{equation}\label{Laplace_first_row_1}
\phi^{(r-1)}=-\xi(r-1)\,\sum_{j=1}^{r-1}\xi(j-1) x_j
x_r^{r-j-1}\phi^{(j-1)}.
\end{equation}

Suppose now $0\leq i\leq r-1$. Taking $x_i$-derivatives of both
sides of (\ref{Laplace_first_row}), for $i$ in this range, yields
\begin{equation}\label{derivatives_of_fr}
\frac{\partial f^{(r)}}{\partial
x_i}=-\xi(r)\left(\xi(i)x_r^{r-i-1}\phi^{(i)}+\sum_{j=1}^{r-1}\xi(j)
x_j x_r^{r-j-1}\frac{\partial \phi^{(j)}}{\partial x_i}\right).
\end{equation}

Similarly, taking $x_{i+1}$-derivatives of both sides of
(\ref{Laplace_first_row_1}) in the range $0\leq i\leq r-2$, yields
\begin{equation}\label{derivatives_of_fr_1}
\frac{\partial \phi^{(r-1)}}{\partial
x_{i+1}}=-\xi(r-1)\left(\xi(i)x_r^{r-i-2}\phi^{(i)}+
\sum_{j=1}^{r-1}\xi(j-1) x_j x_r^{r-j-1}\frac{\partial
\phi^{(j-1)}}{\partial x_{i+1}}\right).
\end{equation}

Thus, by the inductive hypothesis applied to $f^ {(i)}$, with
$i<r$, hence to $\phi^ {(i)}$, with $i<r$, and by the identity
\eqref{identity1}, we find

\begin{align*}
x_r\,\frac{\partial \phi ^{(r-1)}}{\partial
x_{i+1}}&=-\xi(r-1)\left(\xi(i)x_r^{r-i-1}\phi^{(i)}+
\sum_{j=1}^{r-2}\xi(j) x_j
 x_r^{r-j-1}\left[\xi(j)\xi(j-1)x_r\frac{\partial \phi^{(j-1)}}{\partial x_{i+1}}
\right]\right)\\
&=-\xi(r-1)\left(\xi(i)x_r^{r-i-1}\phi^{(i)}+\sum_{j=1}^{r-1}\xi(j) x_j
x_r^{r-j-1} \frac{\partial \phi^{(j)}}{\partial x_{i}}\right)
\end{align*}

Multiplying the last line above by $\xi(r-1)$, using
\eqref{identity2}, and drawing upon \eqref{derivatives_of_fr} as
multiplied by $\xi(r)$, one obtains

\begin{equation}\label{finalmente}
\xi(r)\frac{\partial f^{(r)}}{\partial x_i}=\xi(r-1)\, x_r\, \frac{\partial
\phi^{(r-1)}}{\partial x_{i+1}}.
\end{equation}
 Since $\xi(r)\,\xi(r-1)=(-1)^r$ for every $r\geq 1$ by
\eqref{identity1}, equation \eqref{finalmente}  yields
\eqref{phased_out_r}.
\medskip

(ii) We induct on $r$. Both assertions are readily verified for
$r=2$ since $f^{(2)}=x_0x_2-x_1^2$. Thus, assume that $r\geq 3$.
Note that if $i<r-j$, the form $\phi^{(j)}$ does not involve the
variable $x_i$, hence all its derivatives with respect to $x_i$
vanish, for $0\leq i\leq r-2$ and $j<r-i$. Thus, using
(\ref{derivatives_of_fr}) we immediately see that \eqref {g.c.d.s}
holds. As for the assertion about the g.c.d., this is easy in the
range $0\leq j\leq r-2$, since it follows form the expressions
(\ref{derivatives_of_fr}) and the inductive hypothesis applied to
$f^ {(i)}$, hence to $\phi^ {(i)}$, with $i<r$.

As for $i=r-1$, we still have the
expression

\begin{equation}\label{derivative_of_order_before_the_last}
\frac{\partial f^{(r)}}{\partial
x_{r-1}}=-\xi(r)\left(\xi(r-1)\phi^{(r-1)}+ \sum_{j=1}^{r-1}\xi(j) x_j
x_r^{r-j-1}\frac{\partial \phi^{(j)}}{\partial x_{r-1}}\right),
\end{equation}
coming from  (\ref{derivatives_of_fr}). To prove the assertion
about the $g.c.d.$, replace $\phi^{(r-1)}$ by its Euler expansion in
\eqref{derivative_of_order_before_the_last} and collect the two
terms in $x_{r-1}\,(\partial \phi^{(r-1)}/\partial x_{r-1})$. We get

\begin{eqnarray}\label{derivative_of_order_before_the_last_after_euler}
\frac{\partial f^{(r)}}{\partial
x_{r-1}}&=&\xi(r)\left(\sum_{j=1}^{r-2}\xi(j) x_j
x_r^{r-j-1}\frac{\partial \phi^{(j)}}{\partial x_{r-1}}\right)\\
\nonumber
 &+&(-1)^{r+1}\frac{1}{r-1}\left(\sum_{j=1}^{r-2} x_j\,\frac{\partial
\phi^{(r-1)}}{\partial x_j}\, + \, (r-1)\,x_{r-1}\,\frac{\partial
\phi^{(r-1)}}{\partial x_{r-1}}\, +\, x_r\,\frac{\partial
\phi^{(r-1)}}{\partial x_r}\right)
\end{eqnarray}

 Now the g.c.d. of the derivatives up
to order $i=r-2$ was found to be $x_r^{r-i-1}=x_r$. If the
derivatives up to order $i=r-1$ would have a nonunit g.c.d. then it had
to be $x_r$. Thus, assume as if it were that $x_r$ divides the left
hand side in
(\ref{derivative_of_order_before_the_last_after_euler}). Since $x_r$
divides the first summand in the right hand side of
(\ref{derivative_of_order_before_the_last_after_euler}) and, by the
inductive hypothesis applied to $f^ {(r-1)}$, $x_r$ divides the summands $x_j\,(\partial
\phi^{(r-1)}/\partial x_j)$, for $1\leq j\leq r-2$, then $x_r$ would
divide the derivative $\partial \phi^{(r-1)}/\partial x_{r-1}$, which
would contradict the inductive hypothesis as applied to $f^{(r-1)}$.

(iii) We begin with  \eqref{perfect_linear_relation}.
 The formula is readily verified for $r=2$ so we induct on $r\geq 3$.
Taking $x_r$-derivatives in (\ref{Laplace_first_row}),
multiplying by $x_r$ we get

\begin{eqnarray*}\label{derivative_of_order_last}
\lefteqn{ x_r\,\frac{\partial f^{(r)}}{\partial x_r} -\xi(r)\left(\sum_{j=0}^{r-2}(r-j-1)\xi({j}) x_j
x_r^{r-j-1}\phi^{(j)} +\sum_{j=0}^{r-1}\xi({j}) x_j x_r^{r-j-1}\,
\biggl(x_r\,\frac{\partial \phi^{(j)}}{\partial x_r}\biggr)\right)}\\
&&=-\xi(r)\left(\sum_{i=0}^{r-2}(r-i-1)\xi(i)x_ix_r^{r-i-1}\phi^{(i)}+
    \sum_{j=0}^{r-1}\xi({j}) x_j x_r^{r-j-1}\, \biggl(\sum_{i=0}^{r-2} (r-i-1)x_i\,
    \frac{\partial \phi^{(j)}}{\partial x_i}\biggr) \right)\\
&&=\,-\xi(r)\sum_{i=0}^{r-2}(r-i-1)x_i\left(\xi(i)x_r^{r-i-1}\phi^{(i)}+
    \sum_{j=0}^{r-1}\xi({j}) x_j x_r^{r-j-1}\, \,
    \frac{\partial \phi^{(j)}}{\partial x_i} \right)\\
    && =\,\,\sum_{i=0}^{r-2}(r-i-1)x_i\, \frac{\partial
f^{(r)}}{\partial x_i}
\end{eqnarray*}
where in the second line we applied the inductive hypothesis to $f^ {(\ell)}$,
for every $l=1,\ldots, r-1$, to wit
\begin{equation}\label{exact_linear_relation_for_induction}
x_r\, \frac{\partial \phi^{(l)}}{\partial
x_r}=\sum_{j=r-l}^{r-2}(r-j-1)\,x_j\,\frac{\partial
\phi^{(l)}}{\partial x_j}=\sum_{j=0}^{r-2}(r-j-1)\,x_j\,\frac{\partial
\phi^{(l)}}{\partial x_j}
\end{equation}
and in the fourth line we used the expression
obtained from multiplying \eqref{derivatives_of_fr} both sides by
$x_i$, for $i=0,\ldots, r-1$.

We now prove formula \eqref{basic_linear_relation}. In the range
$0\leq i\leq r-2$ the formula follows from \eqref{phased_out_r}.
Indeed, the formula is easily obtained for $r=2$. Inducting on $r$
in this range, we assume that

\begin{equation}\label{basic_linear_relation_shifted}
x_r\,\frac{\partial \phi^{(r-1)}}{\partial
x_{i+1}}=-\sum_{j=1}^{i}\frac{2i+1-j}{i}\, x_{r-i-1+j}\,
\frac{\partial \phi^{(r-1)}}{\partial x_j}
\end{equation}
holds. Therefore

\begin{align*}
x_r\,\frac{\partial f^{(r)}}{\partial
x_i}&=(-1)^rx_r\,\biggl(x_r\,\frac{\partial \phi^{(r-1)}}{\partial
x_{i+1}}\biggr)\\
&=-\biggl(\sum_{j=1}^{i}\frac{2i+1-j}{i}\,\, x_{r-i-1+j}\,
\biggl[(-1)^rx_r\frac{\partial \phi^{(r-1)}}{\partial x_j}\biggr]\biggr)\\
&=-\biggl(\sum_{j=1}^{i}\frac{2i+1-j}{i}\,\, x_{r-i-1+j}\,\,\frac{\partial f^{(r)}}{\partial x_{j-1}}\biggr)\\
&=-\biggl(\sum_{j=0}^{i-1}\frac{2i-j}{i} \,\,
x_{r-i+j}\,\,\frac{\partial f^{(r)}}{\partial x_{j}}\biggr),
\end{align*}
as was to be shown.

It remains to get the case where $i=r-1$. For this first note that a
repeated use of \eqref{phased_out_r} yields, for every $j=1,\ldots,
r-1,$ the relation

\begin{equation}\label{iterated_phased_out}
\xi(r)\, \frac{\partial f^{(r)}}{\partial x_{j-1}}= \xi(j)\,
x_r^{r-j}\,\frac{\partial \phi^{(j)}}{\partial x_{r-1}};
\end{equation}
Applying \eqref{perfect_linear_relation} to $\phi^{(r-1)}$
and using Euler's formula yields

\begin{align}\label{Euler}\nonumber
(-1)^rx_r\, [(r-1)\phi^{(r-1)}]=&\sum_{k=1}^{r-1}x_k\,
\biggl[(-1)^rx_r\frac{\partial \phi^{(r-1)}}{\partial
x_k}\biggr]+((-1)^rx_r) \biggl[x_r\frac{\partial \phi^{(r-1)}}{\partial
x_r}\biggr]\\ \nonumber =&\sum_{k=1}^{r-1}\,x_k\,\frac{\partial
f^{(r)}}{\partial
x_{k-1}}+\sum_{k=1}^{r-1}(r-1-k)\,x_k\,\biggl[(-1)^rx_r
\frac{\partial \phi^{(r-1)}}{\partial x_k}\biggr]\\
=&\sum_{k=1}^{r-1}(r-k)\,x_k\,\frac{\partial f^{(r)}}{\partial
x_{k-1}}.
\end{align}.

Combining \eqref {derivative_of_order_before_the_last} with \eqref{iterated_phased_out} and \eqref{Euler},  we get

\begin{align*}\label{relation_of_order_last} (r-1)\,x_r\,
\frac{\partial f^{(r)}}{\partial x_{r-1}}&
=-\,\xi(r)\left(\xi(r-1)\,x_r\,[(r-1)\phi^{(r-1)}]+(r-1)\sum_{j=1}^{r-1}
\,x_j\, \biggl[(\xi(j)x_r^{r-j}\,\frac{\partial \phi^{(j)}}{\partial
x_{r-1}}\biggr]\right)\\ \nonumber
&=-\left((-1)^r\,x_r\,[(r-1)\phi^{(r-1)}]+(r-1)\sum_{k=1}^{r-1}\,x_k\,\frac{\partial
f^{(r)}}{\partial x_{k-1}}\right)\\ \nonumber
&=-\,\sum_{k=1}^{r-1}(2r-1-k)\,x_k\,\frac{\partial f^{(r)}}{\partial
x_{k-1}},
\end{align*}
proving  \eqref{basic_linear_relation} also in this case.

This completes the proof of the lemma. \qed

\begin{Proposition}\label{perfection}
Let  $r\geq 2$.  Set $f=f^ {(r)}$. Then upon factoring out the
g.c.d. of $\partial f/\partial x_0,\ldots,\partial f/\partial x_{i}$,
the resulting polynomials generate a codimension two perfect ideal
$J_i\subset k[x_{r-i},\ldots,x_{r}]$ of linear type with linear
presentation.
\end{Proposition}
\demo  The assertion is readily checked for $r=2$ since
$f^{(2)}=x_0x_2-x_1^2$. Thus, we assume henceforth that $r\geq 3$.

Fix an $i\,$ in the range $1\leq i\leq r-1$. Let $J_i$ denote the
ideal of the ring $R^{[i]}=k[x_{r-i},\ldots,x_{r}]$ generated by the
partial derivatives of $f=f^{(r)}$ with respect to $x_0,\ldots,x_i$
divided by $x_r^{r-i-1}$.

We claim that the presentation matrix of $J_i$ is an $(i+1)\times i$
recurrent matrix having the form
\begin{equation}\label{phi_i}
\Phi^{[i]}=\left(
\begin{array}{c|c}
\hbox{2}\, x_{r-i}&\\
\hbox{${\frac{2i-1}{i}}\, x_{r-i+1}$}&\\
\vdots&\Phi^{[i-1]}\\\
\hbox{$\frac{i+1}{i}\, x_{r-1}$}&\\
x_r&\mbox{{\boldmath $0$}}\end{array} \right), \quad \mbox{{\boldmath
$0$}}=\underbrace{(0,\ldots,0)}_{i-1},
\end{equation}
where the first column comes from (\ref{basic_linear_relation}).
 By induction on $i$, one has that the last
$i-1$ columns of $\Phi^{[i]}$ are relations of $J_i$, hence the full
matrix $\Phi^{[i]}$ is a matrix of relations of $J_i$ and, moreover,
its linear part has maximal rank ($=i$).

On the other hand, by a well-known acyclicity criterion (see, e.g.,
\cite{BE}), in this case it suffices to check that the columns of
$\Phi^{[i]}$ are relations of the generators of $J_i$ and the
determinantal ideal $I_i(\Phi^{[i]})$ has codimension at least $2$.
Thus, we are left with finding two $i$-minors of $\Phi^{[i]}$
without nontrivial common factor. Let $\delta_1$ (respectively,
$\delta_2$) denote the minor obtained by deleting the first
(respectively, the last) row of $\Phi^{[i]}$. By induction on $i$,
$\delta_1=\pm x_r^i$ and $\delta_2$ admits a summand of the form
$\pm (i+1) x_{r-1}^i$ that results from multiplying the entries
along the anti-diagonal of the first $i$ rows -- indeed, by
\eqref{basic_linear_relation}, taking $k=i-1$, the coefficient of
the $(i,i-1)$ entry is $(i+1)/i$, hence the product is
$(i+1/i)(i/i-1)\cdots (3/2)(2/1)=i+1$. It follows that $\delta_1$ is
not divisible by $x_r$. Therefore, $I_i(\Phi^{[i]})$ has codimension
at least $2$, as required.

What we have proved so far is that $J_i$ has a Hilbert--Burch
resolution, and since it has codimension at least $2$ then it is a
codimension two perfect ideal. So, it remains to prove the last
contention of this item, namely, that $J_i$ is an ideal of linear
type. By \cite[Corollary 1.4.2 and Theorem 3.1.6]{wolmer} this will
be the case if the inequalities in (\ref{f1}) are fulfilled.

Since cod$(I_i(\Phi^{[i]}))\geq 2=i-i+2= {\rm rank}
(\Phi^{[i]})-i+2$, we only have to check that

\begin{equation}\label{f2}
{\rm cod} (I_t(\Phi^{[i]}))\geq i-t+2,\quad \mbox{\rm for}\quad
1\leq t\leq i-1,
\end{equation}
We proceed by induction on $i$, so ${\rm cod}
(I_t(\Phi^{[i-1]}))\geq i-1-t+2=i-t+1$ holds true in the range
$1\leq t\leq i-2$. Therefore, one needs, for every $t$ in the range
$1\leq t\leq i-1$, an additional $t$-minor of $\Phi^{[i]}$ which is
a nonzero-divisor on the ideal $I_t(\Phi^{[i-1]})$. Since
$\Phi^{[i-1]}$ has entries in the polynomial ring
$R^{[i-1]}=k[x_{r-i+1},\ldots, x_r]$, it suffices to show that there
exists such a minor effectively involving the extra variable
$x_{r-i}$. Supposing this were not the case, the full matrix of
relations $\Phi^{[i]}$ would have entries entirely contained in the
ring $R^{[i-1]}$ and, since the generators of $J_i$ are the maximal
minors of $\Phi^{[i]}$, they would all belong to $R^{[i-1]}$, which
is not the case.

This finishes the proof of the last statement.
\qed

\subsubsection{Sub-Hankel polynomials are homaloidal}\label {shp}

From the previous lemma ensues the following geometric result.

\begin{Theorem}\label{subhankel} Let $r\geq 2$.
Set $f=f^{(r)}$ and $J=(\partial f/\partial
x_0,\ldots,\partial f/\partial x_r)$. Then:

\begin{itemize}
\item[{\rm (i)}] For every value of $\,i$ in the range $1\leq i\leq r-1$,
the partial derivatives $\partial f/\partial x_0,\ldots,\partial
f/\partial x_{i}$ divided by their common g.c.d. define a Cremona
transformation of $\pp^i$; in addition, the base ideal of the
inverse map is also a codimension two perfect ideal of linear type
and both ideals are generated in degree $i${\rm ;}

\item[{\rm (ii)}] The linear part of
the graded presentation matrix of $J$ has maximal rank{\rm ;}

\item[{\rm (iii)}] The Hessian of $f$ has the form $h(f)=c\,x_r^ {(r+1)(r-2)}$,
 $c\in k, c\neq 0${\rm ;}

\item[{\rm (iv)}] $f$ is homaloidal.
\end{itemize}
\end{Theorem}
\demo (i) This follows from Proposition~\ref{perfection} and
\cite[Example 2.4]{cremona} (see also
Proposition~\ref{maxrank_is_cremona}, (iii)).

\medskip

(ii)  Again from Proposition~\ref{perfection}  we know that
$J_{r-1}$ is linearly presented, generated by $\partial f/\partial
x_0,\ldots,\partial f/\partial x_{r-1}$ (since for the value
$i=r-1$, the g.c.d. is $1$), hence it yields a chunk of rank $r-1$
of the linear part of the graded presentation matrix of $J$.

In addition, by (\ref{basic_linear_relation}) there is a linear
relation with last coordinate $x_r$ -- hence, nonzero. Clearly then
the rank of the full presentation matrix of $J$ has rank at least
$r$. Since this is the maximal possible value of the rank, the
linear part of the matrix has maximal rank.

\medskip

(iii) By (\ref{g.c.d.s}) one had $\partial f/\partial x_i\in
k[x_{r-i},\ldots, x_r]$, hence $\partial^2f/\partial x_i\partial
x_j=0$ for every $j<r-i$, or equivalently,  $\partial^2f/\partial
x_i\partial x_j=0$ for all pairs $i,j$ such that $i+j\leq r-1$. This
means that the Hessian matrix is anti-lower triangular (i.e., all
zeros below the anti-diagonal). Therefore the determinant is the
product of the entries along its anti-diagonal, namely,
$\partial^2f/\partial x_i\partial x_{r-i}$, for $i=0,\ldots, r$.

We now see that
$$\frac{\partial^2f}{\partial x_i\partial x_{r-i}}=c_ix_r^{r-2},$$
for suitable nonzero constants $c_i\in k$.
 To calculate these derivatives, we induct on $r$.
 One may assume at the outset that $0\leq i\leq r-1$ as otherwise
 one changes the roles of $i$ and $r-i$ not affecting the result except
for the value of the nonzero coefficient. Applying
(\ref{phased_out_r}) in this range we obtain
\begin{equation}\label{phased_out_for_hessian}
\frac{\partial^2f}{\partial x_i\partial x_{r-i}}=\pm\,
x_r\frac{\partial^2\phi^{(r-1)}(x_1,\ldots,x_r)}{\partial
x_{i+1}\partial x_{r-i}}. \end{equation} By the inductive hypothesis
applied  to $f^{(r-1)}$, we deduce that

$$\frac{\partial^2\phi^{(r-1)}(x_1,\ldots,x_r)}{\partial
x_{i+1}\partial x_{r-i}}=c_ix_r^{r-3}.$$
Substituting in
(\ref{phased_out_for_hessian}), we get the stated values.

It now follows that $h(f)=c\,x_r^ {(r+1)(r-2)}$, $c=\Pi c_i\neq 0$.

\medskip

(iv) We will show that the polar map $\phi_f$ is a Cremona map by
drawing upon results from \cite{Simis2} and \cite{cremona2}. Since
the latter is not yet published, we choose to state the method {\it
ab initio}, in the form of a self-contained proposition adapted to
our present purpose.
\medskip

\begin{Proposition}\label{maxrank_is_cremona}
Let $\phi=(F_0:\cdots :F_r)\colon\p^r\dasharrow \p^r$ be a rational
map where $F_0,\ldots ,F_r$ are forms of the same degree generating
an ideal $I\subset k[\xx]$ of codimension $\geq 2$. Set $k[\xx,\yy]$
for the bihomogeneous coordinate ring of $\p^r\times \p^r$ and
consider the bigraded incidence $k$-algebra
$${\cal A}=k[\xx,\yy]/I_1(\yy\cdot\varphi_1),$$
defined by the ideal of entries of the product matrix
$(\yy)\cdot\varphi_1$, where $\varphi$ denotes a graded presentation
matrix of $I$ over $k[\xx]$. Finally, let ${\cal R}={\cal
R}_{k[\xx]}(I)$ stand for the Rees algebra of the ideal $I\subset
k[\xx]$. Then:
\begin{enumerate}
\item[{\rm (i)}] There is a surjective map of bigraded $k$-algebras
$\rho:{\cal A}\surjects {\cal R}${\rm ;}
\item[{\rm (ii)}] If the Jacobian determinant of $F_0,\ldots ,F_r$
is nonzero and if $\ker(\rho)$ is a minimal prime of ${\cal A}$ then
$\phi$ is a Cremona map{\rm $\,$;}
\item[{\rm (iii)}] If the Jacobian determinant of $F_0,\ldots ,F_r$
is nonzero and if $\varphi_1$ has maximal rank  $r$ then $\phi$ is a
Cremona map.
\end{enumerate}
\end{Proposition}
\demo (i) This is a general algebraic fact: there is a structural
surjection ${\cal S}\surjects {\cal R}$ where ${\cal S}$ stands for
the symmetric algebra of $I$. Since ${\cal S}\simeq
k[\xx,\yy]/I_1(\yy\cdot\varphi)$, where $\varphi$ is the full
presentation matrix of $I$, there is a natural surjection ${\cal
A}\surjects {\cal S}$.

(ii) Let $ V\subset \p^r\times \p^r$ stand for the {subscheme}
whose bihomogeneous coordinate ring is ${\cal A}$ and let $
\Gamma\subset \p^r\times \p^r$ stand for the { irreducible}
subvariety whose bihomogeneous coordinate ring is ${\cal R}$, i.e.
$\Gamma$ is the closure of the graph of $\phi$. { Let $V_1,\ldots
, V_r$ denote the distinct irreducible components of $V_{\red}$
where, say, $V_1=\Gamma$. Let $\pi_2:V\to\p^r$ denote the second
projection restricted to $V$ and let $p_2:\Gamma\to\p^r$ stand for
its restriction  to $\Gamma$. Since
 $p_2(\Gamma)=\p^r$, we have $\pi_2^{-1}(p)\neq\emptyset$ for every
$p\in\p^r$. By the nature of $V$, given any point $p\in \p^ r$,
there is a non--negative integer $s(p)$ such that the scheme
theoretic fiber $\pi_2^{-1}(p)$ is of the form $\p^{s(p)}\times
\{p\}$, linearly embedded in $\p^r\times \{p\}$. Since
$\p^{s(p)}\times \{p\}$ is irreducible and reduced, for every
$p\in\p^r$ one has $\pi_2^{-1}(p)=\p^{s(p)}\times \{p\}\subseteq
V_i$ for some $i=i(p)$. Moreover,
$$p_2^{-1}(p)=\Gamma\cap \pi_2^{-1}(p)\subseteq V_1\cap V_{i(p)},$$
as schemes. On the other hand, for every $i\geq 2$ we have
$\dim(V_1\cap V_i)<\dim(V_1)=r$ so that $\dim(p_2(V_1\cap V_i))<r$
for every $i\geq 2$. Let
$$W=\bigcup_{i\geq 2}p_2(V_1\cap V_i)\subsetneq\p^r.$$
Then for every $p\in\p^r\setminus W$ we have
$$p_2^{-1}(p)=\pi_2^{-1}(p)=\p^{s(p)}\times \{p\},$$
as schemes. By the theorem on the dimension of the fibers of a
morphism, there exists an open subset $U\subseteq\p^r$ such that
$\dim(p_2^{-1}(p))=0$ for every $p\in U$. Thus for every $p\in U\cap
(\p^r\setminus W)$ we get $s(p)=0$ and scheme theoretically
$p_2^{-1}(p)$ reduces to a point, yielding the birationality of
$p_2$ and hence of $\phi$.

(iii) One shows that the maximal rank condition implies the
condition on $\ker(\rho)$ in (ii). Namely, note that the incidence
algebra ${\cal A}$ is isomorphic, as a bigraded algebra, to the
symmetric algebra ${\cal S}(E)$ of the $k[\xx]$-module $E={\rm
coker}(\varphi_1)$. The assumption on the rank of $\varphi_1$ then
says that $I\simeq E/(k[\xx]-{\rm torsion})$. By definition, ${\cal
R}\simeq {\cal R}_{k[\xx]}(E)$, where the latter is understood as
${\cal S}(E)/(k[\xx]-{\rm torsion})$ (cf., e.g., \cite{ram}).
Therefore $\ker(\rho)$ is actually the  $k[\xx]$-torsion $\tau({\cal
S}(E))$ of ${\cal S}(E)$. If we show that the torsion is a minimal
prime of ${\cal S}(E)$, we will be done. Now, one has by definition
$$\tau({\cal S}(E))=\ker({\cal S}(E)\rar {\cal S}(E)\otimes_{k[\xx]} k(\xx)),$$
hence $\tau({\cal S}(E))$ is a prime ideal and moreover it is
annihilated by some nonzero $g\in k[\xx]$. It follows that any
graded prime ideal of ${\cal S}(E)$ whose degree zero part vanishes
must contain $\tau({\cal S}(E))$, because it contains $(0)=(g)\cdot
\tau({\cal S}(E))$ and does not contain $(g)$.
Since $\tau({\cal S}(E))$ itself
is one such prime -- because $k[\xx]$ is a domain -- it cannot
properly contain a minimal prime of ${\cal S}(E)$ (necessarily
graded) whose degree zero part is nonzero. Therefore, $\tau({\cal
S}(E))$ has to be a minimal prime itself.
 \qed

\medskip

To conclude the proof of part (3) of the theorem we apply
Proposition~\ref{maxrank_is_cremona}, (iii), and parts (2) and (3)
of the theorem.   \qed

\begin{Remark}\label{jonquieres}\rm (a)
{The proof of part (ii) of the proposition is a more geometric
formulation of the argument in \cite[Theorem 4.1]{Simis2} which
imprecisely claims that {\em every} fiber of $p_2$ is linear. This
is true if $V$ is irreducible, but not otherwise in general: some
special fibers of $\pi_2$ may cut $\Gamma$ along non-linear
varieties.}

(b) We note that the recurrence ideals $J_i$
$(1\leq i\leq r-1)$ are $(x_{r-1},x_r)$--primary ideals, however the
full Jacobian ideal $J$ is not. Geometrically, it obtains that the
singular locus of the sub--Hankel hypersurface is a multiple
structure over the codimension $2$ linear subspace $x_{r-1}=x_r=0$
off the codimension $3$ subspace $x_{r-2}=x_{r-1}=x_r=0$. One can
see that the generators of $J_i$ define a generalized de
Jonqui\`eres transformation as introduced by Pan under the
designation of {\it stellar Cremona maps\/} (cf. \cite{Pan0}). This
would give a different proof that $J_i$ is the base ideal of a
Cremona map, yet the structure of $J_i$ might not follow immediately
from {\it loc. cit.}, let alone that of $J$.

(c) A sub-Hankel determinant $f^{(r)}$ is irreducible. Indeed, we
can readily see that
$$f^{(r)}=-\,\xi(r)x_r^{r-1}x_0+
g(x_1,\ldots,x_r),$$
where $g(x_1,\ldots,x_r)\subset
k[x_1,\ldots,x_r]$ is monic in $x_{r-1}$. Therefore, as a polynomial
in the ring $(k[x_1,\ldots,x_r])[x_0]$ it is primitive and of degree
one in $x_0$, hence is irreducible.

(d) There is enough evidence to conjecture that the ideal $J$ is
also of linear type. For one thing, the computation for various
values of $r$ corroborates the conjecture. Knowing that $J$ is of
linear type would shorten by quite a bit the proof of part (3) of
the theorem and circumvent the need for the full apparatus of
Proposition~\ref{maxrank_is_cremona} and, moreover, it would give
immediately the dominance of the polar map $\phi_f$. Finally, it
would also imply that the inverse map to the polar map is defined by
forms of degree $r$ generating a codimension two perfect ideal --
this has also been computationally checked for various values of
$r$. For a fuller coverage of the syzygy theoretic properties of $J$
see \cite{cremona2}.
\end{Remark}

\begin{Example}\label{subhankel34} \rm It is interesting to consider,
in particular, the first two cases
$r=3,4$. The sub--Hankel surface $V(f^{(3)})$ has degree $3$ and it
has the double line $L$ defined by $x_3=x_2=0$. Hence it is a
rational scroll which is a projection of $S(1,2)\subset \p^ 4$.

Consider the general plane $\pi_{\mathbf{\lambda}}$, with
$\mathbf{\lambda}=(\lambda_2,\lambda_3)\neq \mathbf{0}$, through
$L$, defined by the equation $\lambda_2x_3=\lambda_3x_2$. By
introducing a parameter $t$ and by taking $x_0,x_1,t$ as homogeneous
coordinates in $\pi_\mathbf{\lambda}$, the equation of the
intersection $R_\mathbf{\lambda}$ of $\pi_\mathbf{\lambda}$ with
$V(f^{(3)})$ off $L$ is:

$$
\det \left(
\begin{matrix}
x_0&x_1&\lambda_2\\
x_1&\lambda_2t&\lambda_3\\
\lambda_2&\lambda_3 & 0\\
\end{matrix}
\right )=0
$$

Hence $R_\mathbf{\lambda}$ is a line which varies linearly with
$\mathbf{\lambda}$. In particular, when $\lambda_3=0$,
$R_\mathbf{\lambda}$ coincides with $L$. Thus we see that $L$ is a
line directrix of multiplicity $e=2$,   but
$\mu=1$ (see { Section}  \ref {multidir}). This shows that
$V(f^{(3)})$ is the projection of $S(1,2)$ from a point which lies
in a plane spanned by the $(-1)$--section $E$ and by a ruling $F$,
precisely the one corresponding to $R_{(1,0)}$ -- see
Remark~\ref{genscroll}.

The threefold $V(f^{(4)})$ has degree $4$ and it has the double
plane $\Pi$ defined by $x_3=x_4=0$.

As above, consider the general hyperplane $\pi_{\mathbf{\lambda}}$
through $\Pi$, defined by $\lambda_4x_3=\lambda_4x_3$. By
introducing a parameter $t$ and by taking $x_0,x_1,x_3,t$ as
homogeneous coordinates in $\pi_\mathbf{\lambda}$, the equation of
the intersection $Q_\mathbf{\lambda}$ of $\pi_\mathbf{\lambda}$ with
$V(f^{(4)})$ off $\Pi$ is:

$$
\det \left (
\begin{matrix}
x_0&x_1&x_2&\lambda_3\\
x_1&x_2&\lambda_3t&\lambda_4\\
x_2&\lambda_3t&\lambda_4&0\\
\lambda_3&\lambda_4&0&0\\
\end{matrix}
\right )=0
$$

\noindent One sees that $Q_\mathbf{\lambda}$ is a quadric cone with
vertex $P_\mathbf{\lambda}=[2\lambda_3,-\lambda_4,0,0]$, thus
$P_\mathbf{\lambda}$ sits on $\Pi$ and linearly moves on a line as
$\mathbf{\lambda}$ varies.
\end{Example}


\begin{thebibliography}{Comes}


\bibitem[Al]{Aluffi} P. Aluffi, {\it Singular schemes of
hypersurfaces}, Duke Math. J. 80 (1995), 325--351.

\bibitem[AN]{AN} M. Artin and M. Nagata, {\em Residual intersections in
Cohen--Macaulay rings}, J. Math. Kyoto Univ.  12 (1972), 307--323.

\bibitem[Br]{Bruno} A. Bruno, {\it On homaloidal polynomials},
Michigan Math. J.,  55 (2007), 347--354.

\bibitem[BE]{BE} D. Buchsbaum and D. Eisenbud, {\it What makes a complex
exact?}, J. Algebra 25 (1973),  259--268.


\bibitem[Ca]{Cayley} A. Cayley, {\it On certain developable surfaces},
Quarterly Math. J.  6 (1864), 108--126.


\bibitem[Ci]{Ciliberto} C. Ciliberto, {\it Ipersuperficie algebriche a punti parabolici
e relative hessiane}, Rend. Acc. Naz. Scienze  98 (1979-80),
25--42.

\bibitem[Co]{Conf} F. Conforto, {\it  Le Superficie Razionali}, Zanichelli,
Bologna, 1939.

\bibitem[Cr]{Corti} A. Corti, {\it Factoring birational maps of threefolds after Sarkisov},
J. Algebraic Geom.   4  (1995),  223--254.

\bibitem[DP]{DP} A. Dimca, S. Papadima, {\it  Hypersurfaces complements,
Milnor fibres and higher homotopy groups of arrangements}, Ann. of Math. 158 (2003), 473-507.


\bibitem[Do]{Dolgachev} I. Dolgachev, {\it Polar Cremona transformations}, Mich. Math. J.
 48 (2000), 191--202.


\bibitem[ESB]{ESB} L. Ein, N. Sheperd Barron, {\it Some special Cremona
transformations}, Amer. J. Math.
 111 (1989), 783-800.

\bibitem[Ei]{Eisenbook} D. Eisenbud, {\it Commutative Algebra with a
view toward Algebraic Geometry}, Springer-Verlag, Berlin Heidelberg
New York, 1995.

\bibitem[EH]{EH} D. Eisenbud, J. Harris,
{\it On varieties of minimal degree}, Algebraic Geometry, Bowdoin 1985,
Proc. Symp. in Pure Math.  46 (1987), 3--13.


\bibitem[EKP]{EKP} P. Etingof, D. Kazhdan, A. Polishuck, {\it When is the
Fourier transform of an elementary function elementary?}, Sel. Math. New Ser. 8 (2002), 27--66.

\bibitem [FP] {FP} T. Fassarella, J. V. Pereira, {\it On the degree of polar transformations. An approach through logarithmic foliations}, Sel. Math. New Ser.  13 (2007),  239--252.

\bibitem[FW]{FW} G. Fischer, H. Wu, {\it Developable complex
analytic submanifolds}, Inter. Jour. Math. 6 (1995), 229--272.


\bibitem[Fr1]{Fr1} A. Franchetta, {\it Forme algebriche
sviluppabili e relative hessiane}, Atti Acc. Lincei  10 (1951),
1--4.

\bibitem[Fr2]{Fr2} A. Franchetta, {\it Sulle forme algebriche
di $S_4$ aventi l'hessiana indeterminata}, Rend. Mat.  13
(1954), 1--6.

\bibitem [GR] {GR} A. Garbagnati, F. Repetto, {\it A geometrical approach to
Gordan--Noether's and Franchetta's contributions to a question posed by Hesse},
Collect. Math. 60 (2009), 27--41.


\bibitem[GN]{GN} P. Gordan, M. Noether,
{\it Ueber die algebraischen Formen, deren Hesse'sche Determinante
identisch verschwindet}, Math. Ann.  10 (1876), 547--568.

\bibitem[GH]{GH} P. Griffiths, J. Harris, {\it Algebraic geometry and
local differential geometry}, Ann. Sci. Ecole Norm. Sup. 12
(1979), 355-432.


\bibitem[He1]{Hesse1} O. Hesse, {\it \"Uber die Bedingung, unter
welche eine homogene ganze Function von $n$ unabh\"angigen Variabeln
durch Line\"are Substitutionen von $n$ andern unabh\"angigen
Variabeln auf eine homogene Function sich zur\"uck-f\"uhren l\"asst,
die eine Variable weniger enth\"alt},
 J. reine angew. Math.  42 (1851), 117--124.

\bibitem[He2]{Hesse2} O. Hesse, {\it
Zur Theorie der ganzen homogenen Functionen}, J. reine angew. Math.
 56 (1859), 263--269.


\bibitem[KS]{KS}  T. Kimura, M. Sato, {\it A classification of
irreducible pre--homogeneous vector spaces and their relative
invariants}, Nagoya Math. J. 65 (1977), 133--176.

\bibitem[Kl]{Kleiman} S. Kleiman, {\it Tangency and duality}, Proc.
Vancouver Conf. in Algebraic Geometry (J. Carrell, A.V. Geramita, P.
Russell, eds.), CMS Conf. Proc. 6, Amer. Math. Soc.,  Providence, RI, 1986,
163--226.

\bibitem[KM]{KM} J. Koll\'ar, S. Mori, {\it Birational geometry of algebraic varieties},
Cambridge Tracts in Math.
 134,  Cambridge
University Press, Cambridge, 1990.

\bibitem[Lo]{Lossen} C. Lossen, {\it When does the Hessian determinant vanish identically?
{\rm (}On Gordan and Noether's Proof of Hesse' s Claim{\rm )}}, Bull. Braz. Math. Soc. 35 (2004), 71--82.


\bibitem[Mk]{Mukai} S. Mukai, {\it Simple Lie algebra and Legendre variety},
preprint 1998, unpublished, available at http://www.math.nagoya-u.ac.jp/\~{}mukai.

\bibitem[Pan]{Pan0} I. Pan, {\it Les transformations de Cremona stellaires}, Proc. Amer. Math.
   Soc.  129 (2001),  1257--1262.

\bibitem[PRV]{Pan} I. Pan, F. Ronga, T. Vust, {\it Transformation
birationelle quadratiques de l'espace projectif complexe a trois
dimensions}, Ann. Inst. Fourier  51 (2001), 1153--1187.

\bibitem[Pe]{Perazzo} U. Perazzo, {\it Sulle variet\' a cubiche la
cui hessiana svanisce identicamente}, Giornale di Matematiche (Battaglini)  38  (1900),
337--354.

\bibitem[Pt1]{Permutti1} R. Permutti, {\it Su certe forme a
hessiana indeterminata}, Ricerche di Mat. 6 (1957), 3--10.

\bibitem[Pt2]{Permutti2} R. Permutti, {\it Sul teorema di Hesse
per forme sopra un campo a caratteristica arbitraria}, Le Matematiche
13 (1963), 115--128.

\bibitem[Pt3]{Permutti3} R. Permutti, {\it Su certe classi di forme a
hessiana indeterminata}, Ricerche di Mat.  13 (1964), 97--105.
\bibitem[PW]{control} A. Pinkus, B. Wajnryb, {\it A problem of
approximation using multivariate polynomials}, Russian Math.
Surveys 50 (1995), 319--340.

\bibitem[Ru]{Ru} F. Russo, {\it Tangents and secants of algebraic
varieties}, Notes of a course, $24^ \circ$ Col\'oquio Brasileiro de
Matem\'atica, IMPA, Rio de Janeiro, 2003.


\bibitem[RS1]{cremona} F. Russo, A. Simis, {\it On birational maps and Jacobian matrices},
Comp. Math. 126 (2001), 335--358.

\bibitem[RS2]{cremona2} F. Russo, A. Simis, {\it Cremona maps, ideals of linear
type and linear syzygies}, ongoing.


\bibitem[Se1]{Segre} B. Segre, {\it Bertini forms and hessian
matrices}, J. London Math. Soc. 26 (1951), 164--176.

\bibitem[Se2]{segp} B. Segre, {\it Sullo scioglimento delle
singolarit\a\
delle variet\' a   algebriche}, Ann. Mat. Pura e Appl.  32
(1952), 5--48.

\bibitem[Se3]{Segrediff} B. Segre, {\it Some Properties of
Differentiable Varieties and Transformations}, Erg. Math.,
Springer Verlag, Berlin, 1957.

\bibitem[Se4]{Segre2} B. Segre, {\it Sull' hessiano di taluni
polinomi {\rm (}determinanti, pfaffiani, discriminanti, risultanti,
hessiani{\rm)}} I, II,  Atti Acc. Lincei 37 (1964), 109--117 e 215--221.

\bibitem[Se5]{Segrepr} B. Segre, {\it Prodromi di Geometria Algebrica},
Ed. Cremonese, Roma, 1971.

\bibitem[SeC1]{hessianoSegre} C. Segre, {\it Sulla forma Hessiana}, Rend. Acc. Lincei
5 (1895), 143--148.

\bibitem[SeC2]{Seg} C. Segre, {\it Preliminari di una teoria delle
variet\`a luogi di spazi}, Rend. Circ. Mat. Palermo 30 (1910),
87--121.


\bibitem[Si1]{Simis} A. Simis, {\it Two differential themes in characteristic zero},
in ``Topics in Algebraic and Noncommutative Geometry",
Proceedings in Memory of Ruth Michler (Eds. C. Melles, J.-P.
Brasselet, G. Kennedy, K. Lauter and L. McEwan),  Cont.
Math. 324, Amer. Math. Soc., Providence, RI, 2003,
195--204.

\bibitem[Si2]{Simis2} A. Simis, {\it Cremona transformations and some related algebras},
J. Algebra 280 (2004), 162--179.


\bibitem[SUV]{ram} A. Simis, B. Ulrich  and W. Vasconcelos, {\it Rees algebras of modules}, Proc.
London Math. Soc.  87 (2003), 610--646.

\bibitem[Va]{wolmer} W. Vasconcelos, {\em Arithmetic of Blowup
Algebras}, Londn Math. Soc. Lecture Notes Series 195, Cambridge
University Press, Cambridge, 1994.

\bibitem[Ve] {ves} E. Vesentini, {\it  Sul comportamento effettivo delle
curve polari nei punti multipli}, Ann. Mat. Pura e Appl.  34
(1953), 219--245.


\bibitem[Za1]{Zak1} F. L. Zak, {\em Tangents and secants of algebraic
varieties}, Translations of Mathematical Monographs  127,
Amer. Math. Soc.,  Providence, RI, 1993.

\bibitem[Za2]{ZakHesse} F. L. Zak, {\it Determinants of projective
varieties and their degrees}, in
 ``Algebraic transformation groups and algebraic varieties", Proc. Conference
 ``Interesting algebraic varieties arising in algebraic transformation group theory",
 Enc. Math. Sci. 132, Springer Verlag, Berlin, 2004,  207--238.

\end{thebibliography}
\end{document}